\documentclass[opre,nonblindrev]{informs3}
\DoubleSpacedXI

\usepackage{caption}
\usepackage{subcaption}
\usepackage{float}
\usepackage{xspace}

\usepackage{epstopdf}
\usepackage{cancel}

\usepackage{endnotes}
\let\footnote=\endnote
\let\enotesize=\normalsize
\def\notesname{Endnotes}%
\def\makeenmark{\hbox to1.275em{\theenmark.\enskip\hss}}
\def\enoteformat{\rightskip0pt\leftskip0pt\parindent=1.275em
  \leavevmode\llap{\makeenmark}}

\TheoremsNumberedThrough     % Preferred (Theorem 1, Lemma 1, Theorem 2)

\def\TheoremsNumberedThrough{% 
\theoremstyle{TH}%

\newtheorem{lemma}{Lemma}

\newtheorem{assumption}{Assumption}
\theoremstyle{EX}
\newtheorem{remark}{Remark}
\newtheorem{example}{Example}

\newtheorem{definition}{Definition}

}

\ECRepeatTheorems

\EquationsNumberedThrough    % Default: (1), (2), ...

\usepackage{natbib}
 \bibpunct[, ]{(}{)}{,}{a}{}{,}%
 \def\bibfont{\small}%
\newcommand{\norm}{\|}
\newcommand{\arrbegin}{\begin{eqnarray}}
\newcommand{\arrend}{\end{eqnarray}}

\newcommand{\del}{\partial}
\newcommand{\RR}{\mathbb{R}}

\newcommand{\PP}{\mathop{\mathbb{P}}}

\newcommand{\rarr}{\Rightarrow}

\newcommand{\1}{\mathop{\mathbb{1}}}
\newcommand{\defeq}{\vcentcolon=}

\newcommand{\eps}{\epsilon}

\newcommand{\bb}{\{0, 1\}}

\newcommand{\iid}{\overset{\mathrm{iid}}{\sim}}

\usepackage{mathtools}

\usepackage{color}

\usepackage{algorithm}
\usepackage{algpseudocode}
\usepackage[hypertexnames=false]{hyperref}
\usepackage{dsfont}

\usepackage{bm}
\newcommand{\bmx}{\bm{x}}
\newcommand{\bmu}{\bm{\mu}}
\newcommand{\bmv}{\bm{v}}
\newcommand{\bme}{\bm{e}}
\newcommand{\bmw}{\bm{w}}

\usepackage{fancyhdr}

\begin{document}
%%%%%%%%%%%%%%%%

\DeclarePairedDelimiter\ceil{\lceil}{\rceil}
\DeclarePairedDelimiter\floor{\lfloor}{\rfloor}
\DeclarePairedDelimiter\abs{\lvert}{\rvert}

% \usepackage[toc,page]{appendix}
% \usepackage{lmodern}
% \begin{document}

\RUNAUTHOR{Jalan, Chakrabarti, and Sarkar}

\RUNTITLE{Incentive-Aware Models of Dynamic Financial Networks}

\newcommand{\vc}{\textrm{vec}}
\newcommand{\uvc}{\textrm{uvec}}
\newcommand{\tr}{\textrm{tr}}

%%%%%%%%%%%%%%%%%%%%%
%%%%% Front Matter
%%%%%%%%%%%%%%%%%%%%%

%\TITLE{Modeling and Analysis of Dynamic Incentive-Aware Financial Networks}
%\TITLE{A Dynamic Incentive-Aware Model for Financial Networks}
%\TITLE{From Beliefs to Financial Networks: Incentive-aware Models and Analysis}
\TITLE{Incentive-Aware Models of Financial Networks}

\newcommand{\at}{\makeatletter @\makeatother}

\ARTICLEAUTHORS{%
\AUTHOR{Akhil Jalan}
\AFF{Department of Computer Science, University of Texas at Austin, Austin, TX 78712, \EMAIL{akhiljalan@utexas.edu}} %, \URL{}}
\AUTHOR{Deepayan Chakrabarti}
\AFF{McCombs School of Business, University of Texas at Austin, Austin, TX 78712, \EMAIL{deepay@utexas.edu}}
\AUTHOR{Purnamrita Sarkar}
\AFF{Department of Statistics and Data Science, University of Texas at Austin, Austin, TX 78712, \EMAIL{purna.sarkar@austin.utexas.edu}}
} % end of the block

%%%%%%%%%%%%%%%%%%%%%
%%%%% Abstract
%%%%%%%%%%%%%%%%%%%%%

\ABSTRACT{
  % !TEX root = ./paper-draft.tex
Financial networks help firms manage risk but also enable financial shocks to spread.
Despite their importance, existing models of financial networks have several limitations.
Prior works often consider a static network with a simple structure (e.g., a ring) or a model that assumes conditional independence between edges.
We propose a new model where the network emerges from interactions between heterogeneous utility-maximizing firms.
Edges correspond to contract agreements between pairs of firms, with the contract size being the edge weight.
We show that, almost always, there is a unique ``stable network.''
All edge weights in this stable network depend on all firms' beliefs.
Furthermore, firms can find the stable network via iterative pairwise negotiations.
When beliefs change, the stable network changes.
We show that under realistic settings, a regulator cannot pin down the changed beliefs that caused the network changes.
Also, each firm can use its view of the network to inform its beliefs.
For instance, it can detect outlier firms whose beliefs deviate from their peers. 
But it cannot identify the deviant belief: increased risk-seeking is indistinguishable from increased expected profits.
Seemingly minor news may settle the dilemma, triggering significant changes in the network.

}


%\HISTORY{This paper was first submitted on April 12, 1922 and has been with the authors for 83 years for 65 revisions.}

\maketitle

% {\em Subject classification}: \textcolor{blue}{TODO subject classification from this link.} 

% \\
% {\em Area of Review}: \textcolor{blue}{TODO pick area of review}
% \end{abstract}

%\tableofcontents

% !TEX root = ./paper-draft.tex

\section{Introduction}
% \txtred{Need to cite Ben Golub's papers: ``Targeting Interventions in Networks'' talks about bottom principal component; ``Networks and Economic Fragility'' is a recent survey; ``Supply Network Formation and Fragility'' has a different network model. Also check Matthew Jackson (Stanford), ``A Strategic Model of Social and Economic Networks'' and papers that reference it. ``On perfect pairwise stable networks''. }
The financial crisis of $2008$ showed the need for mitigating systemic risks in the financial system.
There has been much recent work on categorizing such risks~\citep{elliott-golub-jackson-2014,glasserman-young-2015, glasserman-young-2016, birge-2021,jackson-survey-2021}.
While the causes of systemic risk are varied, they often share one feature.
This shared feature is the network of interconnections between firms via which problems at one firm spread to others.
One example is the weighted directed network of debt between firms.
If one firm defaults on its debt, its creditors suffer losses.
Some creditors may be forced into default, triggering a default cascade~\citep{eiseinberg-noe-2001}.
Another example is the implicit network between firms holding similar assets.
Sales by one firm can lead to mark-to-market valuation losses at other firms. 
These can snowball into fire sales~\citep{caballero-simsek-2013, cont16credit, feinstein-2020, feinstein-sojmark-2021}.

%These can snowball into fire sales of assets~\citep{caballero-simsek-2013, cont16credit, feinstein-2020, feinstein-sojmark-2021}.

% todo cite elliot-golub-jackson 2022 here somewhere 

The structure of inter-firm networks plays a vital role in the financial system.
Small changes in network structure can lead to jumps in credit spreads in Over-The-Counter (OTC) markets~\citep{ eisfeldt18otc}. 
Network density, diversification, and inter-firm cross-holdings can affect how robust the networks are to shocks and how such shocks propagate~\citep{ elliott-golub-jackson-2014, acemoglu-2015}.
%\blue{The depth and spread of inter-firm cross-holdings can also mitigate or exacerbate cascading failures \citep{elliott-golub-jackson-2014}.}
The network structure also affects the design of regulatory interventions~\citep{papachristou-kleinberg-2022,amini-2015-control,calafiore-2022,galeotti-2020}. 

Despite its importance, many prior works use simplistic descriptions of the network structure.
For instance, they often assume that the network is fixed and observable.
But only regulators may have access to the entire network.
Furthermore, shocks or regulatory interventions can change the network.
Others assume that the network belongs to a general class.
For instance, \citet{caballero-simsek-2013} assume a ring network between banks.
\citet{amini-2015-control} derive tractable optimal interventions for core-periphery networks.
But financial networks can exhibit complex structure~\citep{peltonen14network, eisfeldt18otc}.
Leverage levels, size heterogeneity, and other factors can affect the network topology~\citep{glasserman-young-2016}.
Hence, there is a need for models to help reason about financial networks.

%Others assume that the network belongs to a general class, such as ring networks or core-periphery networks~\citep{caballero-simsek-2013,amini-2015-control}.

%~\citep{cont_credit_2016, glasserman-young-2016}.}

In this paper, we design a model for a weighted network of contracts between agents, such as firms, countries, or individuals.
The contracts can be arbitrary, and the edge weights denote contract sizes.
In designing the model, we have two main desiderata.
First, the model must account for heterogeneity between firms.
This follows from empirical observations that differences in dealer characteristics lead to different trade risk exposures in OTC markets~\citep{eisfeldt18otc}.
Second, each firm seeks to maximize its utility and selects its contract sizes accordingly.
In effect, each firm tries to optimize its portfolio of contracts~\citep{markowitz-1952}.
The model must reflect this behavior.
From this starting point, we ask the following questions:

\smallskip
\begin{center}
\begin{enumerate}
\item How does a network emerge from interactions between heterogeneous utility-maximizing firms?
\item How does the network respond to regulatory interventions?
\item How can the network structure inform the beliefs that firms hold about each other?
\end{enumerate}
\end{center}

\noindent
Next, we review the relevant literature.

\smallskip\noindent
{\bf Imputing financial networks.}
We often have only partial information about the structure of a financial network.
For example, we may know the total liability of each bank in a network.
From this, we want to reconstruct all the inter-bank liabilities~\citep{squartini-2018-reconstruction}. 
One approach is to pick the network that minimizes KL divergence from a given input matrix~\citep{upper04estimating}.
\cite{mastromatteo12reconstruction} use message-passing algorithms, while \cite{gandy-verart-2017} use a Bayesian approach. 
But such random graph models often do not reflect the sparsity and power-law degree distributions of financial networks~\citep{upper-2011-simulation}.
Furthermore, these models do not account for the utility-maximizing behavior of firms.

\smallskip\noindent
{\bf General-purpose network models.}
The simplest and most well-explored network model is the random graph model~\citep{gilbert_random_1959, erdos_random_1959}.
Here, every pair of nodes is linked independently with probability $p$.
Generalizations of this model allow for different degree distributions and edge directionality~\citep{aiello00random}. %, duijn_p2:_2004}.
Exponential random graph models remove the need for independence, but parameter estimation is costly~\citep{frank_markov_1986, wasserman_logit_1996, hunter_inference_2006, caimo_bayesian_2011}.
%\citep{wasserman_logit_1996, hunter_inference_2006, caimo_bayesian_2011}.
Several models add node-specific latent variables to model the heterogeneity of nodes.
For example, in the Stochastic Blockmodel and its variants, nodes are members of various latent communities.
The community affiliations of two nodes determine their probability of linkage~\citep{holland_stochastic_1983, chakrabarti04rmat, airoldi08mixed, mao_overlapping_2018}.
Instead of latent communities, \cite{hoff02latent} assign a latent location to each node.
Here, the probability of an edge depends on the distance between their locations.

All the latent variable models assume conditional independence of edges given the latent variables.
But in financial networks, contracts between firms are not independent.
Two firms will sign a contract only if the marginal benefit of the new contract is higher than the cost.
This cost/benefit tradeoff depends on all other contracts signed with other firms.
Unlike our model, existing general-purpose models do not account for such utility-maximization behavior.

\smallskip\noindent
{\bf Network games.}
%Games have been widely studied for dynamic routing, path planning, and so on~\citep{wu-2021, wu-2022, peters-2021}.
Here, the payoffs of nodes are dependent on the actions of their neighbors \citep{tardos-2004}.
%\textcolor{blue}{For instance, we can cast the negotiations between firms to set contract sizes as a game.}
%But the underlying network structure is fixed and exogenous, unlike our setting.}
One well-studied class of network games is linear-quadratic games, with linear dynamics and quadratic payoff functions.
Prior work has explored the stability of Nash equilibria~\citep{guo-2021-lq-games} and algorithms to learn the agents' payoff functions~\citep{leng-2020-learning}.
But our model does not yield a linear-quadratic game except in exceptional cases.
Instead, our process involves non-linear rational functions of the beliefs of firms.
Thus, our setting differs from linear-quadratic games.
Recently, network games have been extended to settings where the number of players tends to infinity \citep{carmona-2022-graphon}.
However, we only consider finite networks.

\smallskip\noindent
{\bf Games to form networks.}
Several works study the stability of networks.
In a pairwise-stable network, no pair of agents want to form or sever edges.
This may be achieved via side-payments between agents, which our model also uses~\citep{jackson-wolinsky-2003}.
Pairwise stability has been extended to strong stability for networks~\citep{jackson-nouweland-2005}, and also to weighted networks with edge weights in $[0,1]$~\citep{bich-morhaim-2020,bich-2022-perfect}.
We introduce an analogous notion called Higher-Order Nash stability against any deviating coalition.
However, the weights in our network are not bounded in $[0,1]$ and can be negative.
Furthermore, our edge weights denote contract size, requiring agreement from both parties.
In contrast, prior works typically interpret edge weights as the engagement level in an ongoing interaction.

\citet{golub-sadler-2021} study a network game with endogenous network formation, whose stable points are both pairwise stable and Nash equilibria.
We show similar results for our model.
But they consider unweighted networks and focus on the case of separable games.
In our setting, this corresponds to the case where all firms are uncorrelated.
But in financial networks, correlations are widespread and help firms diversify their contracts.

Several authors study the effect of exogenous inputs on production networks~\citep{herskovic-2018,elliott-golub-leduc-2022}.
\citet{acemoglu-azar-2020} also model endogenous network formation but differ from our approach.
Prices in their model equal the minimum unit cost of production.
For us, prices are determined by pairwise negotiations between firms.
Also, each firm in their model only considers a discrete set of choices among possible suppliers.
In our model, firms can choose both their counterparties and the contract sizes.

%Several authors study the effect of exogenous inputs on production networks~\citep{herskovic-2018,elliott-golub-leduc-2022}.
%\citet{acemoglu-azar-2020} also model endogenous network formation, but they focus on total factor productivity at equilibrium based on input supplier choice.
%In contrast, we consider how firms' beliefs drive their interactions and lead to network formation. \blue{Our endogenous network formation process therefore differs from \cite{acemoglu-azar-2020} in two key respects. First, in our model contract prices $P_{ij}$ are bespoke and determined by pairwise negotations between $(i, j)$, rather than reflecting the marginal unit cost of a good. Second, in our model firms choose their overall level of network participation, as measured by the overall size of their contracts, by optimizing a non-monotonic utility function (Eq.~\eqref{eq:utility}). By contrast, \cite{acemoglu-azar-2020} consider overall network participation level of each firm to be fixed, and the firm must optimize over a discrete set of choices representing different combinations of supplier nodes.}

%endogenous network formation from the strategic interactions among firms, based on their beliefs.

\smallskip\noindent
{\bf Risk-sharing and exchange economies.} 
The pricing of risk is a well-studied problem~\citep{arrow1954existence,buhlmann1980economic,buhlmann1984general,tsanakas2006risk,banerjee-feinstein-2022}.
Most models typically price risk via a global market.
%An important line of work on financial networks studies the problem of pricing risk~\citep{arrow1954existence,buhlmann1980economic,buhlmann1984general,tsanakas2006risk,banerjee-feinstein-2022}.
%The price of risk in these models is typically a global market price.
However, in our model, all contracts are pairwise, and the contract terms and payments between a buyer and seller are bespoke.
There is no global contract or global market price.
Since contracts are pairwise, each firm under our model must consider counterparty risks and the correlations between them.
A firm $i$ may make large payments and accept a {\em negative} reward for a contract with firm $j$ to diversify the risk from contracts with other firms.
Finally, in our model, agents can hedge their risk by betting against one another.
In contrast, B{\"u}hlmann equilibria always result in comonotonic endowments, which firms cannot use as hedges for each other~\citep{banerjee-feinstein-2022, yaari1987dual}.

%Also, in our model, a firm's utility decreases if the contract grows too large, unlike the increasing utility of~\cite{tsanakas2006risk}.

%In our model too, firms balance risk and reward, and coordinate on payments via pairwise negotiations.
%However, we do not require a global market price for any asset, unlike \cite{tsanakas2006risk}. 
%All contracts and payments between a buyer and seller are bespoke.
%Furthermore, firms in our model consider both risks and rewards.
%A firm $i$ may make large payments and accept a {\em negative} reward for a contract with firm $j$ in order to diversify the risk from contracts with other firms.
%In the risk-sharing setting, firm $i$ would not have this option, since the price of $j$'s assets would be fixed for all participants. % \citep{tsanakas2006risk}. 
%Finally, in our model, agents can hedge their risk by betting against one another.
%In contrast, B{\"u}hlmann equilibria always result in comonotonic endowments which cannot be used for hedging~\citep{banerjee-feinstein-2022, tsanakas2006risk}. 

\smallskip\noindent
{\bf Network valuation adjustment.}
Some recent works price the risk due to exposure to the entire financial network~\citep{banerjee-feinstein-2022,feinstein-sojmark-2022}.
The network is usually treated as exogenous and fully known to all firms.
In contrast, we consider endogenous network formation resulting from pairwise interactions between firms.
The network valuation algorithm of \citet{barucca-2020} works with incomplete information, but is not designed for network formation, and it needs firms to share information not required to form their contracts.

%However, the algorithm is not designed for network formation and requires firms to share information not strictly necessary for contract formation.
%While \citet{barucca-2020} propose an algorithm for network valuation in the absence of complete information, the algorithm is not designed for network formation, and it requires firms to share information not strictly necessary for contract formation.

%\blue{
%\smallskip\noindent{\bf Production networks.} The formation of production networks, or input/output networks, is a special case of network formation that has attracted attention in recent years \citep{herskovic-2018,acemoglu-azar-2020,elliott-golub-leduc-2022}. Unlike most works on production networks, we focus on endogenous network formation from the strategic interactions among firms. \citet{acemoglu-azar-2020} also model endogenous network formation, but they focus on total factor productivity at equilibrium based on input supplier choice.
%}

\smallskip\noindent
{\bf Properties of equilibria.}
Another line of work considers the efficiency or social welfare of equilibria~\citep{jackson-survey-2021, elliott-golub-survey-2022}.
\citet{galeotti-2020} show that welfare-maximizing interventions rely mainly on the top or bottom eigenvectors of the network.
\citet{elliott-golub-leduc-2022} show an efficiency-stability tradeoff for their model of supply network formation.
%We do not postulate a notion of social welfare for financial networks.
%We focus on what regulators and firms can learn from limited network information.
Like prior work, we show that stable equilibria exist and are non-dominated.
But our emphasis is on potentially valuable insights for regulators and firms.
For instance, we show a negative result about the ability of regulators to infer the causes of changes to the network structure.
The linkage between firms' utilities and their beliefs, and its effect on stability, is not considered in prior work.

\subsection{Our Contributions}
We develop a new network model of contracts between heterogeneous agents, such as firms, countries, or individuals.
Each agent aims to maximize a mean-variance utility parametrized by its beliefs.
But for two agents to sign a contract, both must agree to the contract size.
For a stable network, all agents must agree to all their contracts.
We show that such constraints are solvable by allowing agents to pay each other.
By choosing prices appropriately, every agent maximizes its utility in a stable network.

\smallskip\noindent
{\bf Characterization of stable networks (Section~\ref{sec:model}):}
We show that unique stable networks exist for almost all choices of agents' beliefs.
These networks are robust against actions by cartels, a condition that we call Higher-Order Nash Stability.
The agents can also converge to the stable network via iterative pairwise negotiations.
The convergence is exponential in the number of iterations.
Hence, the stable network can be found quickly.
Finally, we show how to infer the agents' beliefs by observing network snapshots over time, under certain conditions.

\smallskip\noindent
{\bf The limits of regulation (Section~\ref{sec:regulators}):}
A financial regulator can observe the entire network but not the agents' beliefs.
Suppose firm $i$ changes its beliefs about firm $j$.
Then the contract size between $i$ and $j$ will change.
Indirectly, other contracts will change too.
We show empirically that in realistic settings, the indirect effects can be as significant as the direct effects.
In such cases, the regulator cannot infer the underlying cause of changes in the network.
Similarly, suppose the regulator intervenes with one firm, affecting its beliefs.
The resulting network changes need not be localized to that firm's neighborhood in the network.
Thus, targeted interventions can have strong ripple effects.
Broad-based interventions aimed at increasing stability can also have adverse effects.
For instance, increasing margin requirements on contracts may even increase some contract sizes.
%networks impossible.

\smallskip\noindent
{\bf Outlier detection by firms (Section~\ref{sec:firms}):}
A firm $i$ can observe its contracts with counterparties but not the entire network.
Suppose another firm $j$ (say, a real-estate firm) has beliefs that are very different from its peers.
Then, we prove that under certain conditions, $j$'s contract size with~$i$ is also an outlier compared to other real-estate firms.
So, firm $i$ can use the network to detect outliers and update its beliefs.
But suppose all real-estate firms change their beliefs.
This changes all their contract sizes without creating outliers.
We show that $i$ cannot determine the cause of this change.
For example, firm $i$ would observe the same change whether all real-estate firms had become more risk-seeking or profitable.
However, firm $i$ may want to increase its exposure if they are more profitable but reduce exposure if they are more risk-seeking.
Since the data cannot identify the proper action, firm $i$ remains uncertain.
Exogenous, seemingly insignificant information may persuade firm $i$ one way or another.
Thus, minor news may trigger drastic changes in the network.

\smallskip\noindent
{\bf Notation.}
We use lowercase letters, with or without subscripts, to denote scalars (e.g., $c, \gamma_i$).
Lowercase bold letters denote vectors ($\bmu_i, \bm{w}$), and uppercase letters denote matrices ($W, P, \Sigma_i$).
We use $\bmu_{i;j}$ to refer to the $j^{th}$ component of the vector $\bmu_i$, and $\Sigma_{i;jk}$ for the $(j,k)$ cell of matrix $\Sigma_i$.
We use $\bmv^T$ to denote the transpose of a vector $\bmv$, and $\norm\cdot\norm_p$ to denote the $\ell_p$ norm of a vector or matrix.
We say $A\succeq 0$ if $A$ is positive semidefinite, $A\succ 0$ if it is positive definite, and $A\succeq B$ if $A-B\succeq 0$.
The vectors $\bme_1, \ldots, \bme_n$ denote the standard basis in $\RR^n$, and $I_n$ is the $n\times n$ identity matrix.
If $A \in \RR^{m \times n}, B \in \RR^{p \times q}$ then $A \otimes B \in \RR^{mp \times nq}$ denotes their tensor product: $(A \otimes B)_{ij, k\ell} = A_{ik}B_{j \ell}$. For an appropriate matrix $M$, $\tr(M)$ calculates its trace, $\vc(M)$ vectorizes $M$ by stacking its columns into a single vector, and $\uvc(M)$ vectorizes the upper-triangular off-diagonal entries of $M$. 
For an integer $r \geq 1$, we use $[r]$ to denote the set of integers $[r] \defeq \{1, 2, \dots, r\}$.

% !TEX root = ./paper-draft.tex

\section{The Proposed Model}
\label{sec:model}

We consider a {\em weighted} network $W\in\mathbb{R}^{n\times n}$ between $n$ agents (such as firms, countries, or individuals).
The element $W_{ij}$ represents the size of a contract between agents $i$ and $j$.
We make no assumptions about the content of the contract.
For instance, the contract could be a interest rate swap, a stock swap, or an insurance contract.
Since contracts need mutual agreement, $W_{ij}=W_{ji}$.
We take $W_{ii}$ to represent $i$'s investment in itself.
Note that a negative contract ($W_{ij}=W_{ji}<0$) is a valid contract that reverses the content of a positive contract.
For example, if a positive contract is a derivative trade between two firms, the negative contract swaps the roles of the two firms.

Let $\bmw_i$ denote the $i^{th}$ column of $W$ (i.e., $\bmw_{i;j}=W_{ji}$ for all $j$).
Each agent $i$ would prefer to set its contract sizes $\bmw_i$ to maximize its utility.
But other agents will typically have different preferences.
So, to achieve an agreement about the contract size $W_{ij}$, agents $i$ and $j$ can agree to adjust the terms of their contract.
For example, $i$ may agree to pay $j$ an amount $P_{ji}\cdot W_{ji}$ in cash at the beginning of the contract.
Since payments are zero-sum and $W_{ji} = W_{ij}$, we must have $P_{ji}=-P_{ij}$.
More generally, we allow any adjustments that both agents consider to be equivalent to a cash transfer.
For instance, in a loan contract, the borrower $i$ could pay extra interest over time, if the expected net present value of such interest payments is $P_{ji}\cdot W_{ji}$.
%For instance, $i$ could agree to pay a higher-than-standard interest rate on a loan contract.
%The new rate should change the expected net present value of the contract by the same amount $P_{ji}\cdot W_{ji}$.
We only consider adjustments where both agents agree on the expected net present value.

%Let $\bmw_i$ denote the $i^{th}$ column of $W$ (i.e., $\bmw_{i;j}=W_{ji}$ for all $j$).
%Each agent $i$ would prefer to set its contract sizes $\bmw_i$ to maximize its utility.
%But other agents will typically have different preferences.
%So, to achieve an agreement about the contract size, agents can make \blue{adjustments to the contract}. 
%Such \blue{adjustments} can be cash transfers made at the beginning of the contract.
%They can also refer to \blue{modifications that change the expected net present value of the contract}, such as an increased interest rate on a loan contract. \blue{Note that parties $i, j$ must agree on the expected net present value of such adjustments.}
%We denote by $P_{ji}$ the price per unit contract that $i$ pays to $j$, for a total payment of $P_{ji}\cdot W_{ji}$. 
%If $P_{ji}<0$, then $j$ pays~$i$.
%Since payments are zero-sum and $W_{ji} = W_{ij}$, we must have $P_{ji}=-P_{ij}$.

%Note that $P_{ij}$ includes both immediate payments and the expected net present value of any adjustments to the contract.

% \blue{Price can refer to an additional payment...}

% \blue{\begin{assumption}
% \label{assume:positive-definite}
% Let $X_{ij}$ be a random variable representing the payout obtained by agent $i$ from a unit-sized contract with agent $j$. Each agent $i$ believes that $Cov(X_{ij}, X_{ik}) = \Sigma_{i;jk}$, where $\Sigma_i \succ 0$ is a positive definite matrix. 
% \end{assumption}}

Each contract yields a stochastic payout, and agents have beliefs about these payouts.
We represent agent $i$'s beliefs by a vector ${\bm\mu}_i$ of expected returns and a covariance matrix $\Sigma_i\succ 0$. 
Thus, $\Sigma_i$ represents firm $i$'s perceived risk of trading with other firms, and includes both contract-specific risk and counterparty risk.
Note that we do {\em not} assume that the contracts are zero-sum or that the beliefs are correct, even approximately.
Thus, the overall expected return from all contracts of $i$ is $\bmw_i^T(\bmu_i-P\bme_i)$, and the variance of the overall return is $\bmw_i^T \Sigma_i \bmw_i$.
We assume that each agent has a mean-variance utility~\citep{markowitz-1952}:
\begin{align}
\text{agent $i$'s utility } g_i(W, P) &:= \bmw_i^T (\bmu_i - P \bme_i) - \gamma_i \cdot \bmw_i^T \Sigma_i \bmw_i,
\label{eq:utility}
\end{align}
where $\gamma_i > 0$ is a risk-aversion parameter.
In practice, we expect the set $\{\gamma_i\}_{i\in[n]}$ to be not too heterogeneous~\citep{metrick-1995, kimball-2008, ang-2014, paravisini-2017}. 
%Homogeneity of risk aversion has been observed in game settings~\citep{metrick-1995}, surveys~\citep{kimball-2008}, and investment portfolios~\citep{paravisini-2017} (see also~\cite{ang-2014}).
Note that Eq.~\eqref{eq:utility} ignores costs for contract formation; we will consider these in Section~\ref{sec:regulators:friction}. 
Also, we assume that the contract adjustments $P_{ji}$ do not change the perceived risk.

\begin{example}[Loan contract]
Suppose borrower~$i$ takes a loan of size $W_{ij}$ from lender~$j$.
Then, ${\bm\mu}_{j;i}\cdot W_{ij}$ represents the lender $j$'s expected value for this loan, if it is under the ``standard'' terms.
We allow each pair $(i,j)$ to define their own set of standard terms.
The expected value ${\bm\mu}_{j;i}\cdot W_{ij}$ depends on repayment schedule, the collateral, $j$'s estimate of the probability of default, the recovery rate in case of default, etc.
The borrower's expected value ${\bm\mu}_{i;j}\cdot W_{ij}$ depends on the planned use of this loan. For example, if the loan is meant to purchase equipment, ${\bm\mu}_{i;j}$ is the net present value of expected extra profits due to that equipment.
Hence, ${\bm\mu}_{i;j}$ may not be a function of ${\bm\mu}_{j;i}$.
Now, to achieve agreement on the contract size, the borrower and lender may choose to deviate from the standard terms.
For example, the loan may be callable, or have a higher interest rate than usual.
Then $P_{ji}\cdot W_{ij}$ is the net present value of these deviations for the lender~$j$.
\end{example}

\begin{example}[Interest rate swap contract]
Suppose firm~$i$ makes fixed-rate payments to firm~$j$, and receives floating-rate payments in return.
Then, ${\bm\mu}_{i;j}$ is the expected net present value of these payments for $i$ from a standard unit-sized contract.
This value depends $i$'s forecast of future interest rates and need for floating-rate income, e.g., to match future liabilities.
Hence, it may be quite different from ${\bm\mu}_{j;i}$.
Now, firms~$i$ and~$j$ may adjust the standard terms. 
For example, firm $i$ could agree to pay firm $j$ a higher fixed rate than the standard rate. This corresponds to extra cash transfers over the contract's lifetime with an expected net present value of $P_{ji} \cdot W_{ij}$.
\end{example}
 % to include interest rate caps, extra payments, etc.
% Then, $P_{ji}\cdot W_{ij}$ is the net present value of these adjustments for firm~$j$.

\begin{example}[Insurance Contract]
Suppose firm~$i$ buys fire insurance from insurer~$j$.
Then, ${\bm\mu}_{i;j}$ is the buyer's expected insurance payout minus the insurance premium.
The expected payout depends on the probability of a fire, for which the buyer and insurer may have different estimates.
Also, the insurance contract is negatively correlated with the buyer's other contracts (reflected in $\Sigma_i$).
This is because the buyer gains a payout from the insurer in case of a fire, but incurs losses on other contracts.
Hence, the buyer~$i$ may be willing to accept a contract with negative expected reward, and even pay a higher-than-usual premium $P_{ji}$ per contract.
\end{example}

The model above allows contracts between all pairs of agents.
But some edges may be prohibited due to logistical or legal reasons.
For each agent $i$, let $J_i\subseteq [n]$ denote the ordered set of agents with whom $i$ can form an edge.
So, if $k\notin J_i$ (and hence $i\notin J_k$), we have $W_{ik}=W_{ki}=P_{ik}=P_{ki}=0$.
Similarly, if $i\notin J_i$, then self-loops are prohibited ($W_{ii}=P_{ii}=0$).
We will encode these constraints in the binary matrix $\Psi_i\in\RR^{|J_i|\times n}$ where $\Psi_{i;jk} = 1$ if $k$ is the $j^{th}$ element of $J_i$, and $\Psi_{i;jk}=0$ otherwise.
In other words, $\Psi_i$ is obtained from $I_n$ by deleting the rows corresponding to the prohibited counterparties of $i$.
Thus, for any $\bmv\in\RR^n$, $\Psi_i\bmv$ selects the elements of $\bmv$ corresponding to $J_i$.
If all edges are allowed, we have $\Psi_i=I_n$ for all $i$.
\begin{definition}[Network Setting]
A {\em network setting} $(\bmu_i, \gamma_i, \Sigma_i, \Psi_i)_{i \in [n]}$ captures the beliefs and constraints of $n$ agents.
When there are no constraints (i.e., all edges are allowed), we drop the $\Psi_i=I_n$ terms to simplify the exposition.
Finally, we will use $M\in \RR^{n \times n}$ to denote a matrix whose $i^{th}$ column is $\bmu_i$, and $\Gamma$ to denote a diagonal matrix with $\Gamma_{ii}=\gamma_i$.
\end{definition}

% !TEX root = ./paper-draft.tex

\subsection{Characterizing Stable Points}
\label{sec:model:stable}

In the above model, every agent tries to optimize its own utility (Eq.\eqref{eq:utility}).
We now characterize the conditions under which selfish utility-maximization leads to a stable network.

\begin{definition}[Feasibility]
A tuple $(W, P)$ is feasible if $W=W^T$, $P=-P^T$, and $W$ and $P$ obey the constraints encoded in $(\Psi_i)_{i\in[n]}$.
\end{definition}
\begin{definition}[Stable point]
A feasible $(W, P)$ is stable if each agent achieves its maximum possible utility given prices $P$:
\begin{align*}
g_i(W, P) = \max_{\text{feasible} (W', P) \text{ under } \{\Psi_i\}} g_i(W', P) \quad \forall i\in [n].
\end{align*}
\end{definition}

\begin{figure}[tbp]
\centering
\begin{subfigure}{0.39\textwidth}
  \includegraphics[width=\textwidth]{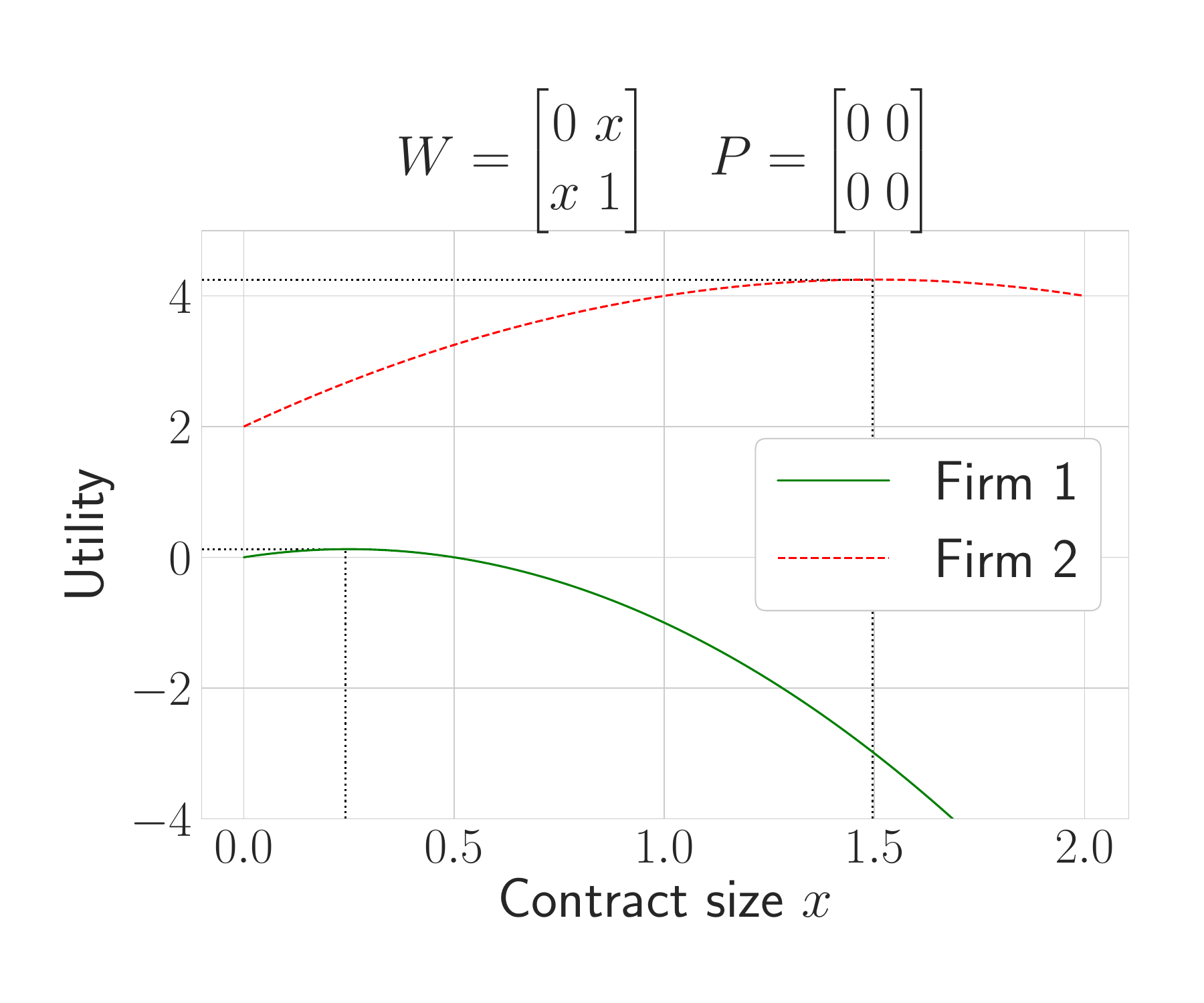}
  \caption{No payments allowed}
  \label{fig:example:noprice}
\end{subfigure}
\begin{subfigure}{0.39\textwidth}
  \includegraphics[width=\textwidth]{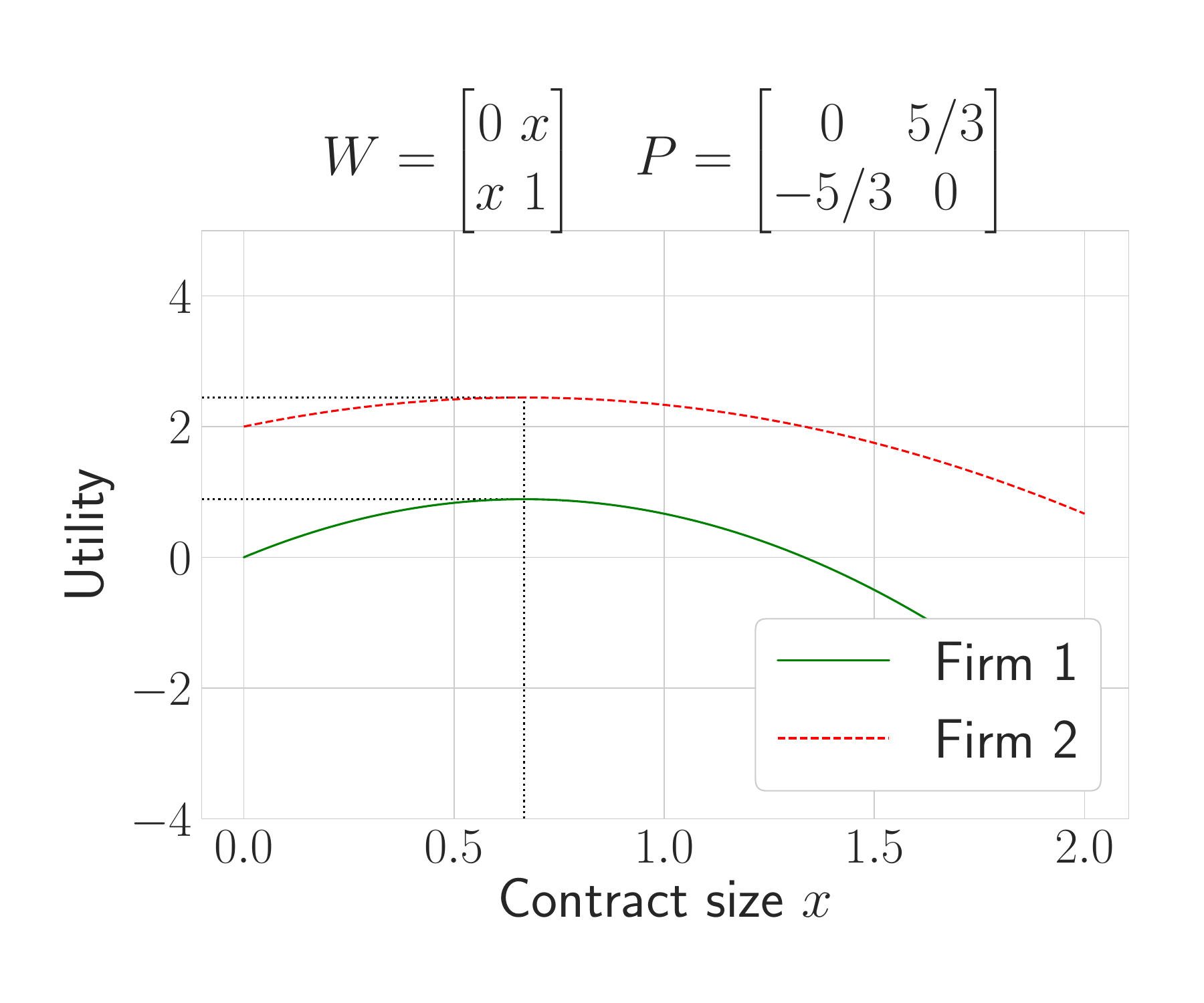}
  \caption{Firm~$2$ pays $5/3$ per contract}
  \label{fig:example:yesprice}
\end{subfigure}
\begin{subfigure}{0.19\textwidth}
\raisebox{1.5em}{
  \includegraphics[width=\textwidth]{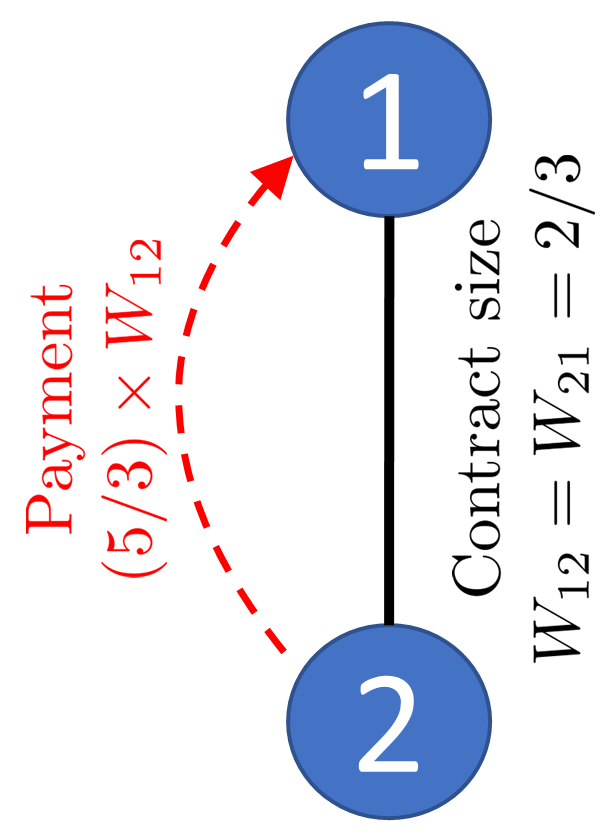}}
  \caption{Network}
  \label{fig:example:network}
\end{subfigure}%
\caption{
{\em Example of a stable point for a borrower (Firm 1) and a lender (Firm 2):} (a) When the borrower cannot pay the lender an additional payment, the firms may be unable to agree to a contract, even if trading improves their utilities.
(b) By allowing for contract-specific payments, both firms can agree on a contract size.
In effect, the borrower (Firm 2) shares its utility with the lender (Firm 1) to achieve agreement.
(c) The stable network is shown.}
\label{fig:example}
\end{figure}

% import matplotlib.pyplot as plt
% import seaborn as sns
% %matplotlib inline

% plt.rcParams['figure.figsize'] = [12, 10]
% sns.set_style('whitegrid')
% plt.rcParams['font.size'] = 36.0
% plt.rcParams['xtick.labelsize'] = 24.0
% plt.rcParams['ytick.labelsize'] = 24.0

% plt.tight_layout()
% plt.savefig('filename.eps', format='eps', dpi=700.0)

\begin{example}
Suppose we only have two firms with the following setting:
\begin{align*}
\text{mean beliefs }M &= \begin{bmatrix}0 & 3\\1 & 4\end{bmatrix} 
& \text{covariance }\Sigma_1=\Sigma_2 &=\begin{bmatrix}1 & 0\\0 & 2\end{bmatrix}
& \text{risk aversion }\gamma_1=\gamma_2&=1.
\end{align*}
So, both firms perceive a benefit from trading ($M_{12}>0, M_{21}>0$).
If trading is disallowed, the optimum $W$ is diagonal with $W_{11}=0$ and $W_{22}=1$ (and $P$ is the zero matrix).
The corresponding utilities are $0$ for firm~$1$ and $2$ for firm~$2$.
Suppose we allow trading but do not allow pricing (Figure~\ref{fig:example:noprice}).
Then, the two firms can each improve their utility by trading, but achieve their optimum utilities at different contract sizes.
Hence, they may be unable to agree to a contract.
In Figure~\ref{fig:example:yesprice}, firm~$2$ pays firm~$1$ a specially chosen price of $5/3$ per unit contract.
At this price, both firms achieve their optimum utilities at the same contract size $W_{12}=W_{21}=2/3$.
Hence, they can agree to a contract.
By paying the price, firm~$2$ shares some of its utility with firm~$1$ to achieve agreement on the contract.
This choice of $W$ and $P$ is a stable point (Figure~\ref{fig:example:network}).
The following results show that this is the {\em only} stable point.$\hfill\Box$
\end{example}\label{example:twofirms}

%\textcolor{blue}{AJ: use the letter $F$ and put into thrm statement. Need to be clear we're not excluding self edges here.}
Define $Q_i = \Psi_i^T (2 \gamma_i \Psi_i \Sigma_i \Psi_i^T)^{-1} \Psi_i$.
When all edges are allowed, $\Psi_i=I_n$ and $Q_i=(2\gamma_i\Sigma_i)^{-1}$.
Let $F = \{(i, j): 1 \leq i < j \leq n, \Psi_i \bm{e}_j \neq \bm{0}\}$ denote the ordered pairs $i<j$ where $P_{ij}$ is allowed to be non-zero.
Note that $|F|\leq n(n-1)/2$. 
%\textcolor{blue}{
%Let $F = \{(i, j): 1 \leq i < j \leq n, \Psi_i \bm{e}_j \neq \bm{0}\}$ denote the ordered set of edges (excluding self-edges) that are permitted according to $\{\Psi_i\}_{i \in [n]}$. We will use $F$ to track nonzero prices. Note that $|F|\leq n(n-1)/2$.}
For any $n\times n$ matrix $X$, let $\uvc(X)_F\in\mathbb{R}^{|F|}$ be a vector whose entries are the ordered set $\{X_{ij}\mid (i,j)\in F\}$.
\begin{theorem}[Existence and Uniqueness of Stable Point]
\label{thm:stable}
Define $n\times n$ matrices $A$, $B_{(i,j)}$, and $C_{(i,j)}$ as follows:
\begin{align*}
A_{ij} &=\bm{e}_i^T Q_j M \bm{e}_j, &
B_{(i,j)} &= \bme_i\bme_j^T Q_i, &
C_{(i,j)} &= (B_{(i,j)} - B_{(j,i)}) - (B_{(i,j)} - B_{(j,i)})^T.
\end{align*}
Let $Z_F$ be the $|F|\times |F|$ matrix whose rows are the ordered sets $\{\uvc(C_{(i,j)})_F \mid (i,j)\in F\}$.
A stable point exists under $\{\Psi_i\}$ if and only if $\uvc(A-A^T)_F$ lies in the column space of $Z_F$.
\end{theorem}

Theorem~\ref{thm:stable} shows that a full-rank $Z_F$ is sufficient for the existence of a stable point (Appendix~\ref{appendix:example-thrm1} shows an example).
%\blue{See Appendix \ref{appendix:example-thrm1} for an example of the matrices $A$, $B_{(i,j)}$, and $C_{(i,j)}$ in a simple setting.}
Furthermore, when the $\Sigma_i$ are random variables, we give a simple sufficient condition that a stable point exists and is unique with probability $1$ (see Appendix \ref{appendix:stable-common}). 
Also, we provide closed-form formulas for the stable point when all agents have the same covariance ($\Sigma_i=\Sigma$ for all $i\in[n]$) (see Appendix \ref{appendix:stability:sharedSigma}).
This occurs when the risk of a contract is primarily counterparty risk (so $\Sigma_{i;jk}$ depends on $j$ and $k$, not $i$) and there is reliable public data on such risks (say, via credit rating agencies).

Next, we consider some properties of the stable point. 
For two feasible tuples $(W_1, P_1)$ and $(W_2, P_2)$, let $(W_2, P_2)$ {\em dominate} $(W_1, P_1)$ if for all $i\in[n], g_i(W_1, P_1)\leq g_i(W_2, P_2)$, with at least one inequality being strict.
\begin{theorem}[Stable points cannot be dominated]
\label{thm:domination}
Suppose a stable point $(W, P)$ exists. Then, there is no feasible $(W', P')$ that dominates $(W, P)$.
\end{theorem}

The stable point obeys a strong form of robustness that we call {\em Higher-Order Nash Stability}.
This strengthens the notions of {\em pairwise stability}~\citep{hellmann-2013} and {\em pairwise Nash}~\citep{calvo-2009-pairwise, golub-sadler-2021} by allowing for agent coalitions, instead of just considering pairs of agents. It is also closely related to the concept of {\em Strong Nash equilibrium}, which strengthens Nash equilibrium by requiring that no subset of agents can deviate at equilibrium without at least one agent being worse off~\citep{mazalov-2019-book}. 
\begin{definition}[Agent Action]
At a given feasible point $(W, P)$, an ``action'' by agent $i$ is the ordered set $(w_{i, j}^\prime, p_{i, j}^\prime)_{j\in J_i}$, where $J_i\subseteq[n]$ is the set of permissible edges for agent $i$.
The action represents a set of proposed changes to $i$'s existing contracts.
Each agent $j\in J_i$ responds as follows:
\begin{enumerate}
    \item If the new $(w_{ij}^\prime, p_{ij}^\prime)$ raises $j$'s utility, then $j$ agrees to the revised contract and price.
    \item Otherwise, $i$ must either keep the existing contract or cancel it ($w_{ij}=p_{ij}=0$).
    We assume that $i$ cancels the contract if and only if this strictly increases $i$'s utility.
\end{enumerate}
We call the shifted $(W^\prime, P^\prime)$ the {\em resulting network}. 
\end{definition}
\begin{definition}[Higher-Order Nash Stability]
A feasible $(W, P)$ is Higher-Order Nash Stable if: 
\begin{enumerate}
    \item {\em Nash equilibrium}: No agent $i$ has an action such that the resulting network $(W^\prime, P^\prime)$ is strictly better for $i$. 
    \item {\em Cartel robustness}: For any proper subset $S \subset [n]$ of agents, there is no feasible point $(W^\prime, P^\prime)$ that differs from $(W, P)$ only for indices $\{i, j\}$ with $i\in S, j\in S$ such that all agents in $S$ have higher utility under $(W^\prime, P^\prime)$ than $(W, P)$. 
\end{enumerate}
\end{definition}
\begin{theorem}[Higher-Order Nash Stability]
\label{thm:pairwise-nash}
Any stable point $(W, P)$ is Higher-Order Nash Stable. 
\end{theorem}

% old utility [0.032781249999999984, 0.03671875000000002, 0.16078125]
% new util [0.052312499999999984, 0.04125000000000001, 0.13671875]

% \blue{AJ: Include an example here of a stable network point, and an example of how players' utilities change after a deviation to a new stable point. Show how some member of the cartel is worse off.}

% !TEX root = ./paper-draft.tex

\subsection{Finding the Stable Point via Pairwise Negotiations}
\label{sec:model:pairwise}

To compute the stable point in Theorem~\ref{thm:stable}, we must know the beliefs of all agents.
But in practice, contracts are set iteratively by negotiations among pairs of agents.
We will now formalize the process of pairwise negotiations and characterize the conditions under which such negotiations can converge to the stable point.

We propose a multi-round pairwise negotiation process.
In round $t+1$, every pair of agents $i$ and $j$ update the price $P_{ij}(t)$ to $P_{ij}(t+1)$ (and hence $P_{ji}(t)$ to $P_{ji}(t+1)$) as follows.
First, they agree to a price $P^\prime_{ij}$ between themselves, {\em assuming optimal contract sizes with all other agents at the current prices $P(t)$.}
In other words, we assume that the other agents will accept the prices in $P(t)$ and the contract sizes preferred by $i$ and $j$. 
Under this condition, $P^\prime_{ij}$ is the price at which $i$'s optimal contract size with $j$ is also $j$'s optimal size with $i$. 
We provide an explicit formula for $P^\prime_{ij}$ in Appendix~\ref{appendix:price-update-rule}.
All pairs of agents calculate these prices {\em simultaneously}.
We create a new price matrix $P^\prime$ from these prices.
Then, we set $P(t+1)= (1-\eta)P(t) + \eta P^\prime$, where $\eta\in (0, 1)$ is a dampening factor chosen to achieve convergence.
Algorithm~\ref{alg:pairwise} shows the details.

%Next, we show how the price $P^\prime_{ij}$ can be computed.
%Slightly abusing notation, we use $P^\prime$ in the next proposition to denote a matrix that is identical to $P(t)$ except for the pair $(i,j)$.

\begin{algorithm}
  \caption{Pairwise Negotiations}
  \label{alg:pairwise}
  \begin{algorithmic}[1]
    \Procedure{Pairwise}{$\eta\in(0,1)$}
      \State $t\gets 0$
      \State $P(0)\gets\text{ any skew-symmetric matrix}$
      \While{$P(t)$ has not converged}
%        \State $W(t) \gets$ matrix with column $i$ being agent $i$'s optimal contracts given $P(t)$ (Prop.~\ref{optimal-contract-given-prices})
%        \State $\forall i,j\in[n], P^\prime_{ij}\gets$ price needed for agreement holding all other prices fixed (Prop.~\ref{price-update-rule})
        \State $\forall i,j\in[n], P^\prime_{ij}\gets$ pairwise-negotiated price for $(i,j)$  (Appendix ~\ref{appendix:price-update-rule})
        \State $P(t+1)\gets (1-\eta) P(t) + \eta P^\prime$
        \State $t\gets t+1$
      \EndWhile
    \EndProcedure
  \end{algorithmic}
\end{algorithm}

\begin{figure}[htbp]
\centering
\includegraphics[width=0.99 \textwidth]{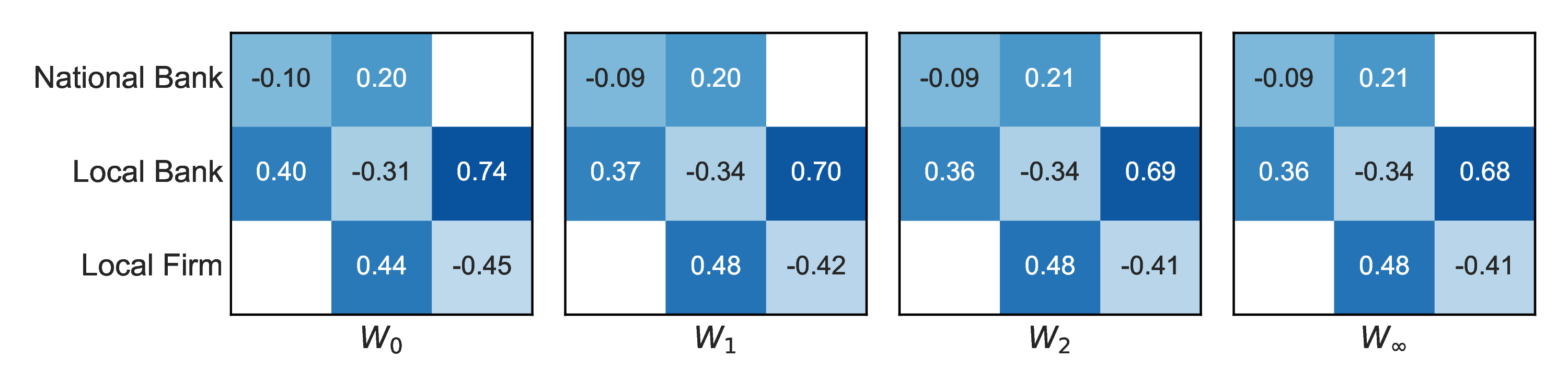} \\
\includegraphics[width=0.99 \textwidth]{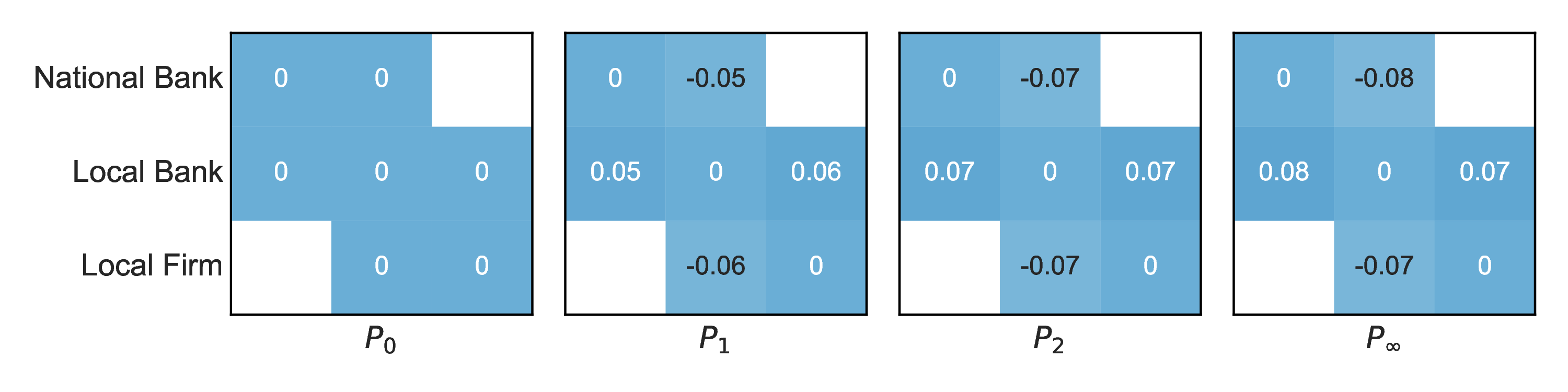}

\caption{
{\em Pairwise negotiations for the setting of Example~\ref{example:pairwise}:}
The contracts matrix $W_t$ and payments matrix $P_t$ after $t = 0, 1, 2$ steps of Algorithm \ref{alg:pairwise} ($\eta=0.5$) converge quickly to the stable point $(W, P) = (W_\infty, P_\infty)$.
Cells corresponding to forbidden edges are empty. 
}
\label{fig:pairwise-example}
\end{figure}

\begin{example}[Pairwise negotations for loan contracts.]
\label{example:pairwise}
Consider a $3$-firm loans network containing a national bank (firm 1), local bank (firm 2), and local firm (firm 3). 
Suppose that the local firm cannot access the national bank, so the edge between firms $1$ and $3$ is prohibited.
The other parameters are:
\begin{align*}
\Sigma_1 = \Sigma_2 = \Sigma_3 &= \begin{bmatrix} 1 & 0.25 & 0.75 \\ 
0.25 & 1 & 0.6 \\ 
0.75 & 0.6 & 1 \end{bmatrix},
& M &= \begin{bmatrix} 0 & 0.9 & 0.9 \\
0.75 & 0 & 0.95 \\
0.5 & 0.8 & 0
\end{bmatrix},
& \gamma_1 = \gamma_2 = \gamma_3 &= 1.
\end{align*}
Figure~\ref{fig:pairwise-example} shows how pairwise negotiations via Algorithm~\ref{alg:pairwise} quickly converge to the stable network.
\end{example}\label{example:pairwise-three-firms}
%We simulate Algorithm \ref{alg:pairwise} for this network, and show that it converges quickly to a stable network in Figure \ref{fig:pairwise-example}. 
% , \bmu_1 = (0, \frac 3 4, \frac 1 2)^T, \bmu_2 = (\frac{9}{10}, 0, \frac{8}{10})^T, \bmu_3 = (\frac{9}{10}, \frac{19}{20}, 0)^T
% todo write M as a matrix 

%{\footnotesize
%\blue{AJ: Here's the same numbers but in matrices.
%\begin{align*}
%W_0 &= \begin{bmatrix}
%  -0.1 & 0.2 & \cancel{0} \\
%  0.4 & -0.31 & 0.74 \\
%  \cancel{0} & 0.44 & -0.45
%\end{bmatrix}
%,
%W_1 &= \begin{bmatrix}
%  -0.09 & 0.2 & \cancel{0} \\
%  0.37 & -0.34 & 0.7 \\
%  \cancel{0} & 0.48 & -0.42
%\end{bmatrix}
%,
%W_2 &= \begin{bmatrix}
%  -0.09 & 0.21 & \cancel{0} \\
%  0.36 & -0.34 & 0.69 \\
%  \cancel{0} & 0.48 & -0.41
%\end{bmatrix},
%W_\infty &= \begin{bmatrix}
%  -0.09 & 0.21 & \cancel{0} \\
%  0.36 & -0.34 & 0.69 \\
%  \cancel{0} & 0.48 & -0.41
%\end{bmatrix} \\
%P_0 &= \begin{bmatrix}
%  0.00 & 0.00 & \cancel{0} \\
%  0.00 & 0.00 & 0.00 \\
%  \cancel{0} & 0.00 & 0.00
%\end{bmatrix}
%,
%P_1 &= \begin{bmatrix}
%  0.00 & -0.05 & \cancel{0} \\
%  0.05 & 0.00 & 0.06 \\
%  \cancel{0} & -0.06 & 0.00
%\end{bmatrix}
%,
%P_2 &= \begin{bmatrix}
%  0.00 & -0.07 & \cancel{0} \\
%  0.07 & 0.00 & 0.07 \\
%  \cancel{0} & 0.00 & 0.00
%\end{bmatrix},
%P_\infty &= \begin{bmatrix}
%  0.00 & -0.08 & \cancel{0} \\
%  0.08 & 0.00 & 0.07 \\
%  \cancel{0} & -0.07 & 0.00
%\end{bmatrix}
%\end{align*}
%}
%}

% \begin{figure}[htbp]
% \centering
% \includegraphics[width=0.99 \textwidth]{figs/examples/P_evol_with_edge_prohibition.eps}
% \caption{
% {\em Pairwise negotiations for the setting of Example~\ref{example:pairwise}:}
% The contracts matrix $W_t$ after $t = 0, 1, 2$ steps of Algorithm \ref{alg:pairwise} ($\eta=0.5$) converges quickly to the stable point $W=W_\infty$.
% Forbidden edges are crossed out.
% }
% \label{fig:pairwise-example}
% \end{figure}

Now, we will show that Algorithm~\ref{alg:pairwise} converges.
First, we define {\em global asymptotic stability} (following \citet{callier-desoer}).

\begin{definition}[Global Asymptotic Stability]
The pairwise negotiation process is globally asymptotically stable for a given network setting and dampening factor $\eta$ if, for any initial price matrix $P(0)$, there exists a matrix $P^\star$ such that the sequence of price matrices $P(t)$ converges to $P^\star$ in Frobenius norm:
$\lim\limits_{t\to\infty} \|P(t) - P^\star\|_F = 0$.
\end{definition}

When pairwise negotiations are globally asymptotically stable, the limiting matrix $P^\star$ must be skew-symmetric since each $P(t)$ is skew-symmetric.
Also, since prices are updated whenever two agents disagree on the size of the contract between them, all agents agree on their contract sizes at $P^\star$.
Hence, $P^\star$ must be a stable point for the given network setting.

Now, we show that for a range of $\eta$, pairwise negotiations are globally asymptotically stable (Appendix~\ref{appendix:example-asymp} presents an example).

\begin{theorem}[Convergence Conditions and Rate]
\label{thm:asymp}
Let $Q_i$ be defined as in Theorem \ref{thm:stable}. 
Define the following $n^2\times n^2$ matrices:
\begin{align*}
K &\defeq \sum\limits_{r=1}^{n} \bm{e}_r \bm{e}_r^T \otimes Q_r + Q_r \otimes \bm{e}_r \bm{e}_r^T & \\
L_{(i - 1)  n + j, (i - 1)  n + j} &= Q_{i;j,j} + Q_{j; i,i} \quad\forall i,j\in[n]  & \text{($L$ is diagonal).}
%$ for all $i,j\in[n]$. 
\end{align*}
%Let $K \defeq \sum\limits_{r=1}^{n} \bm{e}_r \bm{e}_r^T \otimes Q_r + Q_r \otimes \bm{e}_r \bm{e}_r^T$.
%Let $L, R\in\RR^{n^2 \times n^2}$ be diagonal matrices  $L_{(i - 1)  n + j, (i - 1)  n + j} = Q_{i;j,j} + Q_{j; i,i}$, and $R_{(i - 1)  n + j, (i - 1)  n + j} = \mathds{1}_{\{i, j\} \text{ permitted}}$ for all $i,j\in[n]$. 
Let $L^\dagger$ denote the pseudoinverse of $L$, and $(L^\dagger K)\mid_R$ denote the principal submatrix of $L^\dagger K$ containing the rows/columns $(i-1)n+j$ such that the edge $(i,j)$ is not prohibited.
Let $\lambda_{\textrm{max}}, \lambda_{\textrm{min}}$ be the largest and smallest eigenvalues of the matrix $(L^\dagger K)\mid_R$ respectively.
Let $\eta^* = \frac{2}{\lambda_{\textrm{max}}}$.
%Let $\eta^* = \min\{1, \frac{2}{\lambda_{\textrm{max}}}\}$.
%Let $\underline{\lambda}$ and $\overline{\lambda}$ denote the minimum and maximum values in the set $\Lambda\defeq \{Re(\lambda): \lambda\in\sigma(L^{-1} (I_n \otimes C + C \otimes I_n)\}$, where $\sigma(.)$ denotes the set of eigenvalues of a matrix.
Then, we have:
\begin{enumerate}
\item For all $\eta \in (0, \eta^*)$, pairwise negotiations with $\eta$ are globally asymptotically stable. 
\item For such an $\eta$, the convergence is exponential in the number of rounds $t$:
\begin{align*}
\|P(t)-P^\star\|_F &\leq \frac{\alpha^t}{1-\alpha} \cdot \|P(1)-P(0)\|_F, &\text{where }\alpha = \max\{ \abs{1-\eta\lambda_{\textrm{min}}}, \abs{1-\eta\lambda_{\textrm{max}}} \}.
\end{align*}
\end{enumerate}
Here, $P^\star$ is the stable point to which the negotiation converges.
\end{theorem}

\begin{remark}
For clarity of exposition we restrict $\eta \in (0, 1)$ in Algorithm \ref{alg:pairwise}. However, Theorem \ref{thm:asymp} shows that we only need $\eta < \eta^*$ for convergence to the stable point.
\end{remark}

%\blue{See Appendix \ref{appendix:example-asymp} for an example of convergence conditions and rates in the setting of example \ref{example:thrm1}.}

% \blue{AJ: We might want to reword Theorem $4$ to say that $\eta^* \leq 1$ is required. This is implicit in the algorithm.}
% \blue{
% \begin{example}
% In the setting of Example \ref{example:twofirms}, we have 
% $K = L = \begin{bmatrix} 1 & 0 & 0 & 0 \\
% 0 & 3/4 & 0 & 0 \\
% 1 & 0 & 3/4 & 0 \\
% 1 & 0 & 0 & 1/2
% \end{bmatrix}
% $  
% \end{example}
% }
\subsection{Pairwise Negotiations under Random Covariances}
\label{sec:model:randomcov}

So far we have made no assumptions about agents' beliefs. 
In this section, we analyze the convergence of pairwise negotiations for ``data-driven'' agents.
Specifically, each agent $i$ now {\em estimates} its covariance matrix.
For this section only, we will call the covariance matrix $\hat{\Sigma}_i$ instead of $\Sigma_i$ to emphasize that it is an estimated quantity. 

Suppose each agent $i$ observes $m$ independent data samples.
Each sample is a vector of the returns of unit contracts with all $n$ agents.
The samples for agent $i$ are collected in a matrix $X_i\in\RR^{n\times m}$, with one column per sample.
The sample covariance of this data is $\hat{\Sigma}_i$.

%Specifically, $\hat{\Sigma}_i$ is estimated by agent $i$ from data on the returns of unit contracts with all $n$ agents.
%Suppose the return distribution $\mathcal{D}$ is unknown, but each agent observes $m$ independent samples from $\mathcal{D}$.
%The samples for agent $i$ are collected in a matrix $X_i\in\RR^{n\times m}$, with one column per sample, from which it calculates $\hat{\Sigma}_i$. Under a wide range of conditions, $\hat{\Sigma}_i$ converges in probability to the covariance $\Sigma$ of $\mathcal{D}$ \citep{vershynin-book}. 
%
%According to Theorem~\ref{thm:asymp}, this yields a maximum dampening rate $\eta^\star$ that depends on $\Sigma$.

We assume that all agents observe samples from the same return distribution, which has covariance $\Sigma$.
Under a wide range of conditions, $\|\hat{\Sigma}_i - \Sigma\|\to 0$ in probability~\citep{vershynin-book}. 
Hence, at convergence, the maximum allowed dampening rate $\eta^\star$ in Theorem~\ref{thm:asymp} would be a function of $\Sigma$.
But for finite sample sizes, each agent's $\hat{\Sigma}_i$ can be different.
Hence, the maximum dampening $\hat{\eta}^\star$ may be less than $\eta^\star$.
The smaller the $\hat{\eta}^\star$, the worse the rate of convergence of pairwise negotiations.
However, even with a few samples, $\hat{\eta}^\star$ is close to $\eta^\star$, as the next theorem shows.

\begin{theorem}[Small Sample Sizes are Sufficient for Fast Convergence]
\label{thm:asymp:random}
Suppose that $\norm \Sigma \norm, \norm \Sigma^{-1} \norm, \norm \Gamma \norm,$ and $\norm \Gamma^{-1} \norm$ are $O(1)$ with respect to $n$ and all edges are allowed. 
Also, suppose that each sample column of $X_i$ is drawn independently from a $\mathcal{N}(\bm{0}, \Sigma)$ distribution, and let $\hat{\bm{\mu}} = \frac{1}{m} \sum_i X_i$ and ${\hat{\Sigma}_i \defeq \frac{1}{m-1} \sum_i (X_i - \hat{\bm{\mu}}) (X_i - \hat{\bm{\mu}})^T}$.
Let $\hat{\eta}^\star$ be the maximum dampening factor using $(\hat{\Sigma}_i)_{i\in[n]}$ as defined in Theorem~\ref{thm:asymp}.
Let $\eta^\star$ be the dampening factor if $\hat{\Sigma}_i$ were replaced by $\Sigma$ for all $i$. 
If $m = \ceil*{n \log n}$, then for large enough $n$, ${\hat{\eta}^\star \geq (1 - o(1)) \eta^\star}$ with probability at least $1 - \exp(-\Omega(n))$. 
\end{theorem}
%\blue{AJ: Fixed the estimator $\hat\Sigma_i$ to be unbiased. We don't need the biased form at all, it's simply the form that Vershynin used. However his proof for the bound assumed centered data already, so it's equivalent to what we're doing here.}

Theorem~\ref{thm:asymp:random} shows that data-driven agents using a broad range of dampening factors are still likely to find the stable point via pairwise negotiations. 
Furthermore, the amount of data they need is comparable to the number of agents (up to a logarithmic factor). We note that if firms use datasets of fixed sizes $m_1, \dots, m_n$, then the conclusion of Theorem~\ref{thm:asymp:random} still holds, as long as $\min_{i} m_i \geq \ceil*{n \log n}$. For example, firms might use different look-back periods for covariance estimation.

% !TEX root = ./paper-draft.tex
\subsection{Inferring Beliefs from the Network Structure}
\label{sec:model:inference}
Suppose we are given a network that lies at a unique stable point as defined in Theorem~\ref{thm:stable}.
How can we infer the beliefs of the agents?

\paragraph{Non-identifiability of beliefs.}
Suppose we are given a network $W$ that is generated using a single covariance $\Sigma_i=\Sigma\succ 0$.
We want to infer the agents' beliefs $(M, \Gamma, \Sigma)$.
By Corollary~\ref{cor:stability:sharedSigma}, 
\begin{align*}
\frac 1 2 \vc(M+M^T) &= (\Gamma\otimes\Sigma + \Sigma\otimes \Gamma)\vc(W).
\end{align*}
Clearly, the agents' beliefs can only be specified up to an appropriate scaling of $M$, $\Gamma$, and $\Sigma$.
But even if we specify a scale (e.g., $\tr[\Gamma]=\tr[\Sigma]=1$), for any valid choice of $\Gamma$ and $\Sigma$ we can find a corresponding $M$.
Thus, even in the simple setting of identical covariance and fixed scale, the network $W$ cannot be used to select a unique combination of the parameters $(M, \Gamma, \Sigma)$.
By a similar argument, we cannot identify the underlying beliefs even if we observe multiple networks generated using the same $\Sigma$ and $\Gamma$ (but different $M$).
Thus, we need further assumptions in order to infer beliefs.

\begin{assumption}
\label{assume:SDP}
%Consider a sequence of networks $W(t)$ generated over successive timesteps $t\in[T]$. %$t=1, \ldots, T$.
Consider a sequence of networks $W(t)$ over timesteps $t\in[T]$. %$t=1, \ldots, T$.
We assume that (a)~$\Gamma(t)=I$ and $\Sigma_i(t)=\Sigma$ for all $t\in[T]$, (b)~for all $i,j\in[n]$, $M_{ij}(t)$ varies independently according to a Brownian motion with the same parameters for all $(i,j)$, and (c)~$\tr\Sigma = 1$.
\end{assumption}

%We assume $\Gamma=I$ due to the homogeneity of risk aversion noted in Section~\ref{sec:model}.
%The choice of $\Sigma_i(t)=\Sigma$ is motivated by observations that errors in mean estimation are far more significant than covariance estimation errors~\citep{chopra-2013}. 
%So, accounting for variations in $\Sigma$ may be less important than variations in $M$ (but see Remark~\ref{rem:timevaryingSigma} below).
%The Brownian motion assumption is common in the literature on pricing models~\citep{geman-2001,bianchi-2013}. 
%The choice of $\tr\Sigma=I$ fixes the scale, as discussed above.

The first assumption is motivated by the observations in portfolio theory that errors in mean estimation are far more significant than covariance estimation errors~\citep{chopra-2013}. 
So, accounting for variations in $\Sigma$ may be less important than variations in $M$ (but see Remark~\ref{rem:timevaryingSigma} below).
The homogeneity of risk aversion was noted in Section~\ref{sec:model}, and this justifies setting $\Gamma=I$.
The second assumption is common in the literature on pricing models~\citep{geman-2001,bianchi-2013}. 
The third assumption fixes the scale, as discussed above.
%Under these assumptions, we can infer $\Sigma$ (and hence $M(t)$) via a computationally tractable Semidefinite Program (SDP), as shown next.

\begin{proposition}\label{prop:sdp-sigma-recovery}
Finding the maximum likelihood estimator of $\Sigma$ under Assumption~\ref{assume:SDP} is equivalent to the following Semidefinite Program (SDP):
\begin{align*}
\min\limits_{\Sigma} \sum\limits_{t = 1}^{T - 1}
\norm \Sigma (W(t + 1) - W(t)) + 
(W(t + 1) - W(t)) \Sigma \norm_F^2  
\text{$\quad \text{ s. t. }
\Sigma \succeq 0,
\tr(\Sigma) = 1.$}
\end{align*}
\end{proposition}

% \blue{As the number of observations $T$ increases, the maximum likelihood estimator converges to the true $\Sigma$ - see Figure \ref{fig:sigma-inference}}. 

\begin{remark}[Generalization to time-varying $\Sigma$]
\label{rem:timevaryingSigma}
Instead of a constant covariance $\Sigma$, the time range may be split into intervals, with covariance $\Sigma_{(j)}$ in interval $j$.
Then, we can add a regularizer $\nu \cdot \sum_j \|\Sigma_{(j+1)}-\Sigma_{(j)}\|$ for some $\nu>0$ to the objective of the SDP to penalize differences between successive covariances.
This allows the covariance to evolve while keeping the objective convex.
The time intervals can be tuned based on heuristics or prior information.
\end{remark}

\section{Insights for Regulators}
\label{sec:regulators}

A financial regulator can observe the network but does not know the firms' beliefs.
The regulator may ask: what changes in beliefs caused recently observed changes in the network?
What are the side effects of different regulatory interventions?
To answer these questions, we need to know how changes in firms' beliefs or utility functions affect the network.
That is the subject of this section.

% \subsection{Newer budget constraints section}

\subsection{Effect of Friction in Contract Formation}
\label{sec:regulators:friction}

Our model imposes no costs for contract formation.
This is reasonable for large firms where the fixed costs associated with contract negotiations may be small relative to the contract sizes.
However, in an overheating market, a regulator may impose frictions by penalizing large contracts, for example by increasing margin requirements.

We model contract costs via an adding a penalty term $F_i(\bmw_i)$ to the utility of agent $i$ in Eq.~\eqref{eq:utility}:
\begin{align}
\label{eq:utilityFriction-general}
\text{agent $i$'s utility } g_i(W, P) &:= \bmw_i^T (\bmu_i - P \bme_i) - \gamma_i \cdot \bmw_i^T \Sigma_i \bmw_i - F_i(\bmw_i).
\end{align}

\begin{theorem}\label{thm:friction-equilibrium-general}
Consider a network setting where $\Sigma_i=\Sigma$ and all edges are allowed.
Suppose that for each firm $i\in [n]$, the function $F_i: \RR^n \to \RR$ is twice differentiable, and there exist strictly increasing functions $f_{ji}: \RR \to \RR$ such that  for all $\bm{x} \in \RR^n$, $\nabla F_i(\bm{x}) = [f_{1i}(x_1), \dots, f_{ni}(x_n)]^T$.
Then, there exists a unique stable point.
%Then {for a network setting $(\bmu_i, \Sigma, \gamma_i, I)_{i \in [n]}$,} there exists a unique stable point.
\end{theorem}

%\blue{Our theorem as it is assumes shared $\Sigma$ and no prohibited edges. I don't feel like generalizing it, so we need to clarify the assumption.}

\begin{example}
By imposing frictions, the regulator may {\em increase} the sizes of certain contracts. 
For example, let $F_i(\bm{w}_i) = \epsilon\cdot w_{i;i}^2 + \lambda\cdot \sum_{j\neq i}w_{i;j}^2$ for some $\lambda>\epsilon> 0$.
Thus, the cost of inter-firm trades scales with the square of the contract size (we assume $\epsilon\approx 0$).
Consider a network setting with $3$ firms, with $\gamma_i = 1$, $\Sigma_i = \Sigma = \begin{bsmallmatrix} 0.1 & 0.1 & 0.1 \\ 0.1 & 1 & 0.5\\0.1 & 0.5 & 1 \end{bsmallmatrix}$, and
$M = \begin{bsmallmatrix} 0 & 1000 & 111.233\\ 1000 & 1 & 0.1 \\ 1000 & 0.1 & 1 \end{bsmallmatrix}$.
Then, $W_{23}=W_{32}\approx 0$ without frictions (when $F_i(\bm{w}_i)=0$) but $|W_{23}|>0$ for $\lambda>0$.
\end{example}

%We set $\epsilon>0$ to satisfy the conditions of Theorem~\ref{thm:friction-equilibrium-general}.

%\begin{example}
%By imposing frictions, the regulator may {\em increase} the sizes of certain contracts. 
%For example, let $F_i(\bm{w}_i) = \lambda \norm \bm{w}_i \norm_2^2$ for some $\lambda\geq 0$.
%Consider a network setting with $2$ firms, with $\gamma_1 = \gamma_2 = 1$, $M = \begin{bmatrix} 2 & 1 \\ 1 & 1 \end{bmatrix}$, and $\Sigma_1 = \Sigma_2 = \begin{bmatrix} 3 & 2\\ 2 & 3\end{bmatrix}$.
%Then, the size of the contract between the two firms is
%$W_{12}^{\lambda} = \frac{\lambda}{(2\lambda^2 + 12\lambda + 10)}.$
%Hence for any $\lambda > 0$, $W_{12}^{\lambda}$ is strictly positive, whereas for $\lambda = 0$, $W_{12} = 0$.
%\end{example}

\subsection{Effect of Changes in Firms' Beliefs}

%\blue{AJ: Rewrote (and renamed) this section. Tried to keep wording as similar as possible.}. 

%Regulatory actions can change the risk and expected return perceptions of firms, who then update their \blue{ means and } covariances accordingly. 
%If perceived reward increases then contract sizes increase, as expected. \blue{Similarly}, if the risk scales upwards, we find smaller contract sizes. But this does not generalize to non-uniform increases in $\Sigma$. 

Regulatory actions can change the risk and expected return perceptions of firms.
The next theorem shows the effect of such belief changes on the stable point.

%\blue{AJ: Stylistically I like the format of thrm 8.1, 8.2, ... but can change if needed.}
\begin{theorem}\label{thm:effect-changed-sigma-mu}
Suppose $\Sigma_i=\Sigma$ for all firms, and let $M$ be the matrix of expected returns. %Then, we have the following:
\begin{enumerate}
  \item {\bf Change in beliefs about expected returns:}
Let $\Sigma$ have the eigendecomposition $\Sigma = V \Lambda V^T$.
Then for $i, j, k, \ell \in [n]$, 
\begin{align}
\label{eq:delWdelM_general}
\frac{\del W_{ij}}{\del M_{k\ell}} = \frac{1}{2\sqrt{\gamma_i\gamma_j\gamma_k\gamma_\ell}} \sum\limits_{s, t \in [n]} \frac{V_{is} V_{ks} V_{jt} V_{\ell t} + V_{is}V_{\ell s}V_{jt}V_{kt}}{\lambda_s + \lambda_{t}}.
\end{align}
  In particular, $W_{ij}$ is monotonically increasing with respect to $M_{ij}$.
  \item {\bf Risk scaling:} If the covariance $\Sigma$ changes to $c\Sigma$ ($c>0$), then $W$ changes to $(1/c)W$.
  \item {\bf Increase in perceived risk:} Suppose $\gamma_i = \gamma$ for all $i$, and the covariance $\Sigma$ increases to $\Sigma^\prime\succ \Sigma$.
Let $W$ and $W^\prime$ be the stable points under $\Sigma$ and $\Sigma^\prime$ respectively.
Then, $\tr(M^T (W^\prime - W)) < 0.$
\end{enumerate}
\end{theorem}

This shows that, in general, an increase in risk leads to a decrease in the weighted average of the contract sizes.
The weights are given by the expected return beliefs of the firms.
However, individual contracts between firms can increase, as can the norm $\|W\|_F$.
This is because increases in the covariance $\Sigma$ may also increase correlations, which can offer better hedging opportunities.
By hedging some risks, larger contract sizes can be supported.

% \subsection{Effect of Changes in Perceived Expected Returns}

%\subsection{Network Effect of Changes in Perceived Expected Returns}

%\blue{AJ: renamed the section to ``network effect.''}
% \begin{theorem}[$W$ is monotonic with $M$]
% \label{thm:w_monotonic_mu}
% Suppose $\Sigma_i=\Sigma$ for all firms.
% Then, for any $i,j \in [n]$, the value of $W_{ij}$ is monotonic with respect to $M_{ij}$ (the $i^{th}$ component of $\bmu_j$). 
% \end{theorem}

%A change in $M_{ij}$ may affect not only $W_{ij}$ ( but also all other contracts.
Theorem~\ref{thm:effect-changed-sigma-mu} also shows that a change in the perceived expected return $M_{k\ell}$ affects all contracts $W_{ij}$.
Can we trace the changes in $W$ back to the underlying changes in $M$?
For instance, consider the following problem.
\begin{definition}[Source Detection Problem]
Suppose that a financial regulator observes two networks $W$ and $W^\prime$, with the only difference being a small change in a single entry of $M$ (say, $M_{ij}$).
Can the regulator identify the pair $(i, j)$?
\end{definition}

One approach is to try to infer all beliefs of all firms, and then identify the changed belief.
But, as discussed in Section~\ref{sec:model:inference}, the beliefs are only identifiable under extra assumptions and more data.
An alternative approach for the source detection problem is to find the entry $(i,j)$ with the largest change $|W_{ij}-W^\prime_{ij}|$.
The intuition is that a change in $M_{ij}$ has a direct effect on $W_{ij}$ and (hopefully weaker) indirect effects on other contracts.
Thus, the source detection problem is closely tied to the following:
\begin{definition}[Targeted Intervention Problem]
Can a regulator induce a small change in a single entry of $M$ (say, $M_{ij}$) such that the change in $W_{ij}$ is significantly larger than changes in other entries of $W$?
\end{definition}

When all eigenvalues of $\Sigma$ are equal (that is, $\Sigma\propto I_n$), a change in $M_{k\ell}$ only affects $W_{k\ell} (=W_{\ell k})$, as can be seen from Corollary~\ref{cor:stability:sharedSigma}.
But when the eigenvalues are skewed, the terms in Eq.~\eqref{eq:delWdelM_general} corresponding to the smallest eigenvalues have greater weight.
In such circumstances, the indirect effect of a change in $M_{k\ell}$ on other $W_{ij}$ can be significant.
The following empirical results show that this is indeed the case.

\begin{figure}[tbp]
\centering
\begin{subfigure}{0.48\textwidth}
  \includegraphics[width=\textwidth]{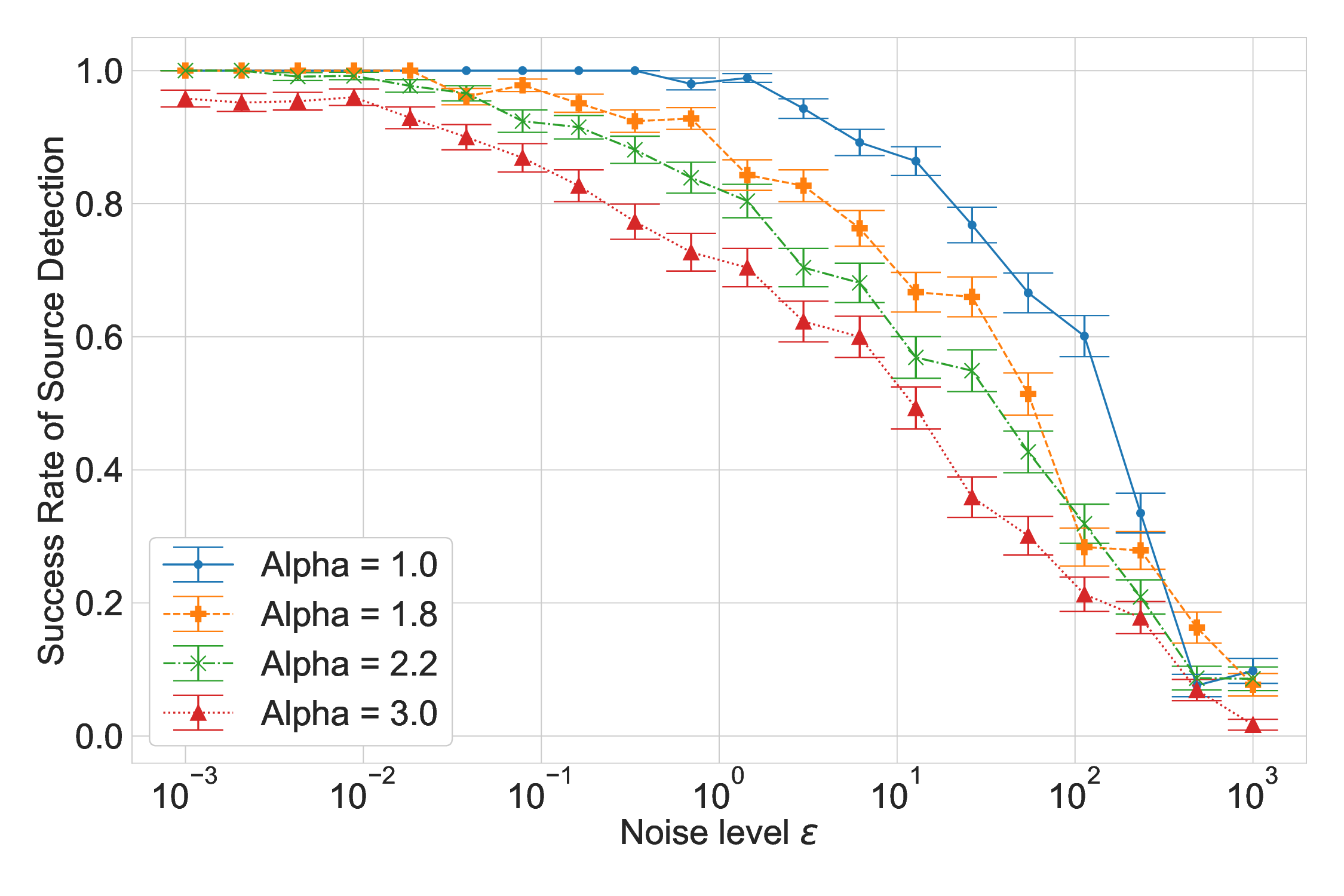}
  \caption{Predict most shifted contract as source}
  \label{fig:equicorr-left}%
\end{subfigure}
\begin{subfigure}{0.48\textwidth}
  \includegraphics[width=\textwidth]{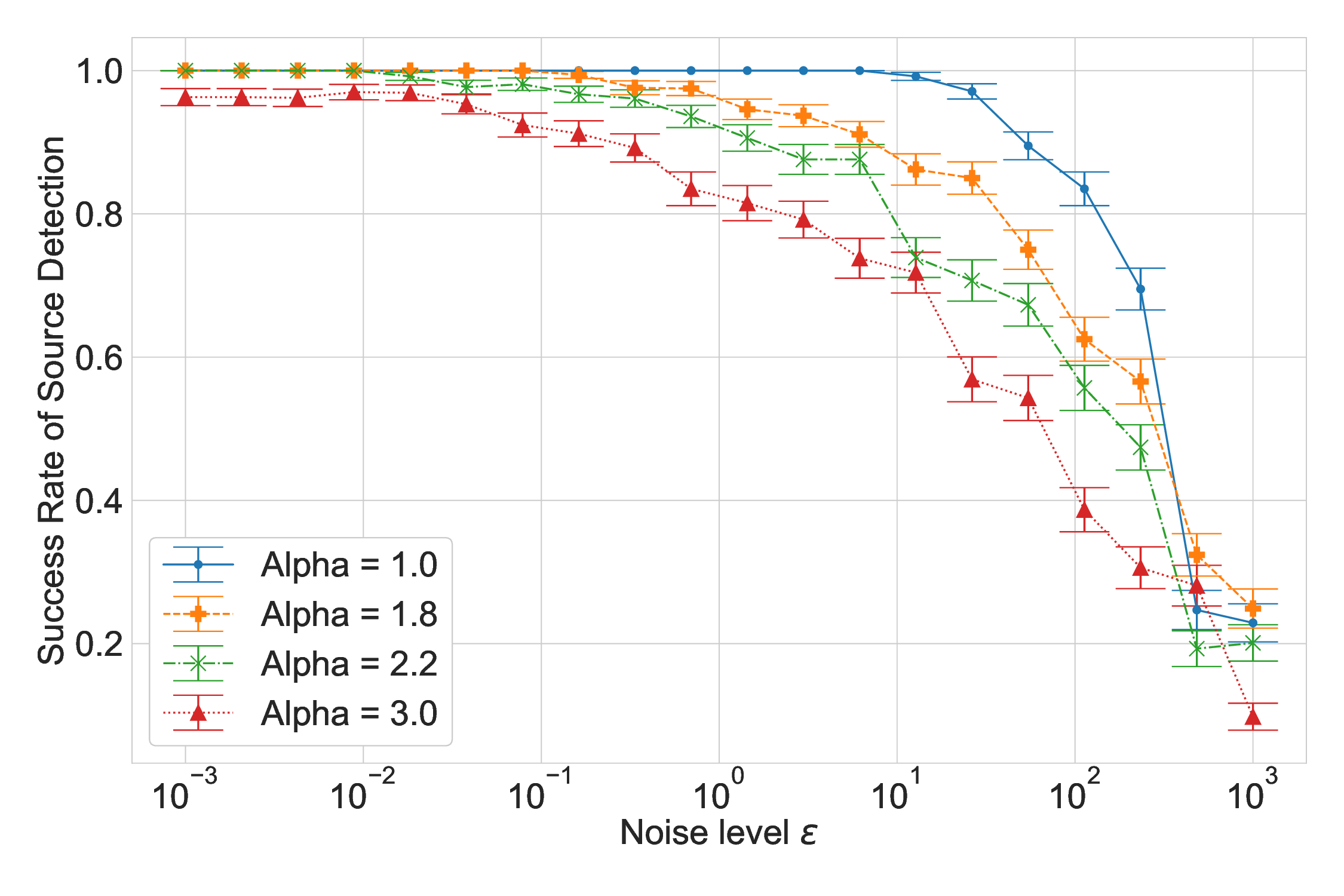}
  \caption{Predict top-$10$ most shifted contracts}
  \label{fig:equicorr-right}
\end{subfigure}
\caption{
{\em Source Detection Problem in a noisy scaled equi-correlation model of $\Sigma$:} 
We rank the entries of $W$ by the magnitude of change induced by a change in one entry of $M$ ($M_{ij}$).
Plot (a) shows the fraction of times $W_{ij}$ is most-changed entry of $W$.
Plot (b) shows the fraction of times $W_{ij}$ is among the top-$10$ most changed entries of $W$.
The success rate goes to zero as $\alpha$ and $\eps$ increase.
%We set $n = 50$ and $\rho = 0.1$. Results are averaged over $1000$ repetitions.
% \txtred{We should say ``Results are averaged over $1000$ repetitions'' in the caption, not in the figure legend.}
}
\label{fig:equicorr}
\end{figure}

\smallskip\noindent
{\bf Empirical Results for the Source Detection Problem (Simulated Data).}
Here, we set the covariance $\Sigma=D^{1/2} (R + \mathcal{E}) D^{1/2}$, where $D$ is a diagonal matrix, $R$ a correlation matrix, and $\mathcal{E}$ a noise matrix.
If $\mathcal{E}=0$, then $D_{ii}$ would be the variance of firm $i$.
We set $D_{ii}$ according to a power law: $D_{ii} = i^{-\alpha}$ for an $\alpha > 0$.
Larger values of $\alpha$ correspond to greater skew in the variances.
We choose $R$ to be an equi-correlation matrix with $1$ along the diagonal and $\rho \in (0, 1)$ everywhere else.
We draw the error matrix $\mathcal{E}$ from a scaled Wishart distribution: $\mathcal{E} = \norm R\norm_2 \cdot \mathcal{W}(\sqrt{\eps}\cdot I_n, n) / n$ for some chosen the noise level $\eps$.
As $\eps$ increases, the noise $\mathcal{E}$ dominates $R$.

Figure \ref{fig:equicorr} shows the success rate of source detection over $1000$ experiments for various values of $(\eps, \alpha)$ for $\rho = 0.1$ and $n=50$. 
As $\alpha$ increases, the variances become more skewed and the source detection can fail even with $\eps=0$ noise. 
When $\eps$ grows, the success rate for the source detection problem goes to zero.
This suggests that skew combined with noise makes source detection difficult. 
These trends occur even if we only test whether the source belongs to the $10$ most changed contracts (Figure \ref{fig:equicorr-right}), as opposed to single largest change (Figure \ref{fig:equicorr-left}).
We observe similar results for real-world choices of $\Sigma$, as we show next.

%This is because the eigenvalues of $\Sigma$ exhibit greater skew at large $\alpha$, and at high $\eps$ the additive noise dominates the signal.
%Next, we show empirical results that source detection and targeted interventions are much harder for realistic choices of $\Sigma$.

\smallskip\noindent
{\bf Empirical Results for the Source Detection Problem (Real-World Data).}
We consider two datasets:
(a) a trade network between $46$ large economies~\citep{oecd-stats}, and
(b) a simulated network between $96$ portfolio managers following various Fama-French strategies~\cite{fama-french-2015}.
For each dataset, we construct a ``ground-truth'' covariance $\Sigma$ using all available data (the details are in in Appendix~\ref{appendix-datasets}). 
Then, using $m$ independent samples $\bmx_i\sim\mathcal{N}(0, \Sigma)$, we build a ``data-driven'' covariance $\hat{\Sigma} = (1/(m - 1)) \sum_{i=1}^m (\bmx_i - \bm{\hat{\mu}}) (\bmx_i - \bm{\hat{\mu}})^T$, where $\bm{\hat{\mu}} = (1/m)  \sum_{i=1}^m \bmx_i$ is the sample mean. 
We use this $\hat{\Sigma}$ to construct the financial network.
% \txtred{Did we convert to the reviewer's version of covariance estimation?}

\begin{figure}[tbp]
\centering
\begin{subfigure}{0.48\textwidth}
  \includegraphics[width=\textwidth]{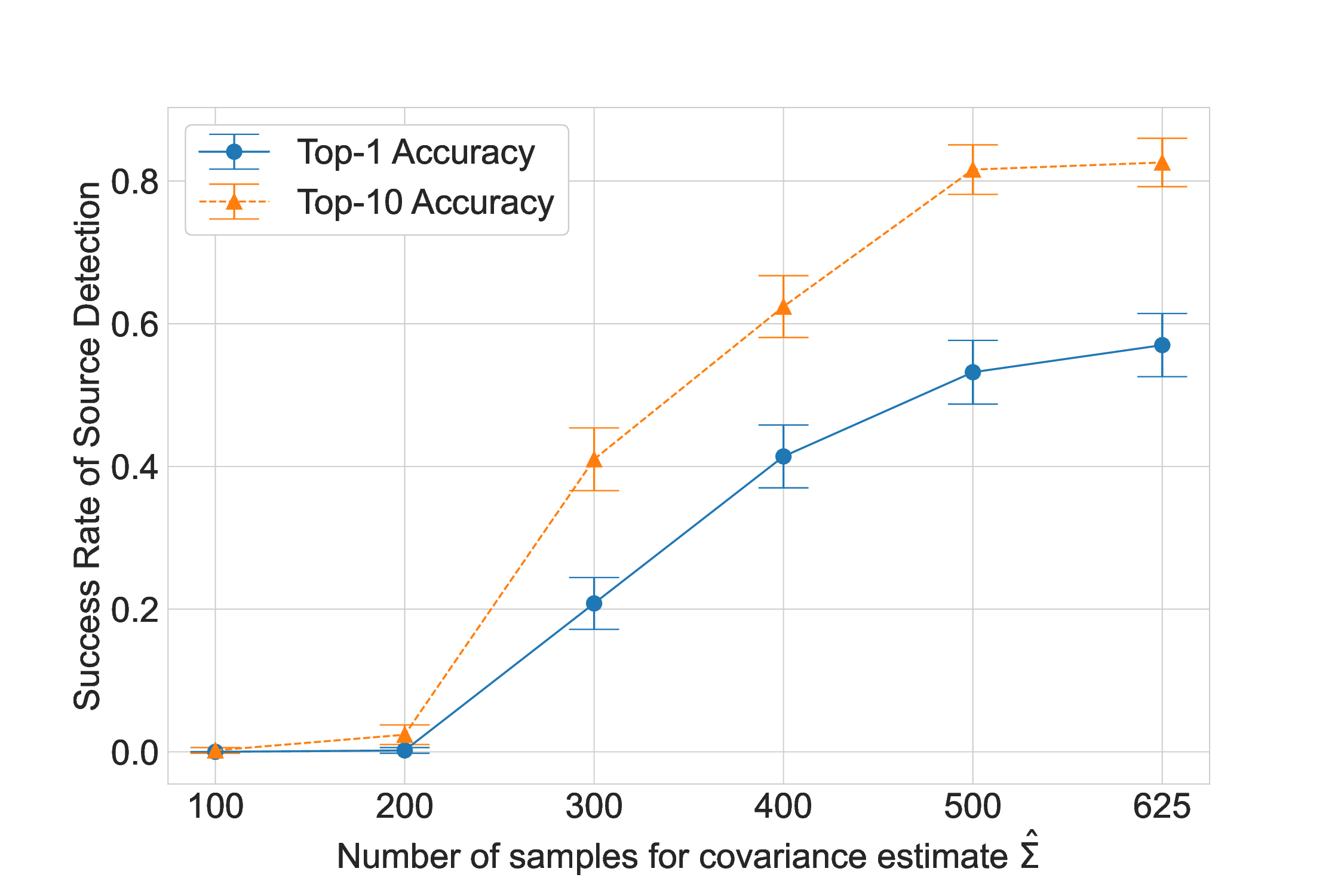}  
  \caption{Simulated network of $96$ portfolio managers.}
\end{subfigure}
\hfill
\begin{subfigure}{0.48\textwidth}
\includegraphics[width=\textwidth]{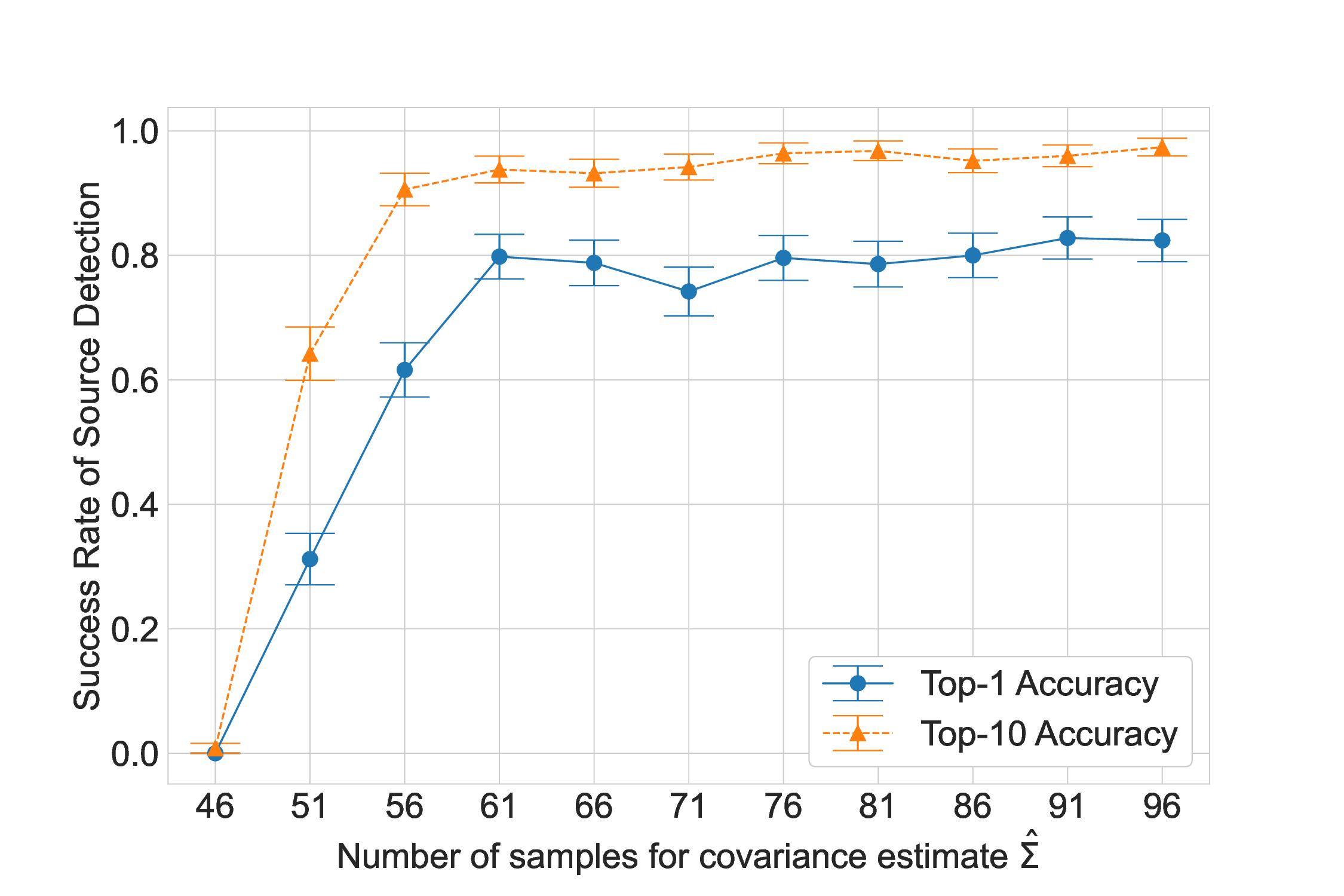}  
\caption{$46$-country (OECD) trade network.}
\end{subfigure}
\caption{
{\em Source Detection Problem on real-world data:}
The success rate scales monotonically with the number of samples used to construct the data-driven covariance matrix $\hat{\Sigma}$. 
%The top row shows the success rate of source detection for two datasets, and the bottom shows the success rate for top $10$ detection.}
% \txtred{Let's have one fig for portfolio and one for OECD. Each fig will have two plots, labeled "Top-1 accuracy" and "Top-10 accuracy", with solid and dashed lines and different colors. If the std-dev bars interfere, we'll remove them and just state the std-dev numbers in the caption. Make sure axes labels and legends are big and legible.}
}
\label{fig:meanShock}
\end{figure}

\begin{figure}[tbp]
\centering
\begin{subfigure}{0.5\textwidth}
  \includegraphics[width=\textwidth]{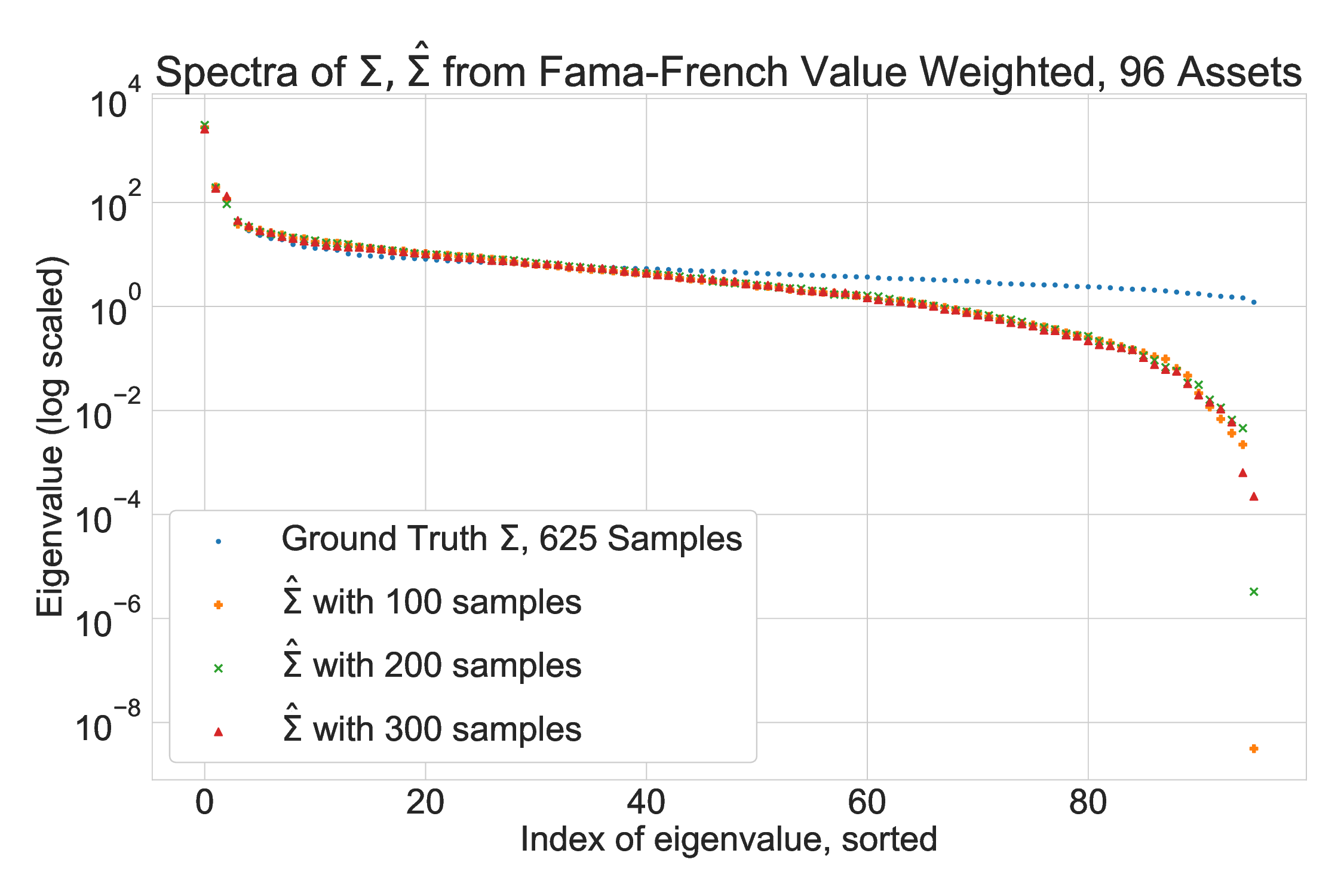}
  \caption{Simulated network of $96$ portfolio managers}
\end{subfigure}%
\begin{subfigure}{0.5\textwidth}
  \includegraphics[width=\textwidth]{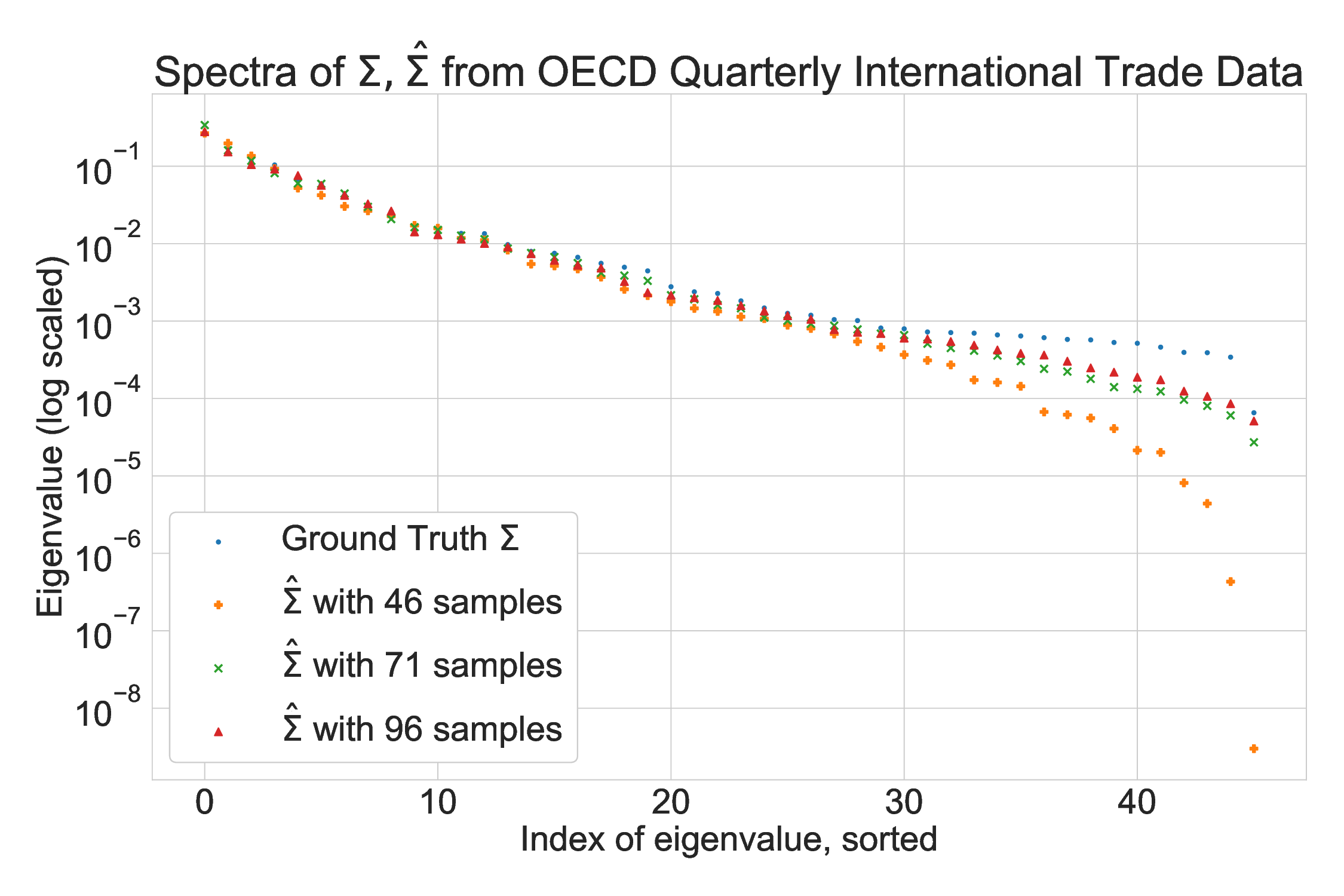}
  \caption{$46$-country (OECD) trade network}
\end{subfigure}
\caption{
  The eigenvalues of estimated covariance matrices are skewed, and the degree of skew depends on the number of samples $m$.
As $m$ decreases, so does the smallest eigenvalue $\lambda_n$ and the ratio $\lambda_n/\lambda_{n-1}$.
}
\label{fig:spectra}
\end{figure}

Figure~\ref{fig:meanShock} shows the success rate over $500$ experiments for various choices of the sample size $m$.
The success rate increases monotonically with $m$. 
%\blue{As in Figure \ref{fig:equicorr}, this trend holds even if we only require the source to lie in the $10$ most changed contracts, albeit with higher success rates for each $m$.}
The reason for this behavior lies in the spectra of $\Sigma$ and $\hat{\Sigma}$.
We find that in both datasets, the largest and smallest eigenvalues of $\Sigma$ are separated by several orders of magnitude.
This gap becomes even more extreme in the data-driven $\hat{\Sigma}$; the fewer the samples $m$, the greater the gap (see Figure \ref{fig:spectra}).
In fact, we observe that the smallest eigenvalue of $\hat{\Sigma}$ is much smaller than the second-smallest eigenvalue: $\lambda_n\ll \lambda_{n-1}$.
\cite{zhao-2019-bounded-noise-portfolio} make similar observations.

In summary, the experiments on both simulated and real-world datasets highlight the difficulty of source detection and targeted intervention in realistic networks.
The reason is the skew in the eigenvalues coupled with noise, which affects the eigenvectors.
Skewed eigenvalues correspond to trade combinations (eigenvectors) that are seemingly low-risk.
Hence, firms use such trades to diversify.
This implies that these eigenvectors have an outsized effect on the network, and how it responds to local changes.
Intuitively, if these eigenvectors are ``random,'' the effect of a changed belief $M_{k\ell}$ affects the rest of the network randomly.
Hence, the direct effects on $W_{k\ell}$ may be less than the indirect effects on other $W_{ij}$.
We explore this theoretically in Appendix~\ref{appendix:dwdm:approx}.

% !TEX root = ./paper-draft.tex

\section{Insights for Firms} % in a Community-Structued Setting}
\label{sec:firms}

Until now, we have treated the beliefs of firms as fixed and exogenous.
In this section, we consider how a firm can use its contracts to gain insights into other firms and update its beliefs.

For instance, suppose a firm $j$ faces a crisis, e.g., a looming debt payment that may make it insolvent.
The firm may then become risk-seeking (i.e., lower its $\gamma_j$), hoping that the risks pay off.
Another firm $i$ may be unaware of the crisis, so $i$'s risk perceptions (perhaps based on historical data) would be outdated.
Can firm $i$ {\em infer} the lower $\gamma_j$, solely from $i$'s contracts $\bmw_i$ with all firms?
What if a group of firms become risk-seeking, and not just one firm?

%\subsection{Distribution of contract sizes}
\subsection{Detecting Outlier Firms}
\label{sec:firms:single}

Intuitively, firm $i$ will try to answer these questions by comparing the behavior of firm $j$ against other similar firms.
We formalize this by assuming that each firm $j$ belongs to a community $\theta_j$, e.g., banking, or real-estate, or insurance, etc.
The community of each firm is publicly known.
Firms in the same community are perceived to have similar return distributions:
\begin{align}
M_{ij} &= f(\theta_i, \theta_j) + \eps^\prime_{\theta_i,j}, &
\Sigma_{ij} &= g(\theta_i, \theta_j), &
\gamma_i &= h(\theta_i) + \eps_i
\label{eq:community}
\end{align}
for some unknown deterministic functions $f(.)$, $g(.)$, and $h(.)$ and random error terms $\eps_i$ and $\eps^\prime_{\theta_i,j}$.
We also assume that all firms use the same covariance $\Sigma$.

Now, suppose one firm $j$ is an outlier, with very different beliefs from other firms in its community. 
For firm $i$ to detect the outlier firm $j$, the contract size $W_{ij}$ should deviate from a cluster of contracts $\{W_{ij'}\mid \theta_{j'} = \theta_j\}$ of other firms from the same community as firm $j$.
%That is, we expect $W_{ij}$ to deviate from a cluster formed by $\{W_{ij'}\}$.
Now, outlier detection methods often assume independent datapoints.
In our model, all contracts are dependent.
But we can still do outlier detection if the contracts are appropriately exchangeable.
We prove below this is the case. 

\begin{definition}
An intra-community permutation is a permutation $\pi: [n] \to [n]$ such that $\pi(i) = j$ implies that $\theta_i = \theta_j$. 
\end{definition}

\begin{proposition}
\label{prop:firms-exchangeable-stronger}
Suppose $M, \Sigma, \Gamma$ exhibit community structure (Eq.~\eqref{eq:community}), and all the error terms $(\eps_i)_{i\in [n]}$ and $(\eps^\prime_{\theta_i, j})_{i,j\in [n]}$ are independent and identically distributed. 
%also assume that if $\theta_r = \theta_s$ then $\eps_{r} \sim \eps_{s}$, and $\eps_{\theta_i, r}^\prime \sim \eps_{\theta_i, s}^\prime$.
Let $\pi: [n] \to [n]$ be any intra-community permutation, and let $\Pi: \RR^n \to \RR^n$ be the corresponding column-permutation matrix: $\Pi(\bm{e}_i) = \bm{e}_{\pi(i)}$. 
Then, $W$ and $\Pi^T W \Pi$ are identically distributed.  % $W\sim \Pi^T W \Pi$.
\end{proposition}

\begin{corollary}\label{cor:firms-exchangeable}
Let $j_1, \dots, j_m \in [n]$ belong to the same community: $\theta_{j_1} = \dots = \theta_{j_m}$. 
Suppose the conditions of Proposition~\ref{prop:firms-exchangeable-stronger} hold.
Then, for any $i \in [n]$, the joint distribution of $(W_{i, j_1}, \ldots, W_{i, j_m})$ is exchangeable.
\end{corollary}

% \begin{figure}[htb]
% \begin{center}
% \includegraphics[width=.8\linewidth]{figs/false-positives-091523/outlier_detection_all_res_top_5_complete.eps}
% \caption{
%   Rate of the outlier risk-seeking firm belonging to the top $5$ counterparties of the testing firm. As in Figure \ref{fig:deviator-gamma-detection}, detection is more successful when there are fewer firms and when the risk-seeking firm's $\gamma_{\textrm{outlier}}$ is more standard deviations away from the $\gamma$ of the normal firms.
% }
% \label{fig:deviator-gamma-detection-top5}
% \end{center}
% \end{figure}

{\bf Empirical Results for Outlier Detection.}
%As long as $\eps_i$ are small enough, 
We generate community-based networks (Eq.~\eqref{eq:community}) such that $\gamma_i\sim N(1,\sigma^2)$ truncated to $[0.5, 1.5]$.
The smaller the $\sigma$, the more closely the $\gamma_i$ values cluster around $1$.
For the outlier risk-seeking firm, we set $\gamma_{\textrm{outlier}}=0.5$.
For clarity of exposition, we set $\eps^\prime=0$ everywhere.

To detect outliers under exchangeability (Corollary~\ref{cor:firms-exchangeable}), we can use methods based on conformal prediction~\citep{guan-2022}.
Here, we use a simpler approach: pick the firm $j$ with the largest contract size as the outlier; $\hat j \defeq \arg\max\limits_{j \in \{j_1, \dots, j_m\}} \abs{W_{i, j}}$. 
To test sensitivity to false negatives, we also test whether the outlier is among the 5 largest contracts in $\{\abs{W_{i, j}}: j = j_1, \dots, j_m \}$. 
We run $500$ experiments for each choice of $\sigma$, and count the frequency with which the outlier firm is detected via its contract size.
Further details are presented in Appendix~\ref{appendix-gamma-exps}.

Figure \ref{fig:deviator-gamma-detection} shows the results.
We characterize the degree of outlierness by how many standard deviations away $\gamma_{\textrm{outlier}}$ is from the baseline of $1$.
The smaller the $\sigma$, the more the outlierness.
The success rate increases with increasing outlierness, as expected.
It also increases when the number of firms $n$ is reduced.
This is because contract sizes depend on the $\gamma$ values of all firms; fewer firms reduces the chances of any one firm attaining large contract sizes due to randomness.

\begin{figure}[tbp]
\centering
\begin{subfigure}{0.5\textwidth}
  \includegraphics[width=\textwidth]{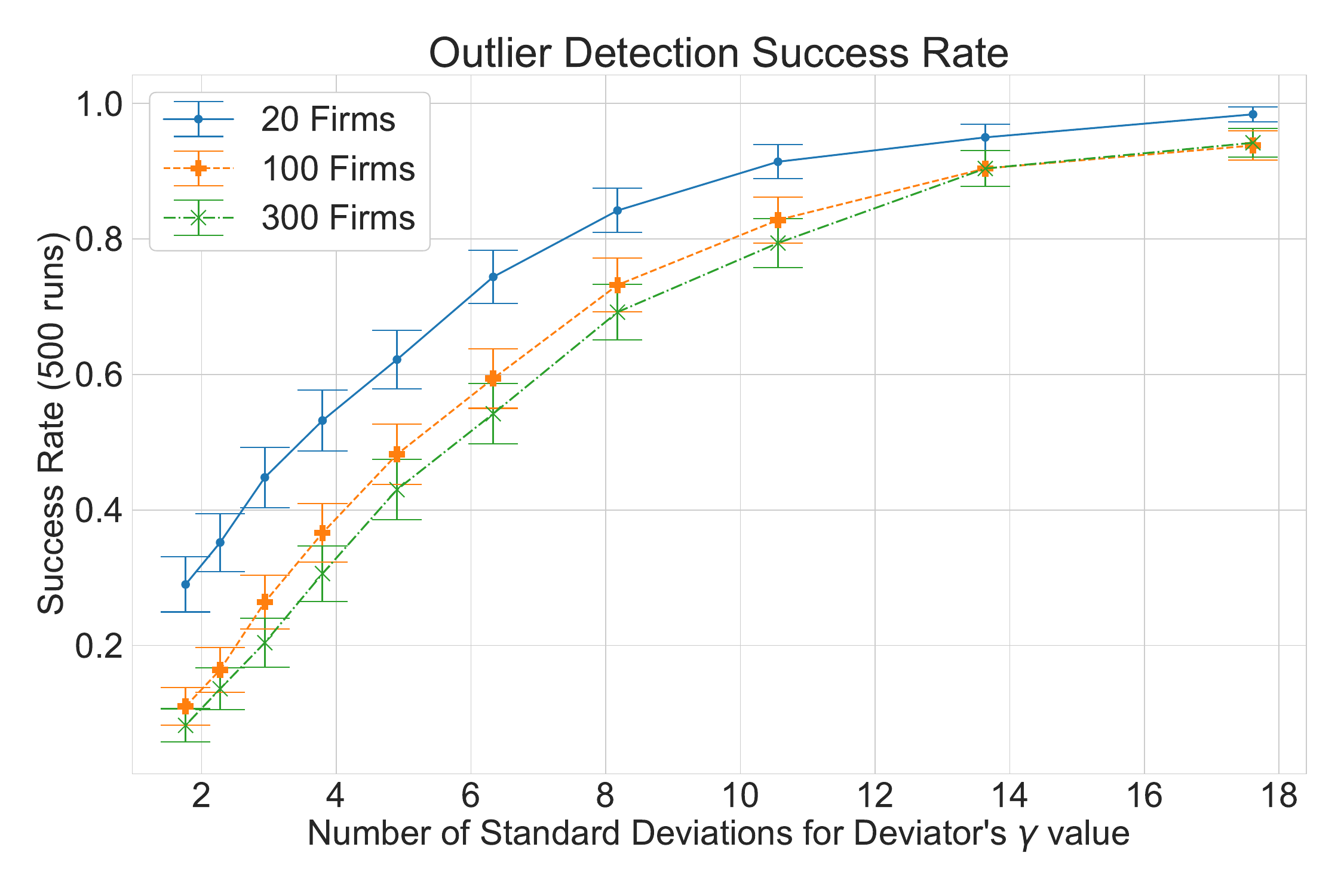}
%  \caption{Outlier has largest contract size.}
  \caption{Predict largest contract as outlier}
\end{subfigure}%
\begin{subfigure}{0.5\textwidth}
  \includegraphics[width=\textwidth]{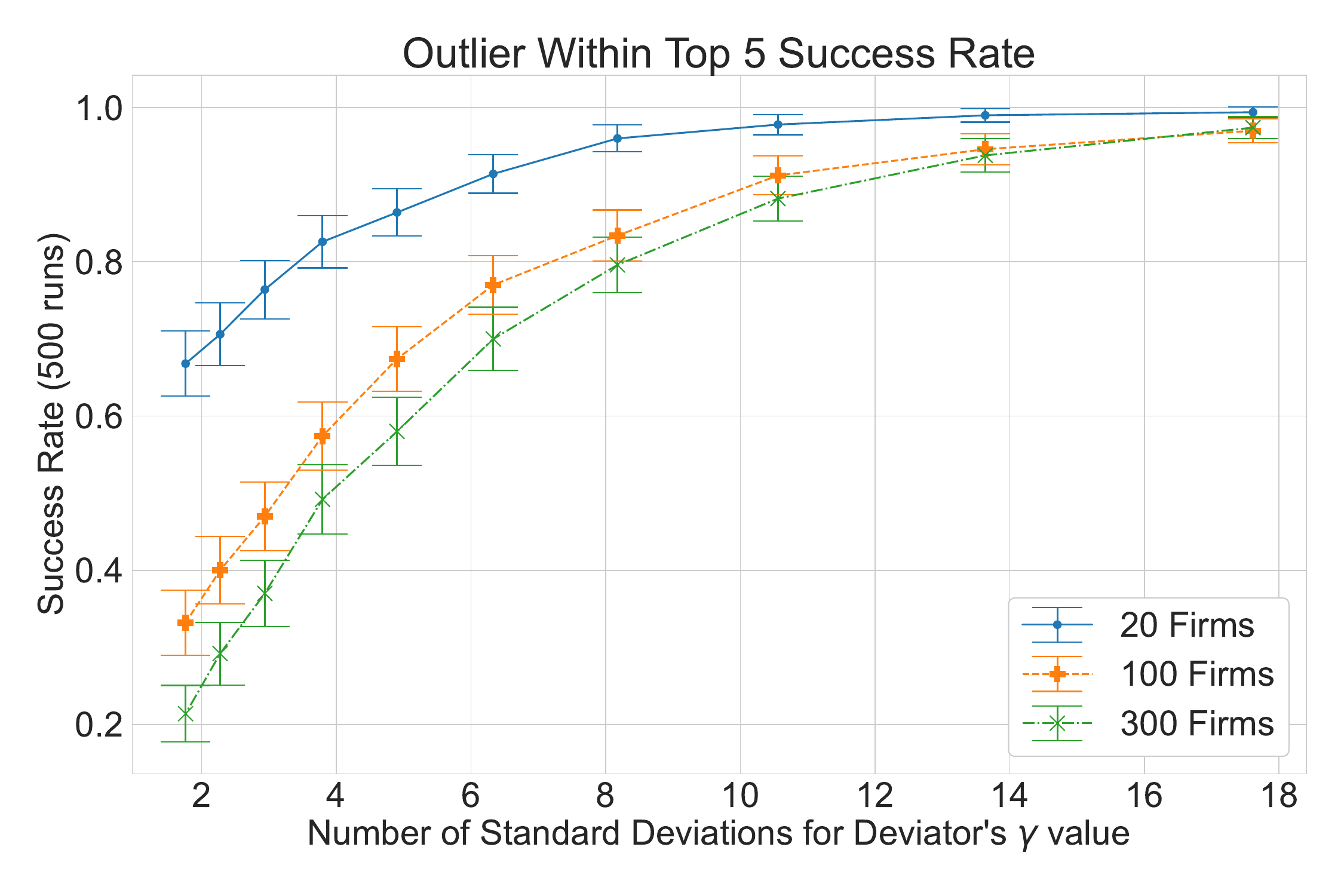}
%  \caption{Outlier is among the $5$ largest contracts.}
  \caption{Predict the top-$5$ largest contracts}
\end{subfigure}
\caption{
{\em Success rate for detecting outlier risk-seeking firms:}
  Detection is easier when there are fewer firms and when the risk-seeking firm's $\gamma_{\textrm{outlier}}$ is more standard deviations away from the $\gamma$ of the normal firms.
}
\label{fig:deviator-gamma-detection}
\end{figure}

\subsection{Risk-Aversion versus Expected Returns}
The discussion above shows that a firm can detect outlier counterparties.
However, the firm cannot determine {\em why} the counterparty is an outlier, as the following theorem shows.

\begin{theorem}[Non-identifiability of risk-aversion versus expected returns]
\label{thm:nonidentGammaM}
Consider two network settings $S=(\bmu_i, \Sigma, \gamma_i)_{i\in[n]}$ and $S^\prime=(\bmu_i, \Sigma, \gamma_i^\prime)_{i\in[n]}$ which differ only in the risk-aversions of firms $J=\{j\mid \gamma_j\neq \gamma^\prime_j\}\subseteq [n]$.
Then, there exists a setting $S^\dagger=(\bmu_i^\dagger, \Sigma, \gamma_i)_{i\in[n]}$ such that $\bmu_i=\bmu^\dagger_i$ for all $i\notin J$ and the stable networks under $S^\dagger$ and $S^\prime$ are identical.
\end{theorem}

Thus, one cannot determine if an outlier is more risk-seeking than its community or expects higher profits.
But risk-seeking behavior may be indicative of stress, while higher profits than similar firms are unlikely.
Hence, in either case, the firm detecting the outlier may choose to reduce its exposure to the outlier.
However, this approach fails if an entire community shifts its behavior.
The following example illustrates the problem.

%The previous section showed that outlier detection algorithms can find risk-seeking firm(s) hiding within a risk-averse community.
%But some conditions of financial stress may affect all firms in a community (e.g., increasing mortgage default rates affecting all real-estate firms).
%Then, the entire community becomes risk-seeking.
%A firm $i$ outside the community only observes a change in the contract sizes $W_{ij'}$ with firms $j'$ from the community.
%However, firm $i$ cannot determine the {\em cause} of the change, as the following theorem shows.

\begin{example}
Consider two communities numbered $1$ and $2$, with $n_1$ and $n_2$ firms respectively.
Let the setting $S$ of Theorem~\ref{thm:nonidentGammaM} correspond to
\begin{align*}
M_{ij} &= \left\{\begin{array}{cl} a & \text{if $\theta_i=\theta_j=1$}\\ b & \text{if $\theta_i=\theta_j=2$}\\ c/2 & \text{otherwise}\end{array} \right.
  & \Sigma_{ij} &= \left\{\begin{array}{cl} 1 & \text{if $\theta_i=\theta_j=1$}\\ 1 & \text{if $\theta_i=\theta_j=2$}\\ 0 & \text{otherwise}\end{array} \right.
  & \gamma_i &= 1.
\end{align*}
%Note that typically $a>0$ and $b>0$ since firms should expect positive rewards from investing in their own business.
%Then, using Corollary~\ref{cor:stability:sharedSigmaGamma}, we find the stable network to be
%\begin{align*}
%W_{ij} &= \left\{\begin{array}{cl} a/(2n_1) & \text{if $\theta_i=\theta_j=1$}\\ b/(2n_2) & \text{if $\theta_1=\theta_2=2$}\\ c/(n_1+n_2) & \text{otherwise}\end{array} \right.
%\end{align*}
Now, suppose that under setting $S^\prime$, $\gamma_i\mapsto\gamma_i+\delta$ for some small $\delta$ for all nodes $i$ in community $1$.
The change in the network would be the same if we had updated the columns corresponding to community~$1$ in the $M$ matrix instead (setting $S^\dagger$):
\begin{align*}
M^\dagger_{ij} &= M_{ij} + \Delta(\theta_i, \theta_j) &
\Delta(\theta_i, \theta_j) + O(\delta^2) &= \left\{\begin{array}{cl} -\delta a/2 & \text{if $\theta_i=\theta_j=1$}\\ -\delta b\cdot n_2/(n_1+n_2) & \text{if $\theta_i=2, \theta_j=1$}\\ 0 & \text{if $\theta_j=2$}\end{array} \right.
%(1/\delta)\cdot(W^\prime - W) &= \left\{\begin{array}{cl} a & \text{if $\theta_i=\theta_j=1$}\\ b/(2n_2) & \text{if $\theta_1=\theta_2=2$}\\ c/(n_1+n_2) & \text{otherwise}\end{array} \right.
\end{align*}
%Recall that the $j^{th}$ column of $M$ denotes the expected return beliefs of firm $j$.
Thus, a firm from community $2$ cannot determine if the network change was due to a change in $(\gamma_i)_{\theta_i=1}$ or $(\bmu_i)_{\theta_i=1}$.
For instance, when $b>0$, an increase in risk-seeking ($\delta<0$) looks the same as an increase in trading benefits ($\Delta(1,2)>0$).
In the former case, firms in community~$2$ should {\em reduce} their exposure to community~$1$ firms. 
But in the latter case, they should {\em increase} exposure.
%But firms in community~$2$ should {\em reduce} their exposure to community~$1$ in the former case, \textcolor{blue}{and} {\em increase} it in the latter case.
Since the data cannot be used to choose the appropriate action, the behaviors of firms may be guided by their prior beliefs or inertia.
When such beliefs change due to external events (e.g., due to news about one firm in community~$1$), the resulting change in the network may be drastic.$\hfill\Box$
\end{example}

\section{Conclusions}
\label{sec:conc}
We have proposed a model of a weighted undirected financial network of contracts.
The network emerges from the beliefs of the participant firms.
The link between the two is utility maximization coupled with pricing.
For almost all belief settings, our approach yields a unique network.
This network satisfies a strong Higher-Order Nash Stability property.
Furthermore, the firms can converge to this stable network via iterative pairwise negotiations.

The model yields two insights.
First, a regulator is unable to reliably identify the causes of a change in network structure, or engage in targeted interventions.
The reason is that firms seek to diversify risk by exploiting correlations.
We find that in realistic settings, there are often combinations of trades that offer seemingly low risk.
Hence, all firms aim to use such trades.
The over-dependence on a few such combinations leads to a pattern of connections between firms that thwarts targeted regulatory interventions.

%some the smallest eigenvector of the perceived risk matrix appears to offer a particularly low-risk solution.
%The over-dependence of the network on this single eigenvector leads to a pattern of connectivity between firms that thwarts targeted interventions.

%In realistic data-driven settings, the smallest eigenvector of the perceived risk matrix appears to offer a particularly low-risk solution.

The second insight is that firms can use the network to update their beliefs. 
For instance, they can identify counterparties that behave very differently from their peers.
However, the cause of the outlierness remains hidden.
If all firms in one line of business become more risk-seeking, the result is indistinguishable from that business becoming more profitable.
Innocuous events (such as a news story) may cause beliefs to change suddenly, leading to drastic changes in the network. 
In addition to identifying risky counterparties, firms may use the network to update their mean and covariance beliefs.
For example, a firm that suffers significant losses on its current trades may be judged by others to be a riskier counterparty for future trades.
We leave this for future work.

%For instance, they can identify counterparties that have become more risk-seeking than usual.
%However, if all firms in one line of business , the result is indistinguishable from that business becoming more profitable.
%Since the network cannot help firms decide, firms must rely on their beliefs.

Our work focuses on mean-variance utility, but some of our results are applicable in other settings too.
A second-order Taylor approximation of a twice-differentiable concave utility matches the form of a mean-variance utility.
Hence, results based on mean-variance utility can be useful guides for small perturbations around a stable point.
Some of our results for pairwise negotiations and targeted interventions are based on such perturbation arguments.

%While our work focuses on mean-variance utility, many of our results hold more generally.
%This is because the mean-variance utility is a good approximation of a perturbation of \blue{the utility function $f$ around the stable point, as long as $f$ is real analytic at the stable point}. 
%Some of our results for pairwise negotiations and targeted interventions are based on such perturbation arguments.
%Hence, these generalize to other choices of utility functions.

\section{Acknowledgments}

The authors thank Stathis Tompaidis, Marios Papachristou, Kshitij Kulkarni, and David Fridovich-Keil for valuable discussions and suggestions.
This work was supported by NSF grant 2217069, a McCombs Research Excellence Grant, and a Dell Faculty Award.

\bibliographystyle{informs2014} 
\bibliography{refs, networks} % if more than one, comma separated

\begin{APPENDICES}
% !TEX root = ./paper-draft.tex
\newcommand{\beginpfbare}{\proof{Proof.}}
% \tableofcontents

\section{Proofs}

\subsection{Proof of Theorem \ref{thm:stable}}

Recall that $Q_i = \Psi_i^T (2 \gamma_i \Psi_i \Sigma_i \Psi_i^T)^{-1} \Psi_i$, $F = \{(i, j): 1 \leq i < j \leq n, \Psi_i \bm{e}_j \neq \bm{0}\}$, and $\uvc(X)_F\in\mathbb{R}^{|F|}$ is a vector whose entries are the ordered set $\{X_{ij}\mid (i,j)\in F\}$.  

Note that $\Psi_i \Sigma_i \Psi_i^T$ is positive definite, since it is a principal submatrix of the positive definite matrix $\Sigma_i$. 
We shall prove an expanded version of Theorem~\ref{thm:stable}.

\begin{repeattheorem}[Expanded Version of Theorem \ref{thm:stable}].
Define $n\times n$ matrices $A$, $B_{(i,j)}$, and $C_{(i,j)}$ as follows:
\begin{align*}
A_{ij} &=\bm{e}_i^T Q_j M \bm{e}_j, &
B_{(i,j)} &= \bme_i\bme_j^T Q_i, &
C_{(i,j)} &= (B_{(i,j)} - B_{(j,i)}) - (B_{(i,j)} - B_{(j,i)})^T.
\end{align*}
Let $Z_F$ be the $|F|\times |F|$ matrix whose rows are the ordered sets $\{\uvc(C_{(i,j)})_F \mid (i,j)\in F\}$.
Then, we have the following:
\begin{enumerate}
\item A stable point $(W, P)$ under $\{\Psi_i\}$ exists if and only if $\uvc(A-A^T)_F$ lies in the column space of $Z_F$.
\item If a stable point $(W, P)$ exists, then $Z_F \uvc(P)_F = \uvc(A - A^T)_F$. 
\item A unique stable point always exists if $Z_F$ is full rank.
\end{enumerate}
\end{repeattheorem}

\proof{Proof.}
For clarity of exposition, we first prove the result when all edges are allowed, and then consider the case of disallowed edges.

\noindent{\bf (1) All edges allowed.}
Here, $E=\{i,j\mid 1 \leq i < j\leq n\}$, and we use $\uvc(.)$ and $Z$ to refer to $\uvc(.)_E$ and $Z_E$ in the theorem statement.
%{\em Item (1).} 
For any price matrix $P$ with $P=-P^T$, consider the matrix $W$ whose $j^{th}$ column has the utility-maximizing contract sizes for agent $j$:
\begin{align*}
W_{ij} &= \bm{e}_i^T \Psi_j^T (2 \gamma_j \Psi_j \Sigma_j \Psi_j^T)^{-1} \Psi_j (M - P)\bm{e}_j = \bm{e}_i^T Q_j (M - P) \bm{e}_j.
\end{align*}
The tuple $(W, P$) is stable if $W=W^T$.
So, for all $i<j$, we require
\begin{align}
W_{ij} &= W_{ji} \\
\Leftrightarrow \bme_i^T Q_j (M - P) \bme_j &= \bme_j^T Q_i (M - P) \bme_i \nonumber\\
\Leftrightarrow \bme_i^T Q_j M \bme_j - \bme_j^T Q_i M \bme_i &= \bme_i^T Q_j P \bme_j - \bme_j^T Q_i P \bme_i \nonumber\\ 
\Leftrightarrow \bme_i^T (A-A^T) \bme_j &= \bme_i^T (Q_j P - (Q_i P)^T) \bme_j.\label{eq:stability1}
\end{align}
Since $P=-P^T$, we must have $P=R - R^T$, where $R$ is upper-triangular with zero on the diagonal.
Hence, using $Q_i=Q_i^T$, we have
\begin{align*}
\bme_i^T (Q_j P - (Q_i P)^T) \bme_j
    &= \bme_i^T (Q_j P + P Q_i) \bme_j\\
    &= tr P (\bme_j\bme_i^T Q_j + Q_i\bme_j\bme_i^T)\\
    &= tr (R - R^T) (B_{(j,i)} + B_{(i,j)}^T)\\
    &= tr R^T C_{(i,j)}\\
    &= \uvc(R)^T \uvc(C_{(i,j)}),
\end{align*}
where we used the upper-triangular nature of $R$ in the last step.
Plugging into Eq.~\eqref{eq:stability1}, a stable point exists if and only if there is an appropriate vector $\bm{p}:=\uvc(R)\in\RR^{n(n-1)/2}$ such that for all $1 \leq i < j \leq n$, $\bme_i^T (A-A^T) \bme_j = \uvc(C_{(i,j)})^T \bm{p}$. This is equivalent to $\uvc(A-A^T) = Z \bm{p}$. 
If such a solution vector $\bm{p}$ exists, then by definition it corresponds to a matrix $P = -P^T$ via $P = R - R^T$ and ${\bm p} = \uvc(R)$.

\smallskip\noindent{\bf (2) Disallowed edges. }
If $\{i, j\}$ is a prohibited edge then $\Psi_i \bm{e}_j = \Psi_j \bm{e}_i = \bm{0}$, so $B_{(i, j)} = B_{(j, i)} = 0$, so $\bm{e}_{ij}^T Z = \bm{0}^T$. Also, $A_{ij} = A_{ji} = 0$ so $\uvc(A - A^T)_{ij} = 0$. Therefore, the equality $\bme_i^T (A-A^T) \bme_j = \uvc(C_{(i,j)})^T \bm{x}$ is achieved for any solution vector $\bm{x}$ if $\{i, j\}$ is a prohibited edge. We can therefore reduce the linear system $Z \bm{p} = \uvc(A-A^T)$ from part (1) by deleting rows of $Z$ corresponding to prohibited edges. 

Similarly, since the system is constrained by $\bm{p}_{ij} = 0$ for prohibited edges $\{i, j\}$, the columns of $Z$ corresponding to such edges have no effect on the solution set. 

We conclude that the linear system in (1) is equivalent to the (unconstrained) reduced system $Z_F \bm{p}_F = \uvc(A - A^T)_F$.
Each solution $\bm{p}_F$ corresponds to a skew-symmetric $P$ by construction.
Finally, if $Z_F$ has full rank then the unique reduced solution is $\bm{p}_F = Z_F^{-1} \uvc(A - A^T)_F$.  
\Halmos
\endproof

\subsection{Example of Stable Network}\label{appendix:example-thrm1}

To illustrate Theorem \ref{thm:stable}, consider the following example.

\begin{example}[Stable points]
Consider a $3$-firm network where the only allowed edges are given by $F = \{(1, 2), (1, 3)\}$.
%contracts between firms $2$ and $3$ are prohibited. Hence the permitted edges are given by $F = \{(1, 2), (1, 3)\}$. 
Suppose firms share the same covariance belief matrix $\Sigma_1 = \Sigma_2 = \Sigma_3 = \Sigma$, but have different mean beliefs $M = \begin{bmatrix} \bmu_1 & \bmu_2 & \bmu_3 \end{bmatrix}$ and risk aversions. The firms' beliefs are: 
\[
M = \begin{bmatrix} 0 & 2/3 & 1/2 \\ 1 & 0 & 0 \\ 1 & 0 & 0 \end{bmatrix}, 
\Sigma = \begin{bmatrix} 1 & 1/2 & 1/2 \\ 1/2 & 1 & 1/2 \\ 1/2 & 1/2 & 1 \end{bmatrix}, 
\gamma_1 = 1, \gamma_2 = 1/2, \gamma_3 = 1/4
\]
Then the $A, B_{(i,j)}$ matrices in Theorem \ref{thm:stable} are given as:
\[
A = \begin{bmatrix} 0 & 2/3 & 3/4 \\ 1/2 & 0 & 0 \\ 1/2 & 0 & 0 \end{bmatrix},
B_{(1, 2)} = \begin{bmatrix}
0 & 3/4 & -1/4 \\
0 & 0 & 0 \\
0 & 0 & 0 
\end{bmatrix}, 
B_{(1, 3)} = \begin{bmatrix}
0 & -1/4 & 3/4 \\
0 & 0 & 0 \\
0 & 0 & 0 
\end{bmatrix}, 
B_{(2, 1)} = \begin{bmatrix}
0 & 0 & 0 \\
1 & 0 & 0 \\
0 & 0 & 0 
\end{bmatrix}, 
B_{(3, 1)} = \begin{bmatrix}
0 & 0 & 0 \\
0 & 0 & 0 \\
3/2 & 0 & 0 
\end{bmatrix}
\]

Hence 
\[
C_{(1, 2)} = \frac 1 4 \begin{bmatrix}
0 & 7 & -1 \\
-7 & 0 & 0 \\
1 & 0 & 0 
\end{bmatrix}, 
C_{(1, 3)} =  \frac 1 4\begin{bmatrix}
0 & -1 & 9 \\
1 & 0 & 0 \\
-9 & 0 & 0 
\end{bmatrix}, 
\]

Therefore, $Z_F = \frac{1}{4} \begin{bmatrix} 7 & -1 \\ -1 & 9 \end{bmatrix}$ and $\uvc(A - A^T)_F = (1/6, 1/4)^T$. Since $Z_F$ is full-rank, there exists a unique stable point for this network setting. 
\label{example:thrm1}
\end{example}

% \subsection{Proof of Theorem \ref{thm:stable-common}}
\subsection{Stable Points are Common}\label{appendix:stable-common}

\begin{lemma}\label{lem:QZF}
Define $F$, $Z_F$ and $Q_i$ as in Theorem~\ref{thm:stable}.
%Let $Q_i^\prime$ be such that $\Psi_i^T \Psi_i Q_i^\prime \Psi_i^T \Psi_i = \Psi_i^T \Psi_i (Q_i + \beta I) \Psi_i^T \Psi_i$.
Let $Q_i^\prime$ be such that
\begin{align*}
(Q_i^\prime)_{j, k} = \begin{cases} (Q_i)_{j, k} + \beta & \text{if } j=k, (i,j)\in F \\ (Q_i)_{j, k} & \text{otherwise} \end{cases}
\end{align*}
Then, the corresponding $Z_F^\prime$ has the form $Z_F^\prime = Z_F + \beta I$.
\end{lemma}
\proof{Proof of Lemma~\ref{lem:QZF}.}
This follows from the form of the matrices $B_{(i,j)}$ and $C_{(i,j)}$ in the statement of Theorem~\ref{thm:stable}.
\Halmos
\endproof

Now, we consider the $\Sigma_i$'s (and hence the $Q_i$'s) to be random variables.
Any distribution of $\{\Sigma_i\}_{i\in [n]}$ induces a distribution on $\{Q_i\}_{i\in [n]}$, where $Q_i\succ 0$.
Define $\tilde{Q}_i:=Q_i-\delta I$, where $\delta>0$ is the minimum of union of the (nonzero) eigenvalues of all the $Q_i$'s.
A distribution over $\{Q_i\}$ corresponds to a distribution over $(\{\tilde{Q}_i\}, \delta)$.

\begin{proposition}\label{prop:perturb-inv-2}
If the distribution of $\delta$ given $\{\tilde{Q}_i\}_{i \in [n]}$ is continuous, then a unique stable point exists with probability $1$.
\end{proposition}
\proof{Proof of Proposition~\ref{prop:perturb-inv-2}.}
Let $\tilde{Z}_F$ be the $|F|\times|F|$ matrix generated from $\{\tilde{Q}_i\}_{i\in [n]}$, and $Z_F$ the corresponding matrix for $\{Q_i\}_{i\in [n]}$.
By Lemma~\ref{lem:QZF}, $Z_F = \tilde{Z}_F + \delta I$.
Hence, $\sigma(Z_F) = \sigma(\tilde{Z}_F) + \delta$, where $\sigma(M)$ denote the set of eigenvalues of $M$.
Since $\sigma(\tilde{Z}_F)$ is a function of $\{\tilde{Q}_i\}$ and $\delta$ is continuous given $\{\tilde{Q}_i\}$, the eigenvalues of $Z_F$ are non-zero with probability $1$.
Hence, by Theorem~\ref{thm:stable}, a unique stable point exists for $\{Q_i\}$ with probability $1$.
\Halmos
\endproof

Note that we require no condition on the distribution of $\{\tilde{Q}_i\}$. The condition of Proposition~\ref{prop:perturb-inv-2} is satisfied if the joint distribution of the $\{\Sigma_{i}\}_{i \in [n]}$ is continuous and all edges are permitted, as shown in the following example.

\begin{example} Fix some $n \geq 2$. 
Suppose the  joint distribution of the $\{\Sigma_{i}\}_{i \in [n]}$ is continuous and all edges are permitted. Then $Q_i = (2 \gamma_i)^{-1} \Sigma_i^{-1}$ so the joint distribution of $\{Q_{i}\}_{i \in [n]}$ is continuous. 
By Bayes' rule, $\PP[\delta | \Tilde{Q}_1, \dots, \Tilde{Q}_n] \propto \PP[\delta, \Tilde{Q}_1, \dots, \Tilde{Q}_1] = \PP[Q_1, \dots, Q_n]$. Since $\PP[Q_1, \dots, Q_n]$ is continuous, we conclude $\PP[\delta | \Tilde{Q}_1, \dots, \Tilde{Q}_n]$ is continuous.

% distribution of $\{\tilde{Q}_i\}_{i \in [n]}$ is parametrized by some $\bm{\theta} \in \Theta$, for a suitable parameter space $\Theta \subseteq \RR^d$. Suppose the distribution $\PP[\Tilde{Q}_1, \dots, \Tilde{Q}_n | \theta]$ is continous. 

% Then the joint distribution $\{Q_i\}$ is parametrized by the pair $(\bm{\theta}, \delta)$. Moreover, since $\delta$ is a continuous function of $Q_1, \dots, Q_n$ (it is a minimum of eigenvalues), we know the distribution $\PP[Q_1, \dots, Q_n | \delta, \theta]$ is continuous. 

% Hence by Bayes' rule, $\PP[Q_1, \dots, Q_n | \delta, \theta] \propto \PP[\Tilde{Q}_1, \dots, \Tilde{Q}_n, \delta | \theta] = \PP[\delta | \Tilde{Q}_1, \dots, \Tilde{Q}_n] \cdot \PP[\Tilde{Q}_1, \dots, \Tilde{Q}_n | \theta]$. Hence $\PP[\delta | \Tilde{Q}_1, \dots, \Tilde{Q}_n]$ is continuous. 
\end{example}

\subsection{Stable Network for the Shared Covariance Case}\label{appendix:stability:sharedSigma}
%Proof of Corollary \ref{cor:stability:sharedSigma}}

In the case of a shared covariance matrix for all agents, we can give a closed form expression for the stable network. 

\begin{corollary}[Shared $\Sigma$, all edges allowed]
Suppose $\Sigma_i=\Sigma$ and $\Psi_i=I_n$ for all $i\in[n]$.
Let $(\lambda_i, \bmv_i)$ denote the $i^{th}$ eigenvalue and eigenvector of $\Gamma^{-1/2}\Sigma\Gamma^{-1/2}$.
Then, the network $W$ can be written in two equivalent ways:
\begin{align*}
\vc(W) &= \frac 1 2 (\Gamma\otimes \Sigma + \Sigma \otimes \Gamma)^{-1}\vc(M+M^T),\\
W &= \Gamma^{-1/2}\left( \sum_{i=1}^n \sum_{j=1}^n \frac{\bmv_i^T \Gamma^{-1/2}(M+M^T)\Gamma^{-1/2} \bmv_j}{2 (\lambda_i + \lambda_j)} \bmv_i \bmv_j^T \right) \Gamma^{-1/2}.
\end{align*}
The prices can be written as:
\begin{align*}
\vc(P) &= (\Gamma^{-1} \otimes \Sigma^{-1} + \Sigma^{-1} \otimes \Gamma^{-1})^{-1} \vc(\Sigma^{-1} M \Gamma^{-1} - \Gamma^{-1} M^T \Sigma^{-1}) \\
P &= \Gamma^{1/2}
\left( \sum_{i=1}^n \sum_{j=1}^n \frac{\bmv_i^T \Gamma^{1/2}(\Sigma^{-1} M \Gamma^{-1} - \Gamma^{-1} M^T \Sigma^{-1})\Gamma^{1/2} \bmv_j}{\lambda_i^{-1} + \lambda_j^{-1}} \bmv_i \bmv_j^T \right) 
\Gamma^{1/2}.
\end{align*}
% P= \sum_{i=1}^n \sum_{j=1}^n 
% \frac{\lambda_j^{-1} \bmv_j^T M \Gamma^{-1} \bmv_i - \lambda_i^{-1} \bmv_j^T \Gamma^{-1} M^T \bmv_i}{\lambda_i + \lambda_j} 
% \bmv_i \bmv_j^T 

\end{corollary}\label{cor:stability:sharedSigma}

\proof{Proof.}
We first prove the identity with $\vc(W)$. 

For each agent $i$ the optimal set of contracts is given as $\bm{w}_i = (2 \gamma_i \Sigma_i)^{-1} (M - P)\bm{e}_i$. Since $\Sigma_i = \Sigma$ for all $i$, we obtain $W = \frac{1}{2} \Sigma^{-1} (M - P) \Gamma^{-1}$. Hence $M - P = 2 \Sigma W \Gamma$. Using $W = W^T$ and $P^T = -P$ for a stable feasible point $(W, P)$, we obtain $\Sigma W \Gamma + \Gamma W \Sigma = \frac{1}{2}(M + M^T)$. 

Vectorization implies $(\Gamma\otimes \Sigma + \Sigma \otimes \Gamma) \vc(W) = \frac 1 2 \vc(M+M^T)$. It remains to show that $(\Gamma\otimes \Sigma + \Sigma \otimes \Gamma)$ is invertible. 

Let $K \defeq (\Gamma\otimes \Sigma + \Sigma \otimes \Gamma)$ for shorthand. Notice $K = (\Gamma^{1/2} \otimes \Gamma^{1/2}) (I \otimes \Gamma^{-1/2} \Sigma \Gamma^{-1/2}  + \Gamma^{-1/2}  \Sigma \Gamma^{-1/2} \otimes I) (\Gamma^{1/2} \otimes \Gamma^{1/2})$. Let $K^\prime = (I \otimes \Gamma^{-1/2} \Sigma \Gamma^{-1/2}  + \Gamma^{-1/2}  \Sigma \Gamma^{-1/2} \otimes I)$. Since $(\Gamma^{1/2} \otimes \Gamma^{1/2})$ is invertible it suffices to show $K^\prime$ is invertible. 

Properties of Kronecker products imply that if a matrix $A \in \RR^{n \times n}$ has strictly positive eigenvalues then $\sigma(I \otimes A + A \otimes I) = \{\lambda + \mu: \lambda, \mu \in \sigma(A)\}$ counting mutiplicities \citep{horn-johnson-topics-2008}. Let $\bm{v} \neq \bm{0}$. Then, since $\Sigma \succ 0$ and $\Gamma^{-1/2} \succ 0$ we obtain $\bm{v}^T \Gamma^{-1/2} \Sigma \Gamma^{-1/2} \bm{v} = (\Gamma^{-1/2} \bm{v})^T \Sigma (\Gamma^{-1/2} \bm{v}) > 0$. Hence $\Gamma^{-1/2} \Sigma \Gamma^{-1/2} \succ 0$, so $K^\prime$ is invertible and hence $K$ is invertible. This proves the first identity. 

Next, we prove the second identity. Properties of Kronecker products imply that $(K^\prime)^{-1}$ has eigendecomposition $(K^\prime)^{-1} = \sum_{i=1}^n \sum_{j=1}^n \frac{1}{\lambda_i + \lambda_j} (\bm{v}_i \otimes \bm{v}_j) (\bm{v}_i \otimes \bm{v}_j)^T$. 

Therefore, since $(\Gamma^{1/2} \otimes \Gamma^{1/2})^{-1} = (\Gamma^{-1/2} \otimes \Gamma^{-1/2})$ we obtain: 
\begin{align*}
\vc(W) &= (\Gamma^{-1/2} \otimes \Gamma^{-1/2})
\sum_{i=1}^n \sum_{j=1}^n \frac{1}{\lambda_i + \lambda_j} (\bm{v}_i \otimes \bm{v}_j) (\bm{v}_i \otimes \bm{v}_j)^T
(\Gamma^{-1/2} \otimes \Gamma^{-1/2}) \vc\big(\frac{M + M^T}{2}\big) \\
&= (\Gamma^{-1/2} \otimes \Gamma^{-1/2}) \sum_{i=1}^n \sum_{j=1}^n \frac{1}{2(\lambda_i + \lambda_j)} (\bm{v}_i \bm{v}_i^T \otimes \bm{v}_j \bm{v}_j^T) 
\vc\big(\Gamma^{-1/2} (M + M^T) \Gamma^{-1/2} \big) \\
&= (\Gamma^{-1/2} \otimes \Gamma^{-1/2}) 
\vc\left( \sum_{i=1}^n \sum_{j=1}^n \frac{\bmv_i^T \Gamma^{-1/2}(M+M^T)\Gamma^{-1/2} \bmv_j}{2 (\lambda_i + \lambda_j)} \bmv_i \bmv_j^T \right) \\
W &= \Gamma^{-1/2} \left( \sum_{i=1}^n \sum_{j=1}^n \frac{\bmv_i^T \Gamma^{-1/2}(M+M^T)\Gamma^{-1/2} \bmv_j}{2 (\lambda_i + \lambda_j)} \bmv_i \bmv_j^T \right) \Gamma^{-1/2} 
\end{align*} 

Finally, the formulas for $\vc(P)$ and $P$ follow from similar reasoning, using $W = W^T$ and $W = \frac{1}{2} \Sigma^{-1} (M - P) \Gamma^{-1}$.
%Finally, the third identity follows from similar reasoning to the first identity, using $W = W^T$ and $W = \frac{1}{2} \Sigma^{-1} (M - P) \Gamma^{-1}$. 
%The fourth identity follows from the third and similar reasoning to the proof of the second.
\Halmos
\endproof

\subsection{Proof of Theorem \ref{thm:domination}}

\begin{repeattheorem}[Restatement of Theorem \ref{thm:domination}].
Let $(W^*, P^*)$ be a stable feasible point. Then there is no feasible $(W, P)$ such that $(W, P) \succ (W^*, P*)$. 
\end{repeattheorem}

\proof{Proof of Theorem \ref{thm:domination}.}
{\em Case 1: $P = P^*$.} First, consider a feasible $(W, P)$ such that $P = P^*$. Then $W \neq W^*$. Since $W^*$ is stable, by definition each agent optimizes contracts with respect to $P^*$, so no agent is worse off under $(W^*, P^*)$ then $(W, P^*)$. Hence $(W, P) \not\succ (W^*, P^*)$. 

{\em Case 2: $P \neq P^*$.} Second, suppose that $P \neq P^*$. Let $\Delta_i \defeq g_i(W, P) - g_i(W, P^*)$. It follows that $\Delta_i = (W \bm{e}_i)^T ((P^* - P) \bm{e}_i)$. Let $A \in \RR^{n \times n}$ be defined as $A_{ij} = W_{ij} (P_{ij}^* - P_{ij})$. Then $\Delta_i = \bm{e}_i^T A \bm{1}$. 

Next, notice that $A_{ji} = - A_{ij}$. Therefore, $\sum_i \Delta_i = \bm{1}^T A \bm{1} = 0$. Hence, either $\Delta_i = 0$ for all $i$, or there exists $k$ such that $\Delta_k < 0$. 

{\em Case 2(i)}. Suppose there exists $k$ such that $\Delta_k < 0$. Then $g_k(W, P) < g_k(W, P^*)$. By case $1$, we have $g_k(W, P^*) \leq g_k(W^*, P^*)$. Therefore agent $k$ is strictly worse off, so $(W, P) \not\succ (W^*, P^*)$. 

{\em Case 2(ii)}. Suppose $\Delta_i = 0$ for all $i$. Then $g_i(W, P) = g_i(W, P^*)$ for all $i$. By case $1$, we have $g_i(W, P^*) \leq g_i(W^*, P^*)$. Therefore no agent is better off, so $(W, P) \not\succ (W^*, P^*)$. \Halmos
\endproof

\subsection{Proof of Theorem \ref{thm:pairwise-nash}}
\begin{repeattheorem}[Restatement of Theorem \ref{thm:pairwise-nash}]
Any stable point $(W, P)$ is Higher-Order Nash Stable. 
\end{repeattheorem}

\proof{Proof of Theorem \ref{thm:pairwise-nash}.}
First, we argue $(W, P)$ is a Nash equilibrium. Suppose that agent $i$ wants to shift some of their contracts at the stable feasible point $(W, P)$. Suppose they propose $(w_{i, j_1}^\prime, p_{i, j_1}^\prime), \dots, (w_{i, j_m}^\prime, p_{i, j_m}^\prime)$ for $j_1, \dots, j_m \in [n]$. Let $(W^\prime, P^\prime)$ denote the new feasible point that occurs if all changes are accepted. By Theorem \ref{thm:domination} we know that $(W^\prime, P^\prime) \not\succ (W, P)$, so at least one agent does not prefer $(W^\prime, P^\prime)$. Since the only changes are to edges $\{i, j_1\}, \dots, \{i, j_m\}$, there must exist a $j \in \{j_1, \dots, j_m\}$ who does not prefer $(W^\prime, P^\prime)$. Therefore, they will reject the proposal of agent $i$ to shift to $(w_{ij}^\prime, p_{ij}^\prime)$. 

Then, agent $i$ can choose to either maintain the existing contract $(w_{ij}, p_{ij})$ or delete the edge $\{i, j\}$. We claim that agent $i$ prefers to keep the edge, since they could have chosen to set $w_{ij} = 0$ during the network formation process, no matter what price was offered. But $w_{ij} \neq 0$ at equilibrium $(W, P)$. By stability of $(W, P)$ we know $w_{ij}$ is the optimal choice for agent $i$ at prices $P$. Therefore, after agent $j$ rejects $(w_{ij}^\prime, p_{ij}^\prime)$, it follows that the edge remains at $(w_{ij}, p_{ij})$. 

Since $(W^\prime, P^\prime)$ was arbitrary, we conclude that at equilibrium, agent $i$ cannot propose any set of changes that result in a strictly better network for them. Therefore, their optimal action at $(W, P)$ is to not deviate from the equilibrium. 

Next, we show cartel robustness. Suppose $S \subset [n]$ is a strict subset and $(W^\prime, P^\prime) \neq (W, P)$ is a feasible point differing only at indices $\{i, j\}$ such that $i, j \in S$. By Theorem \ref{thm:domination}, we know $(W^\prime, P^\prime)$ cannot dominate $(W, P)$, so there is some agent $i \in [n]$ that does not prefer $(W^\prime, P^\prime)$ to $(W, P)$. Since $(W^\prime, P^\prime)$ only changes contracts where both members are in $S$, the utility of agents in $[n]\setminus S$ must be unchanged. Therefore $i \in S$, and hence not all members of the cartel have higher utility under $(W^\prime, P^\prime)$. 
\Halmos
\endproof

\subsection{Price Update Rule for Pairwise Negotiations}\label{appendix:price-update-rule}

% We prove a generalization of Proposition \ref{price-update-rule} to the prohibited edges setting. 

% \begin{proposition}[Price after Pairwise Negotiation]\label{price-update-rule}
% Consider a network setting $(\bmu_i, \gamma_i, \Sigma_i, \Psi_i)_{i \in [n]}$. Let $Q_i$ be as in Theorem \ref{thm:stable}. 
% Given a price matrix $P=-P^T$ and a pair of firms $(i,j)$ that are permitted to trade, let $P^\prime$ be another skew-symmetric price matrix such that (a) $P^\prime$ differs from $P$ only in the cells $(i,j)$ and $(j,i)$, (b) $i$ and $j$ both maximize their utility at the same contract size under $P^\prime$, and (c) $i$ and $j$ can choose their optimal contract sizes with all other agents given these prices. 
% Then,
% \begin{align*}
% P^\prime_{ij} &= \frac{1}{Q_{i;j,j} + Q_{j;i,i}}
% \Big(\bm{e}_i^T Q_j (M - P) \bm{e}_j
% - \bm{e}_j^T Q_i (M - P) \bm{e}_i
% \Big)
% + P_{ij}
% \end{align*}
% \end{proposition}

We give an explicit formula for the updated price of a unit contract after a pairwise negotiation. 

\begin{proposition}[Price after Pairwise Negotiation]\label{price-update-rule}
Consider a network setting $(\bmu_i, \gamma_i, \Sigma_i, \Psi_i)_{i \in [n]}$. Let $Q_i$ be as in Theorem \ref{thm:stable}. 
Given a price matrix $P=-P^T$ and a pair of firms $(i,j)$ that are permitted to trade, let $P^\prime$ be another skew-symmetric price matrix such that (a) $P^\prime$ differs from $P$ only in the cells $(i,j)$ and $(j,i)$, (b) $i$ and $j$ both maximize their utility at the same contract size under $P^\prime$, and (c) $i$ and $j$ can choose their optimal contract sizes with all other agents given these prices. 
Then,
\begin{align*}
P^\prime_{ij} &= \frac{1}{Q_{i;j,j} + Q_{j;i,i}}
\Big(\bm{e}_i^T Q_j (M - P) \bm{e}_j
- \bm{e}_j^T Q_i (M - P) \bm{e}_i
\Big)
+ P_{ij}
\end{align*}
\end{proposition}

\proof{Proof.}
Let $A_i \defeq \gamma_i Q_i$ for $i \in [n]$. Since $\Sigma_i \succ 0$ and $\Psi_i \Sigma_i \Psi_i^T$ is a principal submatrix, we know $\Psi_i \Sigma_i \Psi_i^T$ is real symmetric and positive definite, and hence its inverse is as well. Therefore $A_i$ is real symmetric and PSD. (It is not full rank in general, unless $\Psi_i = I$).
Since $\{i, j\}$ is a permitted edge, $\Psi_i \bm{e}_j \neq \bm{0}$ and $\Psi_j \bm{e}_i \neq \bm{0}$.
Therefore $A_{i;j,j} = \bm{e}_j^T A_i \bm{e}_j = (\Psi_i \bm{e}_j)^T (2 \Psi_i \Sigma_i \Psi_i^T)^{-1} (\Psi_i \bm{e}_j) > 0$ since $(2 \Psi_i \Sigma_i \Psi_i^T)^{-1}$ is positive definite. So, $A_{i;j,j} > 0$ and similarly $A_{j;i,i} > 0$.

Now, the optimal contracts for agent $i$ under prices $P^\prime$ are given by $\bm{w}_i = A_i (M - P^\prime) \Gamma^{-1} \bm{e}_i$.
Note that $P^\prime = P + (P^\prime_{ij}-P_{ij})(\bme_i\bme_j^T - \bme_j\bme_i^T)$.
Since both $i$ and $j$ maximize their utility at the same contract size, we have:
%Given some existing prices $P$, we calculate $p_{ij}^\prime$ as: 
\begin{align*}
\bm{w}_{i;j} &= \bm{w}_{j;i} \\
\Rightarrow \bm{e}_j^T \bm{w}_{i} &= \bm{e}_i^T \bm{w}_{j} \\
\Rightarrow \bm{e}_j^T (A_i (M - P^\prime) \Gamma^{-1}) \bm{e}_i &=
\bm{e}_i^T (A_j (M - P^\prime) \Gamma^{-1}) \bm{e}_j \\
\Rightarrow \gamma_j \bm{e}_j^T A_i M \bm{e}_i 
- \gamma_i \bm{e}_i^T A_j M \bm{e}_j 
&= \gamma_j \bm{e}_j^T A_i P^\prime \bm{e}_i 
- \gamma_i \bm{e}_i^T A_j P^\prime \bm{e}_j \\
&= \gamma_j \bm{e}_j^T A_i P \bm{e}_i 
- \gamma_i \bm{e}_i^T A_j P \bm{e}_j 
- (P^\prime_{ij}-P_{ij})\left( \gamma_j \bme_j^T A_i \bme_j + \gamma_i \bme_i^T A_j \bme_i  \right)\\
\Rightarrow P^\prime_{ij} - P_{ij} &= 
\frac{1}{\gamma_j A_{i;j,j} + \gamma_i A_{j;i,i}}  
\Big(\bm{e}_i^T \Gamma A_j (M - P) \bm{e}_j
- \bm{e}_j^T \Gamma A_i (M - P) \bm{e}_i
\Big) \\
&= \frac{1}{Q_{i;j,j} + Q_{j;i,i}}
\Big(\bm{e}_i^T Q_j (M - P) \bm{e}_j
- \bm{e}_j^T Q_i (M - P) \bm{e}_i
\Big)
\end{align*}
\Halmos
\endproof

\subsection{Proof of Theorem \ref{thm:asymp}\label{dynamics-with-masking-appendix}}

First, we characterize pairwise negotiation dynamics as linear in the price updates. 

\begin{theorem}
Consider a network setting $(\bmu_i, \gamma_i, \Sigma_i, \Psi_i)_{i \in [n]}$. 
Define $Q_i$ as in Theorem~\ref{thm:stable}.
Let $s_{ij}=1$ if $\{i, j\}$ is a permitted edge and $0$ otherwise.
Let $L, R \in \RR^{n^2 \times n^2}$ be diagonal matrices such that $L_{(i-1)n + j, (i-1)n + j} = Q_{i;jj} +  Q_{j;ii}$ and $R_{(i-1)n+j, (i-1)n+j} = s_{ij}$, and $L^\dagger$ be the pseudoinverse of $L$.
Let $\Delta_{(t + 1)} = P(t + 1) - P(t)$, where $P(t)$ is the price matrix at time step $t$ of pairwise negotiations. 
Then,
\begin{align*}
\vc(\Delta_{(t+1)})
&=
R  \Big(I_{n^2} -\eta  
L^\dagger K
\Big)
\vc(\Delta_{(t)}),
& \text{where } K &=
\sum\limits_{r=1}^{n} 
\big(
\bme_r \bme_r^T \otimes Q_r
+ Q_r \otimes \bme_r \bme_r^T
\big).
\end{align*}
\label{thm:linearization-dynamics-general}    
\end{theorem}

\proof{Proof.}
Let $\{i, j\}$ be a permitted edge. From Proposition~\ref{price-update-rule}, we obtain: 
\begin{align*}
(\Delta_{(t + 1)})_{ij} 
&= \frac{\eta}{Q_{i;j,j} + Q_{j;i,i}}
\Big(\bm{e}_i^T Q_j (M - P(t)) \bm{e}_j
- \bm{e}_j^T Q_i (M - P(t)) \bm{e}_i
\Big) \\
\rarr (\Delta_{(t + 1)})_{ij}  - (\Delta_{(t)})_{ij} 
&= \frac{\eta}{Q_{i;j,j} + Q_{j;i,i}}
\Big(\bm{e}_i^T Q_j (-\Delta_{(t)}) \bm{e}_j
- \bm{e}_j^T Q_i (-\Delta_{(t)}) \bm{e}_i
\Big) \\
&= \frac{- \eta}{Q_{i;j,j} + Q_{j;i,i}}
\Big(\bm{e}_i^T Q_j \Delta_{(t)} \bm{e}_j
- \bm{e}_j^T Q_i \Delta_{(t)} \bm{e}_i
\Big) \\
&= \frac{- \eta}{Q_{i;j,j} + Q_{j;i,i}}
\bm{e}_i^T 
\Big(Q_j \Delta_{(t)}
- (Q_i \Delta_{(t)})^T
\Big) 
\bm{e}_j \\
\rarr (Q_{i;j,j} + Q_{j;i,i})  \Big((\Delta_{(t + 1)})_{ij}  - (\Delta_{(t)})_{ij}\Big) 
&= 
-\eta  
s_{ij} \bm{e}_i^T
\Big(Q_j \Delta_{(t)}
+ \Delta_{(t)} Q_i
\Big)
\bm{e}_j,
\end{align*}
We assumed that $\{i, j\}$ was a permitted edge above, but notice the identity is also true for prohibited $\{i, j\}$ since both the numerator and denominator become $0$, and we can define their ratio to be $0$.
%since $s_{ij} = Q_{i;j,j} = Q_{j;i,i} = 0$, so both sides become zero. 
Defining $Y_{ij} = \bm{e}_i^T \Big(Q_j \Delta_{(t)} + \Delta_{(t)} Q_i \Big) \bm{e}_j$, and recalling the definitions of $L$ and $R$ from the theorem statement, the above formula becomes
\begin{align}
L  \vc(\Delta_{(t+1)} - \Delta_{(t)})
&= -\eta R \vc(Y).
\label{eq:LRtmp1}
\end{align}

We show next that $\vc(Y) = K \vc(\Delta_{(t)})$, where $K$ is defined in the theorem statement. 
Let $tr$ denote the trace operator. Then,
\begin{align*}
%(\bm{e}_i \otimes \bm{e}_j)^T U \vc(\Delta_{(t)}) = Y_{ij} \\
(\bme_j^T\otimes \bme_i^T) \vc(Y) = Y_{ij}
&= \bm{e}_i^T \Big(Q_j \Delta_{(t)} + \Delta_{(t)} Q_i \Big) \bm{e}_j \\
&= tr\Big(\bm{e}_i^T Q_j \Delta_{(t)} \bm{e}_j\Big) + tr\Big(\bm{e}_i^T\Delta_{(t)} Q_i \bm{e}_j \Big) \\
&= tr\Big(\bm{e}_j^T \Delta_{(t)}^T Q_j^T \bm{e}_i\Big) + tr\Big(\bm{e}_i^T\Delta_{(t)} Q_i \bm{e}_j \Big) \\
&= tr\Big(\Delta_{(t)}^T Q_j^T \bm{e}_i \bm{e}_j^T 
\Big) 
+ tr\Big(Q_i \bm{e}_j \bm{e}_i^T \Delta_{(t)}
\Big) \\
&= \vc(\Delta_{(t)})^T \vc(Q_j^T \bm{e}_i \bm{e}_j^T + (Q_i \bm{e}_j \bm{e}_i^T)^T ) \\
&= \vc(Q_j \bm{e}_i \bm{e}_j^T + \bm{e}_i \bm{e}_j^T Q_i)^T \vc(\Delta_{(t)}) ,
\end{align*}
where we used $Q_i=Q_i^T$.

Hence we need to show $(\bme_j^T \otimes \bme_i^T) K = \vc(Q_j \bm{e}_i \bm{e}_j^T + \bm{e}_i \bm{e}_j^T Q_i)^T$. Letting $\delta$ denote the Kronecker delta, we obtain:
\begin{align}
(\bme_j^T \otimes \bme_i^T) K
&= (\bme_j^T \otimes \bme_i^T)  (\sum\limits_{r=1}^{n} \bm{e}_r \bm{e}_r^T \otimes Q_r + Q_r \otimes \bm{e}_r \bm{e}_r^T) \nonumber\\
&= \sum\limits_{r=1}^{n} \delta_{jr} (\bme_j^T \otimes \bme_i^T Q_r)
+ \delta_{ir} (\bme_j^T Q_r \otimes \bme_i^T)\nonumber\\
&= (\bme_j^T \otimes \bme_i^T Q_j)
+ (\bme_j^T Q_i \otimes \bme_i^T)\nonumber\\
&= \left( \bme_j\otimes Q_j \bme_i + Q_i\bme_j\otimes \bme_i \right)^T.\label{eq:RK}
\end{align}
Now, we observe that $\bme_j\otimes Q_j\bme_i$ is the vectorization of a matrix whose $j^{th}$ column is $Q_j\bme_i$, i.e., the matrix $Q_j\bme_i\bme_j^T$.
Similarly, $Q_i\bme_j\otimes \bme_i$ is the vectorization of a matrix whose $i^{th}$ row is $(Q_i\bme_j)^T$, i.e., the matrix $\bme_i\bme_j^TQ_i$.
Hence,
$(\bme_j^T \otimes \bme_i^T) K = \vc(Q_j\bme_i\bme_j^T+\bme_i\bme_j^TQ_i)^T$, as desired.

Plugging into Eq.~\eqref{eq:LRtmp1},
\begin{align*}
L  \vc(\Delta_{(t+1)} - \Delta_{(t)})
&= -\eta RK \vc(\Delta_{(t)}) \\
\rarr L  \vc(\Delta_{(t+1)})
&= L  \vc(\Delta_{(t)}) -\eta R K\vc(\Delta_{(t)}) \\
\rarr \vc(\Delta_{(t+1)}) &= 
\Big(L^\dagger L -\eta  L^\dagger R K \Big) \vc(\Delta_{(t)}) \\
\rarr \vc(\Delta_{(t+1)}) &= 
\Big(R -\eta  R L^\dagger K \Big) \vc(\Delta_{(t)}) \\
&= R \Big(I_{n^2} -\eta  L^\dagger K \Big) \vc(\Delta_{(t)}),
\end{align*}
where we used the facts that $(\Delta_{t})_{ij}=(\Delta_{(t+1)})_{ij}=0$ for disallowed edges, and $L^\dagger L = R$ and $LR = RL = L$, which can be easily confirmed by inspection of these diagonal matrices.
\Halmos
\endproof

We use Lyapunov theory to analyze the convergence of pairwise negotiation dynamics. In particular, we need the the discrete Lyapunov equation, also called the Stein equation.
\begin{theorem}[\cite{callier-desoer} 7.d]
\label{thm:callier}
For the discrete-time dynamical system $\bm{x}_{t + 1} = A \bm{x}_{t}$, with $\bm{x}_t \in \RR^n$, the following are equivalent:
\begin{enumerate}
    \item The system is globally asymptotically stable towards $\bm{0}$. 
    \item For any positive definite $R \in \RR^{n \times n}$, there exists a unique solution $X \succ 0$ to the equation 
    $$A X A^T - X = - R$$
    \item For any eigenvalue $\lambda$ of $A$, $\abs{\lambda}<1$.
\end{enumerate}
\end{theorem}

Pairwise negotiation dynamics can be described as a discrete-time linear system in $\vc(\Delta_{t})$, where $\Delta_{t}$ is the price difference at time $t$. Clearly, the system converges iff $\Delta_{t}$ approaches zero. Therefore, we can use the Stein equation to prove global asymptotic stability conditions. 

We will also need the {\em commutation matrix}. 
\begin{lemma}[\cite{horn-johnson-topics-2008}]
Let $\Pi^{(n, n)}: \RR^{n^2} \to \RR^{n^2}$ be a permutation matrix (called the {\em $(n, n)$ commutation matrix}) defined as $\Pi^{(n, n)} = \sum\limits_{i = 1}^{n} \sum\limits_{j = 1}^{n} \bm{e}_i \bm{e}_j^T \otimes \bm{e}_j \bm{e}_i^T$. Then for any $A, B \in \RR^{n \times n}$, we have 
\begin{align*}
A \otimes B = \Pi^{(n, n)} (B \otimes A) (\Pi^{(n, n)})^T
\end{align*}
\label{lemma:commutation-matrix}
\end{lemma} 

Recall that for a linear operator $T$ that $\sigma(T)$ denotes the eigenvalues of $T$. We are ready to prove Part 1 of Theorem \ref{thm:asymp}. 

\begin{proposition}[Part 1 of Theorem \ref{thm:asymp}]\label{real-part-eigs-condition-stability-prop}
%Consider the pairwise negotation dynamics in Theorem \ref{thm:asymp}. 
Let $L, R, K$ be defined as in Theorem \ref{thm:linearization-dynamics-general}. 
For a matrix $X \in \RR^{n^2 \times n^2}$ let $X\mid_R$ denote the principal submatrix of $X$ corresponding to the nonzero rows/columns of $R$. 
Define $\eta^\star = \min\limits_{\lambda \in \sigma((L^\dagger K)\mid_R)} \frac{2}{\lambda}$. 
%\blue{Define $\eta^\star = \min\{1, 2 / \left\|(L^\dagger K)\mid_R\right\|\}$}. 
Then, for any $\eta \in (0, \eta^\star)$, $\vc(\Delta_{(t)})$ is globally asymptotically stable towards $\bm{0}$. 
\end{proposition}

\proof{Proof of Proposition \ref{real-part-eigs-condition-stability-prop}.}
Let $T = R(I - \eta L^\dagger K)$. 
By Theorem~\ref{thm:callier}, the dynamics are globally asymptotically stable towards $\bm{0}$ iff for all $\lambda \in \sigma(T)$, we have $\abs{\lambda} < 1$.

From Eq.~\eqref{eq:RK} for a prohibited edge $(i,j)$, we see that $(\bme_j^T\otimes \bme_i^T)K=\bm{0}^T$, since $Q_i\bme_j=\bm{0}=Q_j\bme_i$.
Hence, $K=RK$.
Taking transposes and noting that both $K$ and $R$ are symmetric, we find $KR=K$.
Hence, $T=R(I-\eta L^\dagger K)=R(I-\eta L^\dagger K)R$, where we used $R^2=R$.
Thus, $T$ is zero except for the principal submatrix corresponding to the nonzero columns of $R$.
% rows and columns 
So, to apply Theorem~\ref{thm:callier}, we only require $\abs{\lambda}<1$ for $\lambda \in \sigma(T\mid_R)$.

For clarity of exposition we will first consider the case where $R = I$ (no prohibited edges).
Then, the eigenvalues of $T\mid_R=T$ equal $1-\eta \lambda$, where $\lambda\in\sigma(L^{-1}K)=\sigma(L^{-1/2}KL^{-1/2})$ by a similarity transformation.
Also, $K=U_1+U_2$, where $U_1=\sum\limits_{r =1}^{n} \big(\bm{e}_r \bm{e}_r^T \otimes Q_r\big)$ and $U_2 \defeq \sum\limits_{r =1}^{n} \big(Q_r \otimes \bm{e}_r \bm{e}_r^T \big)$.
The matrix $U_1$ is block diagonal with positive-definite blocks $Q_r\succ 0$, so $U_1\succ 0$.
By Lemma~\ref{lemma:commutation-matrix}, $U_2$ is similar to $U_1$ via a permutation matrix, so $U_2\succ 0$.
Hence, $K\succ 0$, and $L^{-1/2}KL^{-1/2}\succ 0$.
So, the eigenvalues of $L^{-1}K$ are real and positive.
Hence, we have convergence iff for all $\lambda\in\sigma(L^{-1}K)$, we have
$1 > (1-\eta\lambda)^2 = 1 - 2\eta\lambda + \eta^2\lambda^2$. 
 i.e., $\lambda < 2/\eta$.
Hence, $\eta^\star=2/\|L^{-1}K\|$ as required.

Now we consider the prohibited edges setting ($R\neq I$).
Here, convergence occurs iff $|1-\eta\lambda|<1$ for all $\lambda\in\sigma((L^\dagger K)\mid_R)$.
Since $RL^\dagger R=L^\dagger$ and $RKR=K$, we have $(L^\dagger K)\mid_R = L^\dagger\mid_R K\mid_R = (L\mid_R)^{-1} K\mid_R$.
Arguing as above, it suffices to show that $K\mid_R \succ 0$.
We claim $K\mid_R = V_1 + V_2$ where $V_1$ is a block diagonal matrix with $i^{th}$ block equal to $(2 \gamma_i \Psi_i \Sigma_i \Psi_i^T)^{-1} \succ 0$, and $V_2$ is similar to $V_1$ via Lemma~\ref{lemma:commutation-matrix}.
Hence $K\mid_R \succ 0$ and the expression for $\eta^\star$ follows. 
\Halmos
\endproof

%Next, we prove part $2$ of Theorem \ref{thm:asymp}.
%\proof{Proof of part $2$ of Theorem \ref{thm:asymp}.}

% in the theorem we define alpha = 1 - eta lambda_min, 1 - eta lambda_max

%, &\text{where }\alpha = \max\{ \abs{1-\eta\lambda_{\textrm{min}}}, \abs{1-\eta\lambda_{\textrm{max}}} \}

\begin{proposition}[Part 2 of Theorem \ref{thm:asymp}]\label{prop:expConvergence}
%Consider the pairwise negotation dynamics in Theorem \ref{thm:asymp}. 
We define $\eta^\star$ as in Proposition~\ref{real-part-eigs-condition-stability-prop}, and $L, R, K, \alpha$ as in Theorem \ref{thm:linearization-dynamics-general}. 
Let $\eta\in(0,\eta^\star)$.
Then, 
\begin{align*}
\|P(t)-P^\star\|_F &\leq \frac{\alpha^t}{1-\alpha} \cdot \|P(1)-P(0)\|_F
\end{align*}
Here, $P^\star$ is the stable point to which the negotiation converges.
\end{proposition}
\proof{Proof.}
Let $\beta$ denote the greatest eigenvalue in absolute value of $R(I_{n^2} - \eta L^\dagger K)$. From Theorem~\ref{thm:linearization-dynamics-general}, we have $\|\Delta_{t+1}\|_F \leq \abs{\beta} \|\Delta_t\|_F$. Recall that $\lambda_{\textrm{max}}, \lambda_{\textrm{min}}$ denote largest and smallest eigenvalues of the matrix $(L^\dagger K)\mid_R$ respectively. Since $\norm R \norm = 1$, it follows that $\abs{\beta} = \max\{ \abs{1-\eta\lambda_{\textrm{min}}}, \abs{1-\eta\lambda_{\textrm{max}}} \} = \alpha$. 

Then, 
\begin{align*}
\|P^\star-P(t)\|_F
  &\leq \sum_{i>t} \|\Delta_i\|_F \\
  &\leq \|\Delta_t\|_F(\alpha + \alpha^2 + \ldots) \\
  &\leq \|\Delta_t\|_F \frac{\alpha}{1-\alpha}\\
  &\leq (\alpha^{t-1} \|\Delta_1\|_F)\frac{\alpha}{1-\alpha}  \\
  &= \|\Delta_1\|_F \frac{\alpha^t}{1-\alpha}
\end{align*}
Since $\|\Delta_1\|_F = \| P(1) - P(0) \|_F$ we are done.
\Halmos
\endproof

\subsection{Example of Convergence Conditions and Rate}\label{appendix:example-asymp}

The following example illustrates Theorem \ref{thm:asymp} in the setting of Example \ref{example:thrm1} (Appendix~\ref{appendix:example-thrm1}). 

\begin{example}[Convergence Conditions and Rate]
In the setting of Example \ref{example:thrm1}, we have 
\[
Q_1 = \begin{bmatrix} 0 & 0 & 0 \\ 0 & 2/3 & -1/3 \\ 0 & -1/3 & 2/3 \end{bmatrix}, 
Q_2 = \begin{bmatrix} 1 & 0 & 0 \\ 0 & 0 & 0 \\ 0 & 0 & 0\end{bmatrix},
Q_3 = \begin{bmatrix} 2 & 0 & 0 \\ 0 & 0 & 0 \\ 0 & 0 & 0 \end{bmatrix},
\]

Hence 
\[
K = \begin{bmatrix}
0 & 0 & 0 & 0 & 0 & 0 & 0 & 0 & 0 \\
0 & \frac{2}{3} + 1 & \frac{-1}{3} & 0 & 0 & 0 & 0 & 0 & 0 \\
0 & \frac{-1}{3} & \frac{2}{3} + 2 & 0 & 0 & 0 & 0 & 0 & 0 \\
0 & 0 & 0 & 1 + \frac{2}{3} & 0 & 0 & \frac{-1}{3} & 0 & 0 \\
0 & 0 & 0 & 0 & 0 & 0 & 0 & 0 & 0 \\
0 & 0 & 0 & 0 & 0 & 0 & 0 & 0 & 0 \\
0 & 0 & 0 & \frac{-1}{3} & 0 & 0 & 2 + \frac{2}{3} & 0 & 0 \\
0 & 0 & 0 & 0 & 0 & 0 & 0 & 0 & 0 \\
0 & 0 & 0 & 0 & 0 & 0 & 0 & 0 & 0 \\
\end{bmatrix}
\]

Also, $L$ is the diagonal matrix with $L_{i, i} = K_{i, i}$ for $i \in [9]$. 
Since the permitted edges are $\{(1, 2), (1, 3)\}$, $R = \{2, 3\}$ and so $(L^\dagger K)_{R} = \begin{bmatrix} 1 & \frac{-1}{5} \\ \frac{-1}{8} & 1 \end{bmatrix}$.
Hence $\lambda_{min} = 1 - \frac{1}{2\sqrt{10}}, \lambda_{max} = 1 + \frac{1}{2\sqrt{10}}$, and $\eta^* = \frac{2}{1 + (40)^{-1/2}} \approx 1.727$. 

It follows that pairwise negotiations with $\eta \in (0, \frac{2}{1 + (40)^{-1/2}})$ are globally asymptotically stable. Suppose that $\eta = 0.99$. Then $\alpha = (1 - \eta \cdot (1 - \frac{1}{2\sqrt{10}})) \approx 0.17$. Hence after $t$ rounds, the distance of $P(t)$ to $P^*$ shrinks by a factor of $\approx \frac{0.17^t}{0.83}$. 
\end{example}

\subsection{Proof of Theorem \ref{thm:asymp:random}}

% \textcolor{blue}{Need $\norm \Gamma^{-1} \norm = O(1)$ in the theorem statement.}

We will use a series of Lemmas to reduce the result of Theorem \ref{thm:asymp:random} to a matrix concentration inequality in each of the $\hat{\Sigma}_i$. 

\begin{lemma}
Let $\widehat{\eta^*}, \eta^*$ be as in Theorem \ref{thm:asymp:random}. Suppose all edges are permitted. 

Suppose that for all $i \in [n]$, we have $\|\hat{\Sigma}_i^{-1}-\Sigma^{-1}\|=o(1)$. Then, $\hat{\eta}^\star > \eta^\star(1-o(1))$.
\label{lemma:sigma-inv-sufficient}
\end{lemma}
\proof{Proof.}
Let $\hat{L}, \hat{K} \in \RR^{n^2 \times n^2}$ be as in Theorem \ref{thm:asymp:random}, but built using $\hat{\Sigma}_1, \dots, \hat{\Sigma}_n$ instead of $\Sigma, \dots, \Sigma$. Let $L, K$ be defined similarly to $\hat{L}, \hat{K}$ but using $\Sigma$ in place of all $\hat{\Sigma}_i$. 

Then $\widehat{\eta^*} \defeq \frac{2}{\max\sigma(\hat{L}^{-1} \hat{K})}$ and $\eta^* \defeq  \frac{2}{\max\sigma(L^{-1} K)}$. 

Let $\eps_L, \eps_K \in \RR^{n^2 \times n^2}$ be such that $\hat{L}^{-1} = L^{-1} + \eps_L$ and $\hat{K} = K + \eps_K$. We will bound $\norm \eps_L \norm, \norm \eps_K \norm$. 

Let $Q_i, \hat{Q}_i$ be defined as in Theorem \ref{thm:stable}, so $Q_i \defeq (2 \gamma_i \Sigma)^{-1}$ and $\hat{Q}_i \defeq (2 \gamma_i \hat{\Sigma}_i)^{-1}$. Let $\alpha = \max\limits_{i \in [n]} \norm \hat{Q}_i - Q_i \norm$. Notice $\norm \Gamma^{-1} \norm = O(1)$, so $\alpha = o(1)$. 

First, since $L$ is diagonal, $\norm \eps_L \norm \leq \max\limits_{i, j \in [n]} \big((\hat{Q}_{i;jj} - Q_{i;jj}) + (\hat{Q}_{j;ii} - Q_{j;ii}) \big) \leq 2 \max\limits_{i, j \in [n]} \big(\hat{Q}_{i;jj} - Q_{i;jj} \big) \leq 2 \max\limits_{i \in [n]} \norm \hat{Q}_i - Q_i \norm = 2a$. 

Second, let $\hat K \defeq \hat U_1 + \hat U_2$ where $\hat U_1, \hat U_2$ are defined analogously to $U_1, U_2$ in the proof of Theorem \ref{thm:asymp}. Letting $\Pi$ be the $(n, n)$ commutation matrix of Lemma \ref{lemma:commutation-matrix}, we know $\hat U_2 = \Pi \hat U_1 \Pi^T$, so $\norm \eps_K \norm \leq 2 \norm \hat U_1 - U_1 \norm$. Since $U_1, \hat U_1$ are block diagonal with $i^{th}$ blocks $Q_i, \hat Q_i$ respectively, it follows $\norm \hat U_1 - U_1 \norm = \max\limits_{i \in [n]} \norm \hat{Q}_i - Q_i \norm = \alpha$. Hence $\norm \eps_K \norm \leq 2 \alpha$. 

Third, notice that since $\norm \Sigma \norm$ and $\norm \Gamma \norm$ are assumed to be $O(1)$ that $\norm L^{-1} \norm = O(\max_{i} \norm Q_i \norm) = O(1)$ and $\norm K \norm = O(\max_{i} \norm Q_i \norm) = O(1)$. 
So, 
\begin{align*}
\norm \hat{L}^{-1} \hat{K} - L^{-1} K \norm_2 
&\leq \norm \eps_L \norm \norm K \norm + \norm L^{-1} \norm \norm \eps_K \norm + \norm \eps_L \norm \norm \eps_K \norm \\
&\leq 2\alpha (\norm K \norm + 2 \alpha)
+ 4 \alpha (\norm L^{-1} \norm + \alpha) \\
&= 4 \alpha (\norm K \norm + \norm L^{-1} \norm) + 8 \alpha^2 \\
&\leq o(1)
\end{align*}
We conclude that $\norm \hat{L}^{-1} \hat{K} \norm_2 \leq \norm L^{-1} K \norm_2 + o(1)$, so $\widehat{\eta^*} \geq \frac{\eta^*}{1 + (o(1) / \norm L^{-1} K \norm)} \geq (1 - o(1)) \eta^*$.
\Halmos
\endproof

\begin{lemma}
Suppose for $i \in [n]$, we have $\delta_i \defeq\| \hat{\Sigma}_i - \Sigma\|=o(1)$. Then  $\|\hat{\Sigma}_i^{-1} - \Sigma_i^{-1}\|=o(1)$.
\label{lemma:sigma-inv-reduction}
\end{lemma}
\proof{Proof.}
Weyl's inequality implies that $\lambda_{min}(\hat{\Sigma}_i) \geq \lambda_{min}(\Sigma) - \norm \hat{\Sigma}_i - \Sigma\norm$. Therefore, 
\begin{align*}
\|\hat{\Sigma}_i^{-1}\| &= \frac{1}{\lambda_{min}(\hat{\Sigma}_i)}\\
&\leq \frac{1}{\lambda_{min}(\Sigma)- \delta_i} \\
&=\frac{1}{\lambda_{min}(\Sigma)}
\Big(1+\frac{\delta_i}{\lambda_{min}(\Sigma)} + 
O\Big( \Big(\frac{\delta_i}{\lambda_{min}(\Sigma)}\Big)^2 \Big)
\Big)\\
&=\|\Sigma^{-1}\|(1+o(1))\\
\rarr \|\hat{\Sigma}_i^{-1} - \Sigma^{-1}\| &= \|\Sigma^{-1}(\Sigma_i-\hat{\Sigma}_i)\hat{\Sigma}_i^{-1}\|\\
&\leq (1 + o(1)) \|\Sigma^{-1}\|^2 \delta_i \\
&\leq o(1)
\end{align*}
The last step follows from the fact $\norm \Sigma^{-1} \norm = O(1)$.
\Halmos
\endproof

The hypothesis of Lemma \ref{lemma:sigma-inv-reduction} follows from a standard argument on the concentration of random covariance matrices. 
\begin{theorem}\label{thrm:wishart-vershynin2}
Under the setting of Theorem~\ref{thm:asymp:random}, with probability at least $1-e^{-\Omega(n)}$, we have $\|\hat{\Sigma}_i - \Sigma\|=o(1)$ for all $i\in[n]$.
\end{theorem}

\proof{Proof of Theorem~\ref{thrm:wishart-vershynin2}.}
Let $\bm{X}_1, \dots, \bm{X}_m \iid \mathcal{N}(\bm{0}, \Sigma)$ be the samples.
Let $\hat{\bm{\mu}} = \frac 1 m \sum\limits_{i=1}^{m} \bm{X}_i$, and $\tilde{\Sigma}_i \defeq \frac{1}{m} \sum\limits_{i=1}^{m} \bm{X}_i \bm{X}_i^T$.
%be the biased covariance estimate, and $\hat{\Sigma} \defeq \frac{1}{m-1} \sum\limits_{i=1}^{m} (\bm{X}_i -\hat{\bm{\mu}})(\bm{X}_i - \hat{\bm{\mu}})^T$ be the unbiased covariance estimate. 
Then, $\hat{\Sigma}_i=m/(m-1)\cdot(\tilde{\Sigma}_i-\hat{\bm\mu}\hat{\bm\mu}^T)$.
Hence,
$$\|\hat{\Sigma}_i-\Sigma\| \leq m/(m-1)\cdot\left(\|\tilde{\Sigma}_i-\Sigma\| + \|\hat{\bm\mu}\hat{\bm\mu}^T\|\right) = m/(m-1)\left(\|\tilde{\Sigma}_i-\Sigma\| + \|\hat{\bm\mu}\|^2\right).$$
Now, $\hat{\bm{\mu}} \sim \mathcal{N}(0, \frac 1 m \Sigma)$, so $\sqrt{m} \Sigma^{-1/2} \hat{\bm{\mu}} \sim \mathcal{N}(0, I_n)$.
By \cite{vershynin-book} (4.7.3 and 2.8.3), there exist constants $c, c_2>0$ such that for any $u, \eps>0$,
\begin{align*}
\PP\bigg[\norm \tilde{\Sigma}_i - \Sigma \norm_2 \leq c \norm \Sigma \norm_2  \bigg(\sqrt{\frac{n + u}{m}} + \frac{n + u}{m}\bigg)\bigg] &\geq 1 -2e^{-u},\\
\PP\bigg[\bigg\lvert \frac 1 n \norm \sqrt{m} \Sigma^{-1/2} \hat{\bm{\mu}} \norm_2^2 - 1 \bigg\rvert \leq \eps \bigg] &\geq 1 - 2e^{-c_2 n \min(\eps, \eps^2)}
\end{align*}
Now we set $\eps>1$ and $u=c_3 n$ for some constant $c_3>0$.
Then, when $m=\lceil n\log n\rceil$, we have $(n+u)/m=o(1)$
Then, with probability at least $1-2e^{-c_3 n}-2e^{-c_2\eps n}$, we have
\begin{align*}
& \norm \tilde{\Sigma}_i - \Sigma \norm_2 \leq \norm \Sigma \norm\cdot o(1),\\
\text{and }& \norm \Sigma^{-1/2} \hat{\bm{\mu}} \norm_2^2 \leq \frac{(1+\eps)n}{m} \Rightarrow \norm\hat{\bm\mu}\norm^2 \leq \frac{(1+\eps)n\norm\Sigma\norm}{m}=\norm\Sigma\norm\cdot o(1),\\
\Rightarrow & \norm\hat{\Sigma}_i-\Sigma\| \leq \norm\Sigma\norm\cdot o(1).
\end{align*}
Choosing large enough $c_3$ and $\eps$, this statement holds for all $i\in[n]$ with probability greater than $1-e^{log n - c_4 n}=1-e^{-\Omega(n)}$.
\Halmos
\endproof

Theorem \ref{thm:asymp:random} follows easily.
\proof{Proof of Theorem \ref{thm:asymp:random}}
When all edges are permitted, the proof follows from Theorem~\ref{thrm:wishart-vershynin2}, Lemma~\ref{lemma:sigma-inv-sufficient}, and Lemma~\ref{lemma:sigma-inv-reduction}.

%Suppose $m = \ceil*{n\log n}$. We will first consider the case of all edges permitted; the proof for prohibited edges is almost identical. 
%
%Applying Theorem \ref{thrm:wishart-vershynin} to each $\hat{\Sigma}_i$, a union bound gives $\norm \Sigma - \hat{\Sigma}_i \norm \leq 2 \sqrt{2} c \norm \Sigma \norm \sqrt{\frac{n}{m}}$ for all $i \in [n]$ with probability at least $(1 - 2e^{-n})^n \geq 1 - e^{- n + O(\log n)} \geq e^{-\Omega(n)}$. Since $\norm \Sigma \norm = O(1)$ and $\frac{n}{m} = o(1)$, this simplifies to $\norm \Sigma - \hat{\Sigma}_i \norm  = o(1)$ for all $i$. By Lemmas \ref{lemma:sigma-inv-sufficient} and \ref{lemma:sigma-inv-reduction}, we have $\hat{\eta}^* \geq \eta^* (1 - o(1))$ with probability at least $1 - e^{-\Omega(n)}$. 

If there are prohibited edges, then we must use matrix concentration to bound $\max\sigma(\hat{L}^\dagger \hat K)$ instead of $\max\sigma(\hat{L}^{-1} \hat K)$. Notice that prohibited edges have the effect of simply zeroing out certain rows and columns of $Q_i$, so that $Q_i \defeq \Psi_i (2 \gamma_i \Psi_i^T \Sigma_i \Psi_i)^{-1} \Psi_i^T$, rather than $(2 \gamma_i \Sigma_i)^{-1}$. Therefore, we can use Theorem \ref{thrm:wishart-vershynin2} to bound $\norm \Psi_i^T \hat{\Sigma}_i \Psi_i - \Psi_i^T \Sigma \Psi_i \norm$ for all $i$, and then prove the appropriate analogue of Lemma \ref{lemma:sigma-inv-sufficient}. In particular, the sample size requirement remains the same. 
\Halmos
\endproof

\subsection{Proof of Proposition \ref{prop:sdp-sigma-recovery}}
\label{appendix:inference}

\begin{repeatproposition}[Restatement of Proposition \ref{prop:sdp-sigma-recovery}.]
Finding the maximum likelihood estimator of $\Sigma$ under Assumption~\ref{assume:SDP} is equivalent to the following SDP:
\begin{align*}
\min\limits_{\Sigma} \sum\limits_{t = 1}^{T - 1}
\norm \Sigma (W(t + 1) - W(t)) + 
(W(t + 1) - W(t)) \Sigma \norm_F^2 
\\
\Sigma \succeq 0 \\
\tr(\Sigma) = 1
\end{align*}
\end{repeatproposition}

Recall that in Assumption \ref{assume:SDP} we assumed that $M_{ij}(t)$ varies independently according to a Brownian motion with the same parameters for all $(i,j)$. To avoid ambiguity, we recall the definition of a standard Brownian motion as follows. 

\begin{definition}[Brownian Motion]
For $d \geq 1$, a $d$-dimensional Brownian motion with scale parameter $\sigma > 0$ is a stochastic process $\{\bm{X_t}: t \geq 0\}$ such that $\bm{X_t} \in \RR^d$ for all $t$, the components of $\bm{X_t}$ are independent, and for all $j \in [d]$, 

i) The process $\{(\bm{X_t})_j: t \geq 0\}$ has independent increments. 

ii) For $r > 0$, the increment $(\bm{X_{t+r}})_j - (\bm{X_{t}})_j$ is distributed as $N(0, r \sigma^2)$.  

iii) With probability $1$, the function $t \mapsto \bm{X_t}$ is continuous on $[0, \infty)$. 
\end{definition}

We can derive the SDP of Proposition \ref{prop:sdp-sigma-recovery} as follows. 

\begin{proposition}
\label{small-mean-shift-prop}
Under Assumption~\ref{assume:SDP}, the maximum likelihood estimator for $\Sigma$ is the unique $\Sigma \succ 0$ such that $tr\Sigma = 1$ and
\begin{itemize}
\item {\bf Consistency}: For all $t \in [T]$, 
\begin{align*}
W(t)\Sigma + \Sigma W(t) &= \frac 1 2 (M(t) + M(t)^T) \quad \text{for some $M(1), M(2), \ldots M(t)$} 
\end{align*}
\item {\bf Minimum mean shift}: The resulting $M(1), \dots, M(T)$ minimize the objective 
 \begin{align*}
 \sum\limits_{t=1}^{T-1} \norm M(t+1) - M(t) \norm_F^2
 \end{align*}
\end{itemize}
\end{proposition}

\proof{Proof of Proposition \ref{small-mean-shift-prop}.}

\begin{align*}
P(M(1), \ldots, M(T) &\mid W(1), \ldots, W(T), \Sigma) &\\
&\propto P(W(1), \ldots W(T) \mid M(1), \ldots, M(T), \Sigma)\cdot P(M(1), \ldots, M(T)\mid \Sigma)\\
&= \left(\prod_{t=1}^T \mathds{1}_{W(t)\Sigma + \Sigma W(t) = 0.5(M(t) + M(t)^T)}\right) \left( \prod_{t=1}^{T-1} P(M(t+1)-M(t))\right)\\
&= \left(\prod_{t=1}^T \mathds{1}_{\vc(W(t)) = 0.5(\Sigma\otimes I + I\otimes \Sigma)\vc(M(t) + M(t)^T)}\right) \left( \prod_{t=1}^{T-1} \exp(-\frac{\|\vc(M(t+1)-M(t))\|^2}{2\sigma^2})\right),
\end{align*}
where the first step follows from Bayes' Rule, the second step from Corollary~\ref{cor:stability:sharedSigma}, and the third from Assumption~\ref{assume:SDP}.
The theorem follows from the observation that for any matrix $X$, we have $\|\vc(X)\|^2 = \|X\|_F^2$.
\Halmos
\endproof

The proof of Proposition \ref{prop:sdp-sigma-recovery} follows easily. 

\proof{Proof of Proposition \ref{prop:sdp-sigma-recovery}.}
By Proposition \ref{small-mean-shift-prop}, we obtain the SDP 
\begin{align*}
\min\limits_{\Sigma} \sum\limits_{t=1}^{T-1} \norm M(t+1) - M(t) \norm_F^2
\\
\forall t \in [T]: W(t)\Sigma + \Sigma W(t) &= \frac 1 2 (M(t) + M(t)^T)   
\end{align*}
under the assumptions of $\Sigma \succ 0$ and $\tr(\Sigma) = 1$. % follow from Assumption \ref{assume:SDP}. It remains to analyze the objective. 
Since the Frobenius norm is invariant under transposes, we have
$$\sum\limits_{t=1}^{T-1} \norm M(t+1) - M(t) \norm_F^2 \propto \sum\limits_{t=1}^{T-1} \norm (M(t+1) + M(t+1)^T) - (M(t) + M(t)^T) \norm_F^2.$$
We can replace $M(t) + M(t)^T$ with $2 W(t)\Sigma + 2 \Sigma W(t)$ for all $t \in [T]$ to obtain the equivalent objective function $\sum\limits_{t=1}^{T-1} \norm (W(t + 1) - W(t))\Sigma + \Sigma(W(t + 1) - W(t)) \norm_F^2$ (up to a constant). 
This substitution enforces the fixed point equation $W(t)\Sigma + \Sigma W(t) = \frac 1 2 (M(t) + M(t)^T)$ for all $t \in [T]$, so the conclusion follows.  
\Halmos
\endproof

\begin{remark}[The prohibited edges setting.]
\label{rem:inference_prohibited}
Proposition \ref{prop:sdp-sigma-recovery} generalizes straightforwardly to the setting of prohibited edges. Let $E$ denote the set of permitted edges. Then minimum mean shift assumption is
equivalent to minimizing $\sum\limits_{t=1}^{T - 1} \sum\limits_{\{i, j\} \in E} \big(M(t+1) + M(t+1)^T - M(t) - M(t)^T\big)_{ij}^2$. In words, the objective just zeroes out prohibited edges, since mean estimates for prohibited edges have no effect on the network. 
For a network setting $(\bmu_j, \Sigma, \gamma_j, \Psi_j)_{j \in [n]}$, some algebra gives $M(t)_{ij} = \bm{e}_i^T 2 \gamma_j  (\Psi_j^T \Psi_j) \Sigma (\Psi_j^T \Psi_j) W(t) \bm{e}_j$. Notice $\Psi_j^T \Psi_j \in \RR^n$ is a diagonal matrix with $(\Psi_j^T \Psi_j)_{ii} = 1$ if $\{i, j\} \in E$ and zero otherwise. 
Therefore, it is clear that upon substitution, the objective is an SDP in $\Sigma$ with the same constraints. 
\end{remark}

\subsection{Proof of Theorem \ref{thm:friction-equilibrium-general}}

\begin{repeattheorem}[Restatement of Theorem \ref{thm:friction-equilibrium-general}]
Suppose that for each firm $i\in [n]$, the function $F_i: \RR^n \to \RR$ is twice differentiable, and there exist strictly increasing functions $f_{ji}: \RR \to \RR$ such that  for all $\bm{x} \in \RR^n$, $\nabla F_i(\bm{x}) = [f_{1i}(x_1), \dots, f_{ni}(x_n)]^T$.
Then there exists a unique stable point.
\end{repeattheorem}
\proof{Proof of Theorem \ref{thm:friction-equilibrium-general}.}
Note that the Hessian of $F_i(\bmw_i)$ is a positive diagonal matrix due to the conditions on $F_i(.)$.
So, any stationary point is a local maximum.
Hence, it suffices to show the existence of a unique stationary point.

Let $R(W)$ be an $n\times n$ matrix whose $(i,j)^{th}$ entry $R(W)_{ij} := f_{ij}(W_{ij})$.
If a stable point $(W, P)$ exists, it must satisfy $W=W^T$, $P=-P^T$, and
\begin{align}
M-P = 2\Sigma W \Gamma + R,\label{eq:newgrad}
\end{align}
following the same steps as the proof for Corollary~\ref{cor:stability:sharedSigma}.
Adding this equation to its transpose, the stable point must satisfy
\begin{align*}
(M+M^T)/2 &= (\Sigma W \Gamma + \Gamma W \Sigma) + (R(W)+R(W)^T)/2.
\end{align*}
For a stable point, $[R(W)+R(W)^T]_{ij} = f_{ij}(W_{ij}) + f_{ji}(W_{ji}) = (f_{ij}+f_{ji})(W_{ij})$, using $W=W^T$.
Define $S(W)$ to be an $n\times n$ matrix with $S(W)_{ij} = (1/2)\cdot (f_{ij}+f_{ji})(W_{ij})$.
Hence, the stable point must satisfy
\begin{align}
(M+M^T)/2 &= (\Sigma W \Gamma + \Gamma W \Sigma) + S(W)\label{eq:sandberg}\\
\Leftrightarrow \vc((M+M^T)/2) &= \underbrace{(\Gamma\otimes \Sigma + \Sigma \otimes \Gamma)}_{Q} \vc(W) + \vc(S(W)).\nonumber
\end{align}
Note that $Q$ is positive-definite (from the proof of Corollary~\ref{cor:stability:sharedSigma}), and each entry of $\vc(S(W))$ is a function of the corresponding entry of $\vc(W)$.
By Theorems 1 and 2 of \citet{sandberg1972existence}, Eq.~\eqref{eq:sandberg} has a unique solution if (1) for all diagonal $D \succ 0$, $\det(D + Q) > 0$ and (2) for any $\bm{x}, \bm{y} \in \RR^{n^2}$ such that $\bm{x} = Q \bm{y}$, we have $\bm{x}^T \bm{y} \geq 0$. 
The first condition holds because $\det(D+Q)=\det(D^{1/2}(I+D^{-1/2}QD^{-1/2})D^{1/2}) = \det(D)\cdot\det(I+D^{-1/2} Q D^{-1/2}) > 0$.
The second condition is true because $\bm{x}^T \bm{y}=\bm{y}^T Q \bm{y} \geq 0$.
Hence, Eq.~\eqref{eq:sandberg} has a unique solution $\mathcal{W}$.

We now show that this solution satisfies the conditions of the stable point, that is, $\mathcal{W}=\mathcal{W}^T$, and there exists a skew-symmetric $P$ which satisfies Eq.~\eqref{eq:newgrad}.
Observe that
\begin{align*}
[S(\mathcal{W})^T]_{ij} &= S(\mathcal{W})_{ji} = (1/2)\cdot (f_{ij}+f_{ji})(\mathcal{W}_{ji}) = S(\mathcal{W}^T)_{ij},
\end{align*}
so $S(\mathcal{W})^T = S(\mathcal{W}^T)$.
Taking the transpose of Eq.~\eqref{eq:sandberg} and using $\Sigma=\Sigma^T$, $\Gamma=\Gamma^T$, and $S(\mathcal{W})^T = S(\mathcal{W}^T)$, we find
\begin{align*}
(M+M^T)/2 &= (\Sigma \mathcal{W}^T \Gamma + \Gamma \mathcal{W}^T \Sigma) + S(\mathcal{W}^T).
\end{align*}
But since there is only one solution to Eq.~\eqref{eq:sandberg}, we must have $\mathcal{W}=\mathcal{W}^T$.

Finally, we choose
\begin{align*}
P &= M - 2\Sigma \mathcal{W} \Gamma - R\\
\Rightarrow P+P^T &= (M+M^T) - 2(\Sigma\mathcal{W}\Gamma + \Gamma\mathcal{W}\Sigma) - 2S(W) = 0,
\end{align*}
using the fact that $\mathcal{W}=\mathcal{W}^T$ is a solution for Eq.~\eqref{eq:sandberg}.
Hence, this choice of $P$ is both skew-symmetric and satisfies Eq.~\eqref{eq:newgrad}.
\Halmos
\endproof

\subsection{Proof of Theorem \ref{thm:effect-changed-sigma-mu}}

\begin{repeattheorem}[Restatement of \ref{thm:effect-changed-sigma-mu}.]
Suppose $\Sigma_i=\Sigma$ for all firms, and let $M$ be the matrix of expected returns. Then, we have the following:
\begin{enumerate}
  \item {\bf Change in beliefs about expected returns:}
Let $\Sigma$ have the eigendecomposition $\Sigma = V \Lambda V^T$.
Then for $i, j, k, \ell \in [n]$, 
\begin{align}
%\label{eq:delWdelM_general}
\frac{\del W_{ij}}{\del M_{k\ell}} = \frac{1}{2\sqrt{\gamma_i\gamma_j\gamma_k\gamma_\ell}} \sum\limits_{s, t \in [n]} \frac{V_{is} V_{ks} V_{jt} V_{\ell t} + V_{is}V_{\ell s}V_{jt}V_{kt}}{\lambda_s + \lambda_{t}}.
\end{align}
  In particular, $W_{ij}$ is monotonically increasing with respect to $M_{ij}$.
  \item {\bf Risk scaling:} If the covariance $\Sigma$ changes to $c\Sigma$ ($c>0$), then $W$ changes to $(1/c)W$.
  \item {\bf Increase in perceived risk:} Suppose $\gamma_i = \gamma$ for all $i$, and the covariance $\Sigma$ increases to $\Sigma^\prime\succ \Sigma$.
Let $W$ and $W^\prime$ be the stable points under $\Sigma$ and $\Sigma^\prime$ respectively.
Then, $\tr(M^T (W^\prime - W)) < 0.$
\end{enumerate}
\end{repeattheorem}

\proof{Proof of Theorem \ref{thm:effect-changed-sigma-mu}.}
1. Let $(\lambda_i, \bmv_i)$ denote the $i^{th}$ eigenvalue and eigenvector of $\Gamma^{-1/2}\Sigma\Gamma^{-1/2}$, and let $V_{ij}=\bme_i^T \bmv_j$.
By Corollary~\ref{cor:stability:sharedSigma}, 
\begin{align*}
W &= \Gamma^{-1/2}\left( \sum_{s=1}^n \sum_{t=1}^n \frac{\bmv_s^T \Gamma^{-1/2}(M+M^T)\Gamma^{-1/2} \bmv_t}{2 (\lambda_r + \lambda_s)} \bmv_s \bmv_t^T \right) \Gamma^{-1/2}\\
\rarr \frac{\del W_{ij}}{\del M_{k\ell}}
&= \bme_i^T \Gamma^{-1/2}\left( \sum_{s=1}^n \sum_{t=1}^n \frac{\bmv_s^T \Gamma^{-1/2}(\bme_k \bme_\ell ^T+\bme_\ell \bme_k^T)\Gamma^{-1/2} \bmv_t}{2 (\lambda_s + \lambda_t)} \bmv_s \bmv_t^T \right) \Gamma^{-1/2}\bme_j\\
&= \frac{1}{2 \sqrt{\gamma_i \gamma_j \gamma_k \gamma_\ell}}
\left( \sum_{s=1}^n \sum_{t=1}^n \frac{\bmv_s^T (\bme_k \bme_\ell ^T+\bme_\ell \bme_k^T) \bmv_t}{(\lambda_s + \lambda_t)} (\bme_i^T\bmv_s) (\bmv_t^T \bme_j) \right)\\
&= \frac{1}{2\sqrt{\gamma_i\gamma_j\gamma_k\gamma_\ell}} \sum_{s=1}^n \sum_{t=1}^n 
\frac{V_{is} V_{ks} V_{jt} V_{\ell t} + V_{is}V_{\ell s}V_{jt}V_{kt}}{\lambda_s + \lambda_{t}}
\end{align*}
This proves Eq.~\eqref{eq:delWdelM_general}.
If $i = k, j = \ell$, we have:
\begin{align*}
\frac{\del W_{ij}}{\del M_{ij}}
&= (2\gamma_i\gamma_j)^{-1}\left( \sum_{s=1}^n \sum_{t=1}^n \underbrace{\frac{V_{is}^2 V_{jt}^2 + V_{is}V_{js}V_{jt}V_{it}}{\lambda_s + \lambda_{t}}}_{Z_{st}} \right)\\
&= (4\gamma_i\gamma_j)^{-1}\left( \sum_{s=1}^n \sum_{t=1}^n Z_{st} + \sum_{t=1}^n\sum_{s=1}^n Z_{ts}\right)\\
&= (4\gamma_i\gamma_j)^{-1} \sum_{s=1}^n \sum_{t=1}^n \left(Z_{st} + Z_{ts}\right)\\
&= (4\gamma_i\gamma_j)^{-1} \sum_{s=1}^n \sum_{t=1}^n \frac{\left(V_{is}V_{jt} + V_{js}V_{it}\right)^2}{\lambda_r + \lambda_s} > 0.
\end{align*}
Hence, $W_{ij}$ is monotonically increasing with respect to $M_{ij}$.

%If $i = k, j = \ell$ then it follows that $\frac{\del W_{ij}}{\del M_{ij}} = (2\gamma_i\gamma_j)^{-1} \sum_{r=1}^n \sum_{s=1}^n \frac{\left(V_{ir}V_{js} + V_{jr}V_{is}\right)^2}{2 (\lambda_r + \lambda_s)} > 0$.
%}

% &= (\gamma_i\gamma_j)^{-1}\left( \sum_{r=1}^n \sum_{s=1}^n \underbrace{\frac{V_{ir}V_{js} + V_{jr}V_{is}}{2 (\lambda_r + \lambda_s)} (V_{ir}V_{js})}_{Z_{rs}} \right)\\
% &= (2\gamma_i\gamma_j)^{-1}\left( \sum_{r=1}^n \sum_{s=1}^n Z_{rs} + \sum_{s=1}^n\sum_{r=1}^n Z_{sr}\right)\\
% &= (2\gamma_i\gamma_j)^{-1} \sum_{r=1}^n \sum_{s=1}^n \left(Z_{rs} + Z_{sr}\right)\\
% &= (2\gamma_i\gamma_j)^{-1} \sum_{r=1}^n \sum_{s=1}^n \frac{\left(V_{ir}V_{js} + V_{jr}V_{is}\right)^2}{2 (\lambda_r + \lambda_s)} > 0.

2. This follows from Corollary~\ref{cor:stability:sharedSigma}. 

3. By Corollary~\ref{cor:stability:sharedSigma}, $\vc(W) = \gamma^{-1} (\Sigma \otimes I + I \otimes \Sigma)^{-1} \vc(\frac{M + M^T}{2})$. Let $K = \gamma (\Sigma \otimes I + I \otimes \Sigma)$ and $K^\prime = \gamma (\Sigma^\prime \otimes I + I \otimes \Sigma^\prime)$. 
Since $\Sigma^\prime \succ \Sigma$ it follows that $K^\prime \succ K$. Therefore $K^{-1} \succ (K^\prime)^{-1}$. 

So, since $\vc(W^\prime - W) = ((K^\prime)^{-1} - K^{-1}) \vc(\frac{M + M^T}{2})$, we immediately obtain $\frac 1 2 \vc(M + M^T)^T \vc(W^\prime - W) < 0$. Since $W, W^\prime$ are symmetric it follows that $\vc(M^T)^T \vc(W^\prime - W) = \vc(M)^T \vc(W^\prime - W)$. So we have $\vc(M)^T \vc(W^\prime - W) < 0$. 

Since $\vc(M)^T \vc(W^\prime - W) = \tr(M^T (W^\prime - W))$, the conclusion follows. 
\Halmos
\endproof

\subsection{Hardness of Source Detection}\label{appendix:dwdm:approx}

% In Figure~\ref{fig:meanShock}, we observe that the smallest eigenvalue of $\hat{\Sigma}$ is much smaller than the second-smallest eigenvalue: $\lambda_n\ll \lambda_{n-1}$.

We begin by defining
\begin{align}
\label{eq:delWdelMapprox}
%\left|\frac{\del W_{ij}}{\del M_{k\ell}}\right| 
%&\approx 
\left|\frac{\del W_{ij}}{\del M_{k\ell}}\right|_{approx}
&\defeq \frac{\left|V_{in} V_{kn} V_{jn} V_{\ell n}\right|}{2 \lambda_n}.
\end{align}
This approximates the right hand side of Eq.~\eqref{eq:delWdelM_general} when the term corresponding to the smallest eigenvalue $\lambda_n$ dominates the sum.
We now show that if the corresponding eigenvector $\bmv_n$ is random, source detection becomes difficult. % impossible.
%We are ready to prove our main claim. 

\begin{proposition}[Hardness of Source Detection]
Suppose $\bmv_n$ is Haar-distributed, that is, $\bmv_n$ is distributed uniformly on the unit sphere $S^{n-1}$. Then, if $\Sigma = V \Lambda V^T$ and $\Gamma = I$, 
\begin{align*}
\PP\left[\max\limits_{i, j \in [n]: (i, j) \neq (k, \ell)} \left|\frac{\del W_{ij}}{M_{k \ell}}\right|_{approx} < 
\left|\frac{\del W_{k \ell}}{M_{k \ell}}\right|_{approx} \right] \leq O\left(\frac{1}{n}\right).
\end{align*}
\label{prop:mean-shock-single-column}
\end{proposition}

\proof{Proof of Proposition \ref{prop:mean-shock-single-column}.}
Without loss of generality we can set $k = 1, \ell = 2$ (the analysis of $k = \ell$ is identical). Notice that $\left|\frac{\del W_{ij}}{M_{k \ell}}\right|_{approx}$ 
is maximized at the $(i, j)$ that maximizes $\abs{V_{in} V_{jn}}$. 

Now, consider $(i, j) \in \{(1, 2), (3, 4), \dots, (n-1, n)\}$. Notice the distribution of $\bm{v}_n$ is permutation-invariant by assumption. Hence the joint distribution of $(V_{in}, V_{jn})$ is the same for all such pairs $(i, j)$. Hence the distribution of $\abs{V_{in} V_{jn}}$ is also the same for all such $(i, j)$. Therefore, 

\begin{align*}
\PP \left[
\arg\max\limits_{(i, j) \in \{(1, 2), (3, 4), \dots, (n-1, n)\}} \left|\frac{\del W_{ij}}{M_{12}}\right|_{approx} = (1, 2)
\right] \leq O(1/n)
\qquad\Halmos
\end{align*}
\endproof

\subsection{Proof of Proposition \ref{prop:firms-exchangeable-stronger}}
\label{appendix:exchangeable}

\begin{repeatproposition}[Restatement of Proposition \ref{prop:firms-exchangeable-stronger}]
Suppose $M, \Sigma, \Gamma$ exhibit community structure (Eq.~\eqref{eq:community}), and all the error terms $(\eps_i)_{i\in [n]}$ and $(\eps^\prime_{\theta_i, j})_{i,j\in [n]}$ are independent and identically distributed. 
%also assume that if $\theta_r = \theta_s$ then $\eps_{r} \sim \eps_{s}$, and $\eps_{\theta_i, r}^\prime \sim \eps_{\theta_i, s}^\prime$.
Let $\pi: [n] \to [n]$ be any intra-community permutation, and let $\Pi: \RR^n \to \RR^n$ be the corresponding column-permutation matrix: $\Pi(\bm{e}_i) = \bm{e}_{\pi(i)}$. 
Then, $W$ and $\Pi^T W \Pi$ are identically distributed.  % $W\sim \Pi^T W \Pi$.
\end{repeatproposition}

\proof{Proof.}
Let $H = \frac 1 2 (M + M^T)$. The fixed point equation for $W$ is given by Corollary \ref{cor:stability:sharedSigma} as $\Sigma W \Gamma + \Gamma W \Sigma = H$. Vectorization implies
$(\Gamma \otimes \Sigma + \Sigma \otimes \Gamma) \vc(W) = \vc(H)$. Let $X \sim Y$ denote that a pair of random variables $X, Y$ are identically distributed. We want to show $\Pi^T W \Pi \sim W$. Vectorization gives $\vc(\Pi^T W \Pi) = (\Pi^T \otimes \Pi^T) \vc(W)$. Let $P = (\Pi^T \otimes \Pi^T)$ and $K = (\Gamma \otimes \Sigma + \Sigma \otimes \Gamma)$ for shorthand. 

In this notation, we want to show that $P K^{-1} \vc(H) \sim K^{-1} \vc(H)$. Since $P$ is a permutation, we have $P K^{-1} \vc(H) = P K^{-1} P^T P \vc(H) = (PKP^T)^{-1} P \vc(H)$. Since the collections of random variables $\{\eps_i\}_i$ and $\{\eps_{\theta_i, j}^\prime \}_{i, j}$ are independent, we know $\vc(H)$ and $K$ are independent. So to show $(PKP^T)^{-1} P \vc(H) \sim K^{-1} \vc(H)$ it suffices to show that $P \vc(H) \sim \vc(H)$ and $PKP^T \sim K$. 

Notice $P \vc(H) = \vc(\Pi^T H \Pi)$. Hence, we want to show $\Pi^T H \Pi \sim H$, which holds iff $\Pi^T (M + M^T) \Pi \sim M + M^T$. Notice that $\Pi^T M^T \Pi = (\Pi^T M \Pi)^T$, so if $\Pi^T M \Pi \sim M$ then we obtain $\Pi^T M^T \Pi \sim M^T$ as well. It suffices to show $\Pi^T M \Pi \sim M$. 

Similarly, we can simplify $PKP^T = \Pi^T \Sigma \Pi \otimes \Pi^T \Gamma \Pi + \Pi^T \Gamma \Pi \otimes \Pi^T \Sigma \Pi$. It suffices to show $\Pi^T \Gamma \Pi \sim \Gamma$ and $\Pi^T \Sigma \Pi = \Sigma$. 

We are left to show that $\Pi^T \Sigma \Pi = \Sigma$ and $\Pi^T A \Pi \sim A$ for $A \in \{\Gamma, M\}$. 

{\em Proof of $\Pi^T \Sigma \Pi = \Sigma$.} Let $i, j \in [n]$. Then $(\Pi^T \Sigma \Pi)_{ij} = \Sigma_{\pi(i), \pi(j)} = g(\theta_{\pi(i)}, \theta_{\pi(j)})$. Since $\pi$ only commutes members within communities, $g(\theta_{\pi(i)}, \theta_{\pi(j)}) = g(\theta_i, \theta_j) = \Sigma_{ij}$. So $\Pi^T \Sigma \Pi = \Sigma$. 

{\em Proof of $\Pi^T \Gamma \Pi \sim \Gamma$.} Notice $\Pi^T \Gamma \Pi$ and $\Gamma$ are both diagonal. Let $i \in [n]$. Then $(\Pi^T \Gamma \Pi)_{ii} = \Gamma_{\pi(i), \pi(i)} = h(\theta_{\pi(i)}) + \eps_{\pi(i)} = h(\theta_{i}) + \eps_{\pi(i)}$. Since $\theta_{i} = \theta_{\pi(i)}$, we know $\eps_i\sim \eps_{\pi(i)}$. The conclusion follows. 

{\em Proof of $\Pi^T M \Pi \sim M$.} Let $i, j \in [n]$. Then $(\Pi^T M \Pi)_{ij} = M_{\pi(i), \pi(j)} = f(\theta_{\pi(i)}, \theta_{\pi(j)}) + \eps_{\theta_{\pi(i)}, \pi(j)}^\prime = f(\theta_{i}, \theta_{j}) + \eps_{\theta_{i}, \pi(j)}^\prime$. Since $\theta_{j} = \theta_{\pi(j)}$, we know that $\eps_{\theta_{i}, \pi(j)}^\prime\sim \eps_{\theta_{i}, j}^\prime$, and the conclusion follows. 
\Halmos
\endproof

% We now prove the corollary. 

% \begin{repeatcorollary}[Restatement of Corollary \ref{cor:firms-exchangeable}]
% Let $j_1, \dots, j_m \in [n]$ belong to the same community: $\theta_{j_1} = \dots = \theta_{j_m}$. 
% Suppose the conditions of Proposition~\ref{prop:firms-exchangeable-stronger} hold.
% Then, for any $i \in [n]$, the joint distribution of $(W_{i, j_1}, \ldots, W_{i, j_m})$ is exchangeable.
% \end{repeatcorollary}

% \proof{Proof.}
% Let $\tau: \{j_1, \dots, j_m\} \to \{j_1, \dots, j_m\}$ be an arbitrary permutation and $\pi: [n] \to [n]$ be its extension to $[n]$. Let $\Pi: \RR^n \to \RR^n$ be the corresponding permutation matrix. By Proposition \ref{prop:firms-exchangeable-stronger}, $\Pi^T W \Pi \sim W$. 

% \Halmos
% \endproof

% \begin{corollary}[Restatement of Proposition~\ref{prop:firms-exchangeable}.]
% Suppose that if $\theta_r = \theta_s$ then $\eps_{r} \sim \eps_{s}$ and $\eps_{\theta_i, r}^\prime \sim \eps_{\theta_i, s}^\prime$.
% Let $j_1, \dots, j_m \in [n]$ belong to the same community: $\theta_{j_1} = \dots = \theta_{j_m}$. 
% For any $i \in [n]$, the joint distribution of $(W_{i, j_1}, \ldots, W_{i, j_m})$ is exchangeable.
% \end{corollary}

% The proof follows from Proposition~\ref{prop:firms-exchangeable-stronger} below.

\subsection{Proof of Theorem \ref{thm:nonidentGammaM}}
\begin{repeattheorem}[Restatement of Theorem \ref{thm:nonidentGammaM}]
Consider two network settings $S=(\bmu_i, \Sigma, \gamma_i)_{i\in[n]}$ and $S^\prime=(\bmu_i, \Sigma, \gamma_i^\prime)_{i\in[n]}$ which differ only in the risk-aversions of firms $J=\{j\mid \gamma_j\neq \gamma^\prime_j\}\subseteq [n]$.
Then, there exists a setting $S^\dagger=(\bmu_i^\dagger, \Sigma, \gamma_i)_{i\in[n]}$ such that $\bmu_i=\bmu^\dagger_i$ for all $i\notin J$ and the stable networks under $S^\dagger$ and $S^\prime$ are identical.
\end{repeattheorem}

\proof{Proof of Theorem \ref{thm:nonidentGammaM}.}
First, consider the network settings $S$ and $S^\prime$. 
Let $\Gamma \in \RR^{n \times n}$ be a diagonal matrix with $\Gamma_{i, i} = \gamma_i$; define $\Gamma^\prime$ similarly under $S^\prime$.
Let the corresponding networks be $W$ and $W^\prime$, and let $\Delta_W=W^\prime-W$ and $\Delta_\Gamma=\Gamma^\prime-\Gamma$.
By Corollary~\ref{cor:stability:sharedSigma}, we have
\begin{align}
\Sigma W \Gamma + \Gamma W \Sigma &= \frac{M + M^T}{2} = \Sigma W^\prime \Gamma^\prime + \Gamma^\prime W^\prime \Sigma\nonumber\\
\rarr 
\frac{M + M^T}{2} &= \Sigma (W + \Delta_W)(\Gamma + \Delta_\Gamma) 
+ (\Gamma + \Delta_\Gamma) (W + \Delta_W) \Sigma \nonumber \\
&= \Sigma W \Gamma + \Gamma W \Sigma 
+ \Sigma \Delta_W \Gamma + \Gamma \Delta_W \Sigma 
+ \Sigma W \Delta_\Gamma + \Delta_\Gamma W \Sigma
+  \Sigma \Delta_W \Delta_\Gamma + \Delta_\Gamma \Delta_W \Sigma \nonumber \\
\rarr \Sigma \Delta_W \Gamma + \Gamma \Delta_W \Sigma 
&= - (\Sigma W \Delta_\Gamma + \Delta_\Gamma W \Sigma
+  \Sigma \Delta_W \Delta_\Gamma + \Delta_\Gamma \Delta_W \Sigma)\nonumber\\
&= - (\Sigma W^\prime \Delta_\Gamma + \Delta_\Gamma W^\prime \Sigma)\label{eq:dG}
\end{align}

%Let $\Gamma^\prime \in \RR^{n \times n}$ be the analogous matrix for $S^\prime$, and let $\Delta_\Gamma = \Gamma^\prime - \Gamma$. Let $W^\prime$ be the fixed point network under $S^\prime$, given by $\Sigma W^\prime \Gamma^\prime + \Gamma^\prime W^\prime \Sigma = \frac{M + M^T}{2}$. Let $\Delta_W = W^\prime - W$. Then, 
%
%
%Recall the fixed point equation $\Sigma W \Gamma + \Gamma W \Sigma = \frac{M + M^T}{2}$. 
%
%First, consider network setting $S$ versus $S^\prime$. Let $\Gamma \in \RR^{n \times n}$ be a diagonal matrix with $\Gamma_{i, i} = \gamma_i$. Let $\Gamma^\prime \in \RR^{n \times n}$ be the analogous matrix for $S^\prime$, and let $\Delta_\Gamma = \Gamma^\prime - \Gamma$. Let $W^\prime$ be the fixed point network under $S^\prime$, given by $\Sigma W^\prime \Gamma^\prime + \Gamma^\prime W^\prime \Sigma = \frac{M + M^T}{2}$. Let $\Delta_W = W^\prime - W$. Then, 
%
%\begin{align}
%\frac{M + M^T}{2} &= \Sigma (W + \Delta_W)(\Gamma + \Delta_\Gamma) 
%+ (\Gamma + \Delta_\Gamma) (W + \Delta_W) \Sigma \nonumber \\
%&= \Sigma W \Gamma + \Gamma W \Sigma 
%+ \Sigma \Delta_W \Gamma + \Gamma \Delta_W \Sigma 
%+ \Sigma W \Delta_\Gamma + \Delta_\Gamma W \Sigma
%+  \Sigma \Delta_W \Delta_\Gamma + \Delta_\Gamma \Delta_W \Sigma \nonumber \\
%\rarr \Sigma \Delta_W \Gamma + \Gamma \Delta_W \Sigma 
%&= - (\Sigma W \Delta_\Gamma + \Delta_\Gamma W \Sigma
%+  \Sigma \Delta_W \Delta_\Gamma + \Delta_\Gamma \Delta_W \Sigma) \label{eq:dG}
%\end{align}

Next, consider $S$ versus $S^\dagger$. Suppose that $M^\dagger$ has columns $\bm{\mu}_1^\dagger, \dots, \bm{\mu}_n^\dagger$ and let $\Delta_M = M^\dagger - M$. Let $W^\dagger$ be the fixed point network under $S^\dagger$, given by $\Sigma W^\dagger \Gamma + \Gamma W^\dagger \Sigma = \frac{M^\dagger + (M^\dagger)^T}{2}$. Let $\Delta_W^\dagger = W^\dagger - W$. Then a similar argument gives: 

\begin{align}
\frac{\Delta_M + \Delta_M^T}{2} = \Sigma \Delta_W^\dagger \Gamma + \Gamma \Delta_W^\dagger \Sigma \label{eq:dM}
\end{align}

Therefore, from Eq \eqref{eq:dG} and \eqref{eq:dM}, it follows that $W^\prime = W^\dagger$ if 
$$\frac{\Delta_M + \Delta_M^T}{2} = - (\Sigma W^\prime \Delta_\Gamma + \Delta_\Gamma W^\prime \Sigma).$$
%$\frac{\Delta_M + \Delta_M^T}{2} = - (\Sigma W \Delta_\Gamma + \Delta_\Gamma W \Sigma
%+  \Sigma \Delta_W \Delta_\Gamma + \Delta_\Gamma \Delta_W \Sigma) = - (\Sigma W^\prime \Delta_\Gamma) - (\Sigma W^\prime \Delta_\Gamma)^T$. 
Hence, $W^\prime=W^\dagger$ if we set $\Delta_M = - \Sigma W^\prime \Delta_\Gamma$.

It remains to show that $M^\dagger$ differs from $M$ only in columns corresponding to $J$. Suppose that $i \not \in J$. Then $\gamma_i = \gamma_i^\prime$, so $\Delta_\Gamma \bm{e}_i = \bm{0}$. We conclude that $\Delta_M \bm{e}_i = \bm{0}$ and hence $M \bm{e}_i = M^\dagger \bm{e}_i$. \Halmos
\endproof

\section{Experimental Details}\label{appendix-datasets}

\subsection{Fama-French Stock Market Data}

We use the Fama-French value-weighted asset returns dataset, for 96 assets over 625 months \citep{fama-french-2015}.

\subsection{OECD International Trade Data}

We use international trade statistics from the OECD to get quarterly measurements of bilateral trade between 46 large economies, including the top 15 world nations by GDP \cite{oecd-stats}. 
The data are available at the OECD Statistics webpage (\url{https://stats.oecd.org/}).
The data are measured quarterly from Q1 2010 to Q2 2022. 
We take the sum of trade flows $i \to j$ and $j  \to i$ to measure the weight of an edge $\{i, j\}$. 

To obtain the corresponding $\Sigma$, we run our inference procedure (Section~\ref{sec:model:inference}).
Since there is no data for within-country trade, the network has no self-loops ($W_{ii}=0$).
So we modify the inference according to Remark~\ref{rem:inference_prohibited} in Appendix~\ref{appendix:inference}.

%``International Trade and Balance of Payments -- International Merchandise Trade Statistics (Monthly and Quarterly) -- Quarterly (by partner country) -- Quarterly Intl Trade Statistics.''

%The full list of countries is: Argentina,
%Australia,
%Austria,
%Belgium,
%Brazil,
%Canada,
%Chile,
%China (People's Republic of),
%Colombia,
%Costa Rica,
%Czech Republic,
%Denmark,
%Estonia,
%Finland,
%France,
%Germany,
%Greece,
%Hungary,
%Iceland,
%India,
%Indonesia,
%Ireland,
%Israel,
%Italy,
%Japan,
%South Korea,
%Latvia,
%Lithuania,
%Luxembourg,
%Mexico,
%Netherlands,
%New Zealand,
%Norway,
%Poland,
%Portugal,
%Russia,
%Saudi Arabia,
%Slovak Republic,
%Slovenia,
%South Africa,
%Spain,
%Sweden,
%Switzerland,
%Türkiye,
%United Kingdom,
%United States. 

\subsection{Outlier Detection Simulation}\label{appendix-gamma-exps}

The experiments in Figure \ref{fig:deviator-gamma-detection} proceed as follows. Fix a number of communities $k$ and number of firms $n$. Fix a value of $\sigma > 0$. For us, $k = 2$, $n \in \{20, 100, 300\}$, and $\sigma \in \{\sigma_1, \dots, \sigma_{10}\}$, where the $\sigma_i$ are logarithmically spaced on the interval $[0.1, 1]$, so that 
\begin{align*}
\sigma \in \{0.1       , 0.12915497, 0.16681005, 0.21544347, 0.27825594, \\
0.35938137, 0.46415888, 0.59948425, 0.77426368, 1.0\}
\end{align*}

For a setting of $n, k, \sigma$, we perform the following simulation $m = 500$ times. 

{\em Generate communities.} Generate the community membership matrix $\Theta \in \bb^{n \times k}$ with rows independently and uniformly at random from $\{\bm{e_1}, \dots, \bm{e_k}\}$. 

{\em Generate the network setting.} The deterministic functions $f, g, h$ for $M, \Sigma, \Gamma$ respectively are as follows. First $f(\theta_1, \theta_2) = f(\theta_2, \theta_1) = 1$ and $f = 0$ otherwise. Next, let $G \in \RR^{k \times k}$ be the matrix $G_{ij} = g(\theta_i, \theta_j)$. Then $G$ is generated from a normalized Wishart distribution centered at $I_k$ and with $5$ degrees of freedom. Finally, $h(\theta_i) = 1$ for all $i$. 

The noise variables for agent beliefs are as follows. Sample i.i.d. $\eps_i$ according to a $N(0, \sigma^2)$ distribution truncated to $[-0.5, 0.5]$ for all $i$. Sample $\eps_{\theta_{i}, j}^\prime \iid N(0, \sigma^2)$ for all $i, j$. 

{\em Designate an outlier.} Set the the noise parameter $\eps_1 = -0.5$ for firm $1$ (the risk-seeker), so as $\sigma \to 0$, $\gamma_1$ gets further separated from all other $\gamma_i$. 

{\em Outlier detection simulation.} Then for a random firm $i$ such that $\theta_i \neq \theta_1$, we test whether the outlier $\hat j \defeq \arg\max\limits_{j: \theta_j = \theta_1} \abs{W_{i, j}}$ is equal to the true outlier firm $1$. 

{\em Collate results.} Once the $m = 500$ runs are completed for a single setting of $n, k, \sigma$, we obtain an estimate $\hat p$ for the probability of successful deviator detection at this setting of parameters. We plot a confidence interval $[p - 2 \sqrt{\frac{\hat p (1 - \hat p)}{m}}, p + 2 \sqrt{\frac{\hat p (1 - \hat p)}{m}}]$. This is plotted on the $y$-axis. The $x$-axis quantifies how much $\gamma_1$ deviates from the mean, in terms of the number of standard deviations of the truncated normal distribution $\eps_i$.

% !TEX root = ./paper-draft.tex

% Appendix A.6 (price update rule)
% Appendix A.8 (proof of thrm 5, the one about random covariances)
% Appendix A.9 (proof that SDP for Sigma gives the MLE)
% Appendix A.11 (stable point after a shock)
% Appendix A.12 (gradient of W wrt M)
% Appendix A.13 (hardness of source detection, the Haar argument)
% A.14 (equidistribution of W after intra-community permutation)
% A.15 (indistinguishability between gamma shift and M shift under community structure)
% Appendix B (experimental details) 

\end{APPENDICES}

\end{document}